\documentclass[a4paper,reqno]{amsart}
\usepackage{cmap} 
\usepackage[british]{babel}
\usepackage[T1]{fontenc}
\usepackage{amssymb,amsthm,bm}
\usepackage{mathtools}
  \mathtoolsset{mathic}
  \allowdisplaybreaks
\usepackage{textcomp}
\usepackage{courier,mathptmx}
\usepackage[spacing]{microtype}
\usepackage[numbers,sort&compress]{natbib} \renewcommand{\citeyearpar}{\cite}
\usepackage{paralist}
\usepackage{hyperref,url}
  \hypersetup{  
    bookmarksnumbered
  }
  \hypersetup{  
    breaklinks=false,
    pdfborderstyle={/S/U/W 0},  
    citebordercolor=.235 .702 .443,
    urlbordercolor=.255 .412 .882,
    linkbordercolor=.804 .149 .19,
  }
  \hypersetup{ 
    pdfauthor={Gert de Cooman \& Erik Quaeghebeur},
    pdftitle={Exchangeability and sets of desirable gambles},
    pdfkeywords={sets of desirable gambles, coherence, desirability, representation}
  }
\usepackage[all]{hypcap} 
\usepackage{fixltx2e} 
\usepackage{xcolor}
\usepackage{tikz}
  \usetikzlibrary{calc,matrix,arrows}
\usepackage{wrapfig}
  \newlength{\columnsepsave}
\usepackage{pdfsync}
\DeclarePairedDelimiter{\group}{(}{)} 
\pgfdeclarelayer{background}
\pgfdeclarelayer{foreground}
\pgfsetlayers{background,main,foreground}
\newcommand{\xtick}{node[fill,inner sep=0pt,minimum width=.5pt,minimum height=4pt] {}}
\newcommand{\ytick}{node[fill,inner sep=0pt,minimum width=4pt,minimum height=.5pt] {}}
\newcommand{\values}{\mathcal{X}}
\newcommand{\cnts}{\mathcal{N}}
\newcommand{\pspace}{\varOmega}
\newcommand{\rv}{X}
\newcommand{\event}{C}

\newcommand{\nats}{\mathbb{N}}
\newcommand{\reals}{\mathbb{R}}

\newcommand{\gambles}{\mathcal{G}}
\newcommand{\poly}{\mathcal{V}}
\newcommand{\pospoly}{\mathcal{V}^+_0}
\newcommand{\npospoly}{\mathcal{V}^-}
\newcommand{\nnegpoly}{\mathcal{V}^+}
\newcommand{\desirs}{\mathcal{R}}
\newcommand{\cdesirs}{\mathcal{S}}
\newcommand{\tdesirs}{\mathcal{H}}
\newcommand{\atdesirs}{\mathcal{F}}
\newcommand{\natex}[1][]{\mathcal{E}_{#1}\group}
\newcommand{\exnatex}[2][N]{\mathcal{E}_{\mathrm{ex}}^{#1}({#2})}
\newcommand{\bernnatex}{\mathcal{E}_{\mathrm{Be}}\group}
\newcommand{\weakly}[1]{\mathcal{D}_{#1}}
\newcommand{\marginally}[1]{\mathcal{M}_{#1}}
\newcommand{\identity}[1]{\iden_{#1}}
\newcommand{\assessment}{{\mathcal{A}}}
\newcommand{\passessment}{{\mathcal{A}}}
\newcommand{\cone}{{\mathcal{C}}}
\newcommand{\subspace}{{\mathcal{K}}}
\newcommand{\spacecone}{{(\subspace,\cone)}}
\newcommand{\alldesirs}[1][]{\mathbb{D}_{#1}}
\newcommand{\maxdesirs}[1][]{\mathbb{M}_{#1}}
\newcommand{\allexdesirs}{\mathbb{D}_{\mathrm{ex}}}
\newcommand{\allberndesirs}{\mathbb{D}_{\mathrm{Be}}}
\newcommand{\chain}{\mathbb{K}}
\newcommand{\permuts}{\mathcal{P}}
\newcommand{\average}{\mathcal{U}}
\newcommand{\set}[2]{\left\{#1\colon#2\right\}}
\newcommand{\biggset}[2]{\biggl\{#1\colon#2\biggr\}}
\newcommand{\abs}[1]{\lvert#1\rvert}
\newcommand{\pr}{P}
\newcommand{\lpr}{{\underline{\pr}}}
\newcommand{\upr}{{\overline{\pr}}}
\newcommand{\cpr}{Q}
\newcommand{\clpr}{\underline{\cpr}}

\newcommand{\sample}[1]{#1}
\newcommand{\osample}[1]{\obs{\sample{#1}}}
\newcommand{\rsample}[1]{\rest{\sample{#1}}}
\newcommand{\cnt}[1]{#1}
\newcommand{\ocnt}[1]{\obs{\cnt{#1}}}
\newcommand{\rcnt}[1]{\rest{\sample{#1}}}
\newcommand{\bcntf}{\cntf}
\newcommand{\btheta}{\theta}
\newcommand{\cntf}{T}
\newcommand{\tuple}[2]{{({#1}_1,\dots,{#1}_{#2})}}
\newcommand{\atom}[1]{{[#1]}}
\newcommand{\batom}{\atom}
\newcommand{\oatom}[1]{\atom{\obs{#1}}}
\newcommand{\ratom}[1]{\atom{\rest{#1}}}
\newcommand{\simplex}[1]{\varSigma_{#1}}
\newcommand{\bern}[1]{B_{\cnt{#1}}}
\newcommand{\bexp}[2]{b_{#1}^{#2}}
\newcommand{\update}[2]{{{#1}\vert{#2}}}
\newcommand{\walleyupdate}[2]{{{#1}\Vert{#2}}}
\newcommand{\restrict}[2]{{{#1}\rfloor{#2}}}
\newcommand{\obs}[1]{\check{#1}}
\newcommand{\rest}[1]{\hat{#1}}
\newcommand{\then}{\Rightarrow}
\newcommand{\asa}{\Leftrightarrow}
\newcommand{\vacuous}[1][N]{\desirs^{#1}_{\mathrm{ex,v}}}

\newcommand{\nonpositive}[1][\subspace]{{#1}_{\preceq0}}
\newcommand{\positive}[1][\subspace]{{#1}_{\succ0}}

\DeclareMathOperator{\ex}{ex}
\DeclareMathOperator{\muhy}{Hy}
\DeclareMathOperator{\enl}{enl}
\DeclareMathOperator{\mult}{Mn}
\DeclareMathOperator{\cmult}{CoMn}
\DeclareMathOperator{\ocntf}{Co}
\DeclareMathOperator{\iden}{id}
\DeclareMathOperator{\exten}{ext}
\DeclareMathOperator{\proj}{proj}
\DeclareMathOperator{\opspanning}{span}
\DeclareMathOperator{\posi}{posi}

\newcommand{\imuhy}[2]{\muhy_{{#1}}^{{#2}}}
\newcommand{\vmuhy}[1]{\imuhy{}{#1}}
\newcommand{\iex}[2]{\ex_{{#1}}^{{#2}}}
\newcommand{\vex}[1]{\iex{}{#1}}
\newcommand{\vmult}[1]{\mult^{#1}}
\newcommand{\vcmult}[1]{\cmult^{#1}}
\newcommand{\vbcntf}[1]{\bcntf^{#1}}
\newcommand{\vocntf}[1]{\ocntf^{#1}}
\newcommand{\vcnts}[1]{\cnts^{#1}}
\newcommand{\vproj}[2]{\proj_{#1}^{#2}}
\newcommand{\vexten}[2]{\exten_{#1}^{#2}}
\newcommand{\venl}[2]{\enl_{#1}^{#2}}
\newcommand{\vassessment}[1]{\assessment^{#1}}
\newcommand{\vdesirs}[1]{\desirs^{#1}}
\newcommand{\vcdesirs}[1]{\cdesirs^{#1}}
\newcommand{\vtdesirs}[1][]{\tdesirs^{#1}}
\newcommand{\vatdesirs}[1][]{\atdesirs^{#1}}
\newcommand{\fdesirs}[2]{\desirs^{#1}_{#2}}

\newcommand{\ftdesirs}[2][]{\tdesirs^{#1}_{#2}}
\newcommand{\vpoly}[1][\!]{\poly^{#1}(\vsimplex)}
\newcommand{\vpospoly}{\pospoly(\vsimplex)}
\newcommand{\vnpospoly}[1][]{\npospoly_{#1}(\vsimplex)}
\newcommand{\vnnegpoly}{\nnegpoly(\vsimplex)}
\newcommand{\irestrict}[4]{#1_{{#2}}^{{#3}}\rfloor{{#4}}}
\newcommand{\vrestrict}[2]{\irestrict{\desirs}{}{\rest{#1}}{\ocnt{#2}}}
\newcommand{\vcrestrict}[2]{\irestrict{\cdesirs}{}{\rest{#1}}{\ocnt{#2}}}
\newcommand{\vtrestrict}[2][]{\irestrict{\tdesirs}{}{#1}{\ocnt{#2}}}

\newcommand{\vcntsrest}[1]{\cnts^{\rest{#1}}}
\newcommand{\vextenrest}[2]{\exten_{\rest{#1}}^{\rest{#2}}}

\newcommand{\vcmultrest}[1]{\cmult^{\rest{#1}}}
\newcommand{\vsimplex}{\simplex{\values}}

\newtheorem{theorem}{Theorem}
\newtheorem{proposition}[theorem]{Proposition}
\newtheorem{lemma}[theorem]{Lemma}
\newtheorem{corollary}[theorem]{Corollary}
\newtheorem{definition}{Definition}
\theoremstyle{definition}
\newtheorem{example}{Example}

\begin{document}
\title{Exchangeability and sets of desirable gambles}
\author{Gert de Cooman}
\address{Ghent University\\SYSTeMS Research Group\\Technologiepark--Zwijnaarde 914\\9052 Zwijnaarde\\Belgium}
\email{Gert.deCooman@UGent.be}
\author{Erik Quaeghebeur}
\address{Ghent University\\SYSTeMS Research Group\\Technologiepark--Zwijnaarde 914\\9052 Zwijnaarde\\Belgium}
\email{Erik.Quaeghebeur@UGent.be}

\begin{abstract}
  Sets of desirable gambles constitute a quite general type of uncertainty model with an interesting geometrical interpretation.
  We give a general discussion of such models and their rationality criteria.
  We study exchangeability assessments for them, and prove counterparts of de Finetti's finite and infinite representation theorems.
  We show that the finite representation in terms of count vectors has
  a very nice geometrical interpretation, and that the representation in terms of frequency vectors is tied up with multivariate Bernstein (basis) polynomials.
  We also lay bare the relationships between the representations of updated exchangeable models, and discuss conservative inference (natural extension) under exchangeability and the extension of exchangeable sequences.
\end{abstract}


\maketitle

\section{Introduction}
In this paper, we bring together desirability, an interesting approach to modelling uncertainty, with exchangeability, a structural assessment for uncertainty models that is important for inference purposes.
\par
Desirability, or the theory of (coherent) sets of desirable gambles, has been introduced with all main ideas present---so far as our search has unearthed---by \citet{williams1976,williams1975,williams1975lmps}.
Building on \citeauthor{finetti19745}'s betting framework \citeyearpar{finetti19745}, he considered the `acceptability' of \emph{one-sided} bets instead of \emph{two-sided} bets.
This relaxation leads one to work with cones of bets instead of with linear subspaces of them.
The germ of the theory was, however, already present in \citeauthor{smith1961}'s work \citeyearpar[p.~15]{smith1961}, who used a (generally) open cone of `exchange vectors' when talking about currency exchange.
Both authors influenced \citeauthor{walley1991} \citeyearpar[Section~3.7 and App.~F]{walley1991}, who describes three variants (almost, really, and strictly desirable gambles) and emphasises the conceptual ease with which updated (or posterior) models can be obtained in this framework \citeyearpar{walley2000}.
\citeauthor{moral2003} \citeyearpar{moral2003,moral2005b} then took the next step and applied the theory to study epistemic irrelevance, a structural assessment.
He also pointed out how conceptually easy extension, marginalisation, and conditioning are in this framework.
\Citet{cooman2005c} made a general study of transformational symmetry assessments for desirable gambles.
Recently, \citet{couso2009isipta} discussed the relationship with credal sets, computer representation, and maximal sets of desirable gambles.
\par
The structural assessment we are interested in here, is exchangeability.
Conceptually, it says that the order of the samples in a sequence of them is irrelevant for inference purposes.
The first detailed study of this concept was made by \citet{finetti1937}, using the terminology of `equivalent' events.
He proved the now famous Representation Theorem, which is often interpreted as stating that a sequence of random variables is exchangeable if it is conditionally independent and identically distributed.
Other important work---all using probabilities or previsions---was done by, amongst many others, \citet{hewitt1955}, \citet{heath1976}, and \citet{diaconis1980}.\footnote{See, e.g., \citet{kallenberg2002,kallenberg2005} for a measure-theoretic discussion of exchangeability.}
Exchangeability in the context of imprecise-probability theory---using lower previsions---was studied by \citeauthor{walley1991} \citeyearpar[Section~9.5]{walley1991} and more in-depth by \Citeauthor{cooman2006d} \cite{cooman2005c,cooman2007e,cooman2006d}.
The first embryonic study of exchangeability using desirability was recently performed by \citeauthor{Quaeghebeur2009phd} \citeyearpar[Section~3.1.1]{Quaeghebeur2009phd}.
\par
Here we present the results of a more matured study of exchangeability using sets of desirable gambles.
First, in Section~\ref{sec:gendesir}, we give a general discussion of desirability, coherence---the criteria that define which sets of desirable gambles are rational uncertainty models---and the smallest and maximal such sets compatible with some assessment.
Next, in Section~\ref{sec:desir}, we introduce the special case of this general theory that will form the basis of our analysis of exchangeability using the theory of desirable gambles.
Then, in Section~\ref{sec:finite-xch}, we give a desirability-based analysis of finite exchangeable sequences, presenting a Representation Theorem---both in terms of count and frequency vectors---and treating the issues of natural extension and updating under exchangeability.
Building on these results, we extend our scope to countable exchangeable sequences in Section~\ref{sec:infinite-xch}, where we present a second Representation Theorem---in terms of frequency vectors---and again also treat updating and natural extension.
Finally, in Section~\ref{sec:seq-ext-finite}, we see if and how finite exchangeable sequences can be extended to longer, finite or even infinite exchangeable sequences.
\par
Proofs of this paper's results are included in Appendix~\ref{app:proofs}.
Appendix~\ref{app:bernstein} collects a few relevant facts about multivariate Bernstein basis polynomials.

\section{A general discussion of desirability and coherence}\label{sec:gendesir}
Consider a non-empty set $\pspace$ describing the possible and mutually exclusive outcomes of some experiment.
We also consider a subject who is uncertain about the outcome of the experiment.

\subsection{Sets of desirable gambles}\label{sec:sets-of-desirs}
A \emph{gamble} $f$ is a bounded real-valued map on $\pspace$, and it is interpreted as an uncertain reward.
When the actual outcome of the experiment is $\omega$, then the corresponding (possibly negative) reward is $f(\omega)$, expressed in units of some pre-determined linear utility.
This is illustrated for $\pspace=\{\omega,\varpi\}$.
$\gambles(\pspace)$ denotes the set of all gambles on $\pspace$, $\gambles^+_0(\pspace)$ the non-negative non-zero ones, and $\gambles^-(\pspace)$ the non-positive ones.
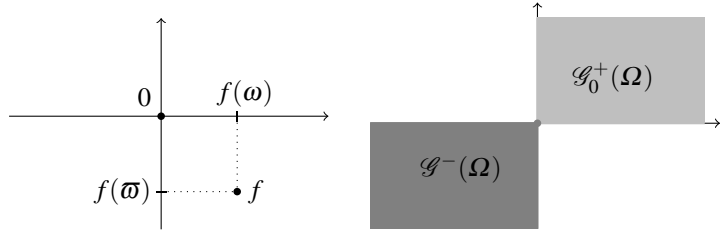
\begin{center}
  \begin{tikzpicture}[scale=2]
    \draw[->] (-1,0) -- (1.1,0);
    \draw[->] (0,-.75) -- (0,.65);
    \node[fill,circle,inner sep=1pt,label={[xshift=-1pt]right:$f$}] (f) at (.5,-.5) {};
    \draw[dotted] (f) -- (f |- 1,0) coordinate (fomega) \xtick node[above,xshift=3pt] {$f(\omega)$};
    \draw[dotted] (f) -- (f -| 0,1) coordinate (fomegaprime) \ytick  node[left] {$f(\varpi)$};
    \node[fill,circle,inner sep=1pt,label={[shift={(1pt,-1pt)}]above left:$0$}] (zero) at (0,0) {};
  \end{tikzpicture}
  \hspace{1em}
  \begin{tikzpicture}[scale=2]
    \draw[->] (-1.1,0) -- (1.2,0);
    \draw[->] (0,-.7) -- (0,.8);
    \fill[gray] (0,0) rectangle (-1.1,-.7);
    \draw[gray,thick] (0,-.7) |- (-1.1,0);
    \fill[lightgray] (0,0) rectangle (1.1,.7);
    \draw[lightgray,thick] (1.1,0) -| (0,.7);
    \node[fill=gray,circle,inner sep=1pt] (zero) at (0,0) {};
    \node at (.5,.3) {$\gambles^+_0(\pspace)$};
    \node at (-.5,-.3) {$\gambles^-(\pspace)$};
  \end{tikzpicture}
\end{center}
We also use the following notational conventions throughout: subscripting a set with zero corresponds to \emph{removing} zero (or the zero gamble) from the set, if present.
For example $\reals^+$ ($\reals^+_0$) is the set of non-negative (positive) real numbers including (excluding) zero.
Furthermore, ${f\geq g}$ iff ${f(\omega)\geq g(\omega)}$ for all~$\omega$ in~$\pspace$; $f>g$ iff $f\geq g$ and $f\neq g$.
\par
We say that a non-zero gamble $f$ is \emph{desirable} to a subject if he accepts to engage in the following transaction, where:
\begin{inparaenum}[(i)]
  \item the actual outcome $\omega$ of the experiment is determined, and
  \item he receives the reward $f(\omega)$, i.e., his capital is changed by $f(\omega)$.
\end{inparaenum}
The zero gamble is not considered to be desirable.\footnote{The nomenclature in the literature regarding desirability is somewhat confusing, and we have tried to resolve some of the ambiguity here. Our notion of desirability coincides with \citeauthor{walley2000}'s later \citeyearpar{walley2000} notion of desirability, initially (and quite recently \cite{couso2009isipta}) also used by \citet{moral2003}. \citeauthor{walley1991} in his book \citep[App.~F]{walley1991} and \citeauthor{moral2005b} in a later paper \citeyearpar{moral2003} use another notion of desirability. The difference between the two approaches resides in whether the zero gamble is assumed to be desirable or not. We prefer to use the non-zero version here, because it is better behaved in conjunction with our notion of weak desirability in Definition~\ref{def:weak-desirability}.}
\par
We try and model the subject's beliefs about the outcome of the experiment by considering which gambles are desirable for him.
We suppose the subject has some set $\desirs\subseteq\gambles(\pspace)$ of desirable gambles.

\subsection{Coherence}
Not every such set should be considered as a reasonable model, and in what follows, we give an abstract and fairly general treatment of ways to impose `rationality' constraints on sets of desirable gambles.
\par
We begin with a few preliminary definitions involving ordered linear spaces.
\par
The set $\gambles(\pspace)$ of all gambles on $\pspace$ is a linear space with respect to the (point-wise) addition of gambles, and the (point-wise) scalar multiplication of gambles with real numbers.
The \emph{positive hull operator} $\posi$ generates the set of \emph{strictly} positive linear combinations of elements of its argument set:  for any subset $\assessment$ of $\gambles(\pspace)$,
\begin{equation}
  \posi(\assessment)
  \coloneqq\set{\sum_{k=1}^n\lambda_kf_k}
  {f_k\in\assessment,\lambda_k\in\reals^+_0,n\in\nats_0}.
\end{equation}
A subset $\cone$ of $\gambles(\pspace)$ is a \emph{convex cone} if it is closed under (strictly) positive linear combinations, or in other words, if $\posi(\cone)=\cone$.
\par
Consider a linear subspace $\subspace$ of the linear space $\gambles(\pspace)$.
With any convex cone $\cone\subset\subspace$ such that $0\in\cone$ we can always associate a \emph{vector ordering} $\succeq$ on $\subspace$, defined as follows:\footnote{We require that $\cone$ should be \emph{strictly included} in $\subspace$ ($\cone\neq\subspace$) because otherwise the ordering $\succeq$ would be trivial: we would have that $f\succeq g$ for all $f,g\in\subspace$.}
\begin{equation}
  f\succeq g\asa f-g\in\cone\asa f-g\succeq0.
\end{equation}
The partial ordering $\succeq$ turns $\subspace$ into an ordered linear space \cite[Section~11.44]{schechter1997}.
We also write $f\succ g$ if $f-g\in\cone_0$, or in other words, if $f\succeq g$ and $f\neq g$.
As usual, $f\preceq g$ means $g\succeq f$ and similarly, $f\prec g$ means $g\succ f$.
Finally, we let
\begin{equation}
  \nonpositive
  \coloneqq\set{f\in\subspace}{f\preceq0}
  =-\cone
  \quad\text{ and }\quad
  \positive
  \coloneqq\set{f\in\subspace}{f\succ0}
  =\cone_0.
\end{equation}

\begin{definition}[Avoiding non-positivity and coherence]
  Let\/ $\subspace$ be a linear subspace of\/ $\gambles(\pspace)$ and let $\cone\subset\subspace$ be a convex cone containing the zero gamble $0$.
  We say that a set of desirable gambles~$\desirs\subseteq\subspace$ \emph{avoids non-positivity relative to $\spacecone$} if $f\not\preceq0$ for all gambles~$f$ in~$\posi(\desirs)$, or in other words if $\nonpositive\cap\posi(\desirs)=\emptyset$.
 \par
  We say that a set of desirable gambles~$\desirs\subseteq\subspace$  is \emph{coherent relative to $\spacecone$} if it satisfies the following requirements, for all gambles $f$, $f_1$, and $f_2$ in $\subspace$ and all real $\lambda>0$:
  \begin{compactenum}[\upshape {D}1.]
    \item if $f=0$ then $f\notin\desirs$;\label{item:Dnonzero}
    \item if $f\succ0$ then $f\in\desirs$, or equivalently $\positive\subseteq\desirs$;\label{item:Dapg}
    \item if $f\in\desirs$ then $\lambda f\in\desirs$ [scaling];\label{item:Dscl}
    \item if $f_1,f_2\in\desirs$ then $f_1+f_2\in\desirs$ [combination].\label{item:Dcmb}
  \end{compactenum}
  We denote by $\alldesirs[\spacecone](\pspace)$ the set of sets of desirable gambles that are coherent relative to $\spacecone$.
\end{definition}
\noindent
Requirements~D\ref{item:Dscl} and~D\ref{item:Dcmb} make~$\desirs$ a \emph{cone}: $\posi(\desirs)=\desirs$.
Due to~D\ref{item:Dapg}, it includes $\positive$; due to~D\ref{item:Dnonzero}, D\ref{item:Dapg} and~D\ref{item:Dcmb}, it excludes $\nonpositive$:
\begin{compactenum}[\upshape {D}1.]
\addtocounter{enumi}{4} 
  \item if $f\preceq0$ then $f\notin\desirs$, or equivalently $\nonpositive\cap\desirs=\emptyset$ .\label{item:Dnonzeroapl}
\end{compactenum}
The non-triviality requirement $\cone\neq\subspace$ makes sure that $\subspace$ is \emph{never} coherent relative to $\spacecone$.
On the other hand, $\positive$ is \emph{always} coherent relative to $\spacecone$, and it is the smallest such subset of $\subspace$.

\subsection{Natural extension}
If we consider an arbitrary non-empty family of sets of desirable gambles $\desirs_i$, $i\in I$ that are coherent relative to $\spacecone$, then their intersection $\bigcap_{i\in I}\desirs_i$ is still coherent relative to $\spacecone$.
This is the idea behind the following result.
If a subject gives us an \emph{assessment}, a set $\assessment\subseteq\subspace$ of gambles on $\pspace$ that he finds desirable, then it tells us exactly when this assessment can be extended to a coherent set, and how to construct the smallest such set.

\begin{theorem}[Natural extension]\label{theo:natex}
  Let\/ $\subspace$ be a linear subspace of\/ $\gambles(\pspace)$ and let $\cone\subset\subspace$ be a convex cone containing the zero gamble $0$.
  Consider an assessment $\assessment\subseteq\subspace$, and define its $\spacecone$-\emph{natural extension}:\footnote{As usual, in this expression, we let $\bigcap\emptyset=\subspace$.}
  \begin{equation}
    \natex[\spacecone]{\assessment}
    \coloneqq\bigcap\set{\desirs\in\alldesirs[\spacecone](\pspace)}{\assessment\subseteq\desirs} 
  \end{equation}
  Then the following statements are equivalent:
  \begin{compactenum}[\upshape (i)]
  \item\label{item:Aapl} $\assessment$ avoids non-positivity relative to $\spacecone$;
  \item\label{item:Asubsetcoh} $\assessment$ is included in some set of desirable gambles that is coherent relative to~$\spacecone$;
  \item\label{item:natextotAnonempty} $\natex[\spacecone]{\assessment}\neq\subspace$;
  \item\label{item:natextAiscoh} the set of desirable gambles $\natex[\spacecone]{\assessment}$ is coherent relative to~$\spacecone$;
  \item\label{item:natextAsmallestcoh} $\natex[\spacecone]{\assessment}$ is the smallest set of desirable gambles that is coherent relative to~$\spacecone$ and includes $\assessment$.
  \end{compactenum}
  When any (and hence all) of these equivalent statements hold, then
  \begin{equation}
    \natex[\spacecone]{\assessment}
    = \posi\bigl(\positive\cup\assessment\bigr).\label{eq:posi-natex}
  \end{equation}
\end{theorem}

\noindent This shows that if we have an assessment $\assessment$ with a finite description, we can represent its natural extension on a computer by storing a finite description of its extreme rays.

\subsection{Maximal coherent sets}
We see that the set $\smash[b]{\alldesirs[\spacecone]}(\pspace)$ of all sets of desirable gambles that are coherent relative to~$\spacecone$ can be partially ordered by set inclusion~$\subseteq$.\footnote{This structure is a complete meet-semilattice, where intersection plays the role of infimum.}
Following in the footsteps of \citet{couso2009isipta}, let us now investigate the maximal elements of this poset in more detail.
\par
An element $\desirs$ of $\alldesirs[\spacecone](\pspace)$ is called \emph{maximal} if it is not strictly included in any other element of  $\alldesirs[\spacecone](\pspace)$, or in other words, if adding any gamble $f$ to $\desirs$ makes sure we can no longer extend the result $\desirs\cup\{f\}$ to a set that is still coherent relative to  $\spacecone$:
\begin{equation}
  \maxdesirs[\spacecone](\pspace)
  =\set{\desirs\in\alldesirs[\spacecone](\pspace)}
  {(\forall\desirs'\in\alldesirs[\spacecone](\pspace))
  (\desirs\subseteq\desirs'\then\desirs=\desirs')}
\end{equation}
is the set of all maximal elements of $\alldesirs[\spacecone](\pspace)$.
\par
The following proposition provides a characterisation of such maximal elements.

\begin{proposition}\label{prop:maximality}
  Let\/ $\subspace$ be a linear subspace of\/ $\gambles(\pspace)$ and let $\cone\subset\subspace$ be a convex cone containing the zero gamble $0$.
  Let $\desirs\in\alldesirs[\spacecone](\pspace)$, then $\desirs$ is a maximal coherent set relative to $\spacecone$ iff
  \begin{equation}\label{eq:maximality}
    (\forall f\in\subspace_0)
    (f\notin\desirs\then-f\in\desirs).
  \end{equation}
\end{proposition}

For the following important result(s), a constructive proof can be given in case $\pspace$ is finite, based on the same ideas as in \cite{couso2009isipta}.
They guarantee that $\alldesirs[\spacecone](\pspace)$ has all the useful properties of a \emph{strong belief structure} \cite{cooman2003a}.
In Appendix~\ref{app:proofs}, we give a non-constructive proof (based on Zorn's Lemma) for $\pspace$ that may also be infinite.

\begin{theorem}\label{theo:maximality}
  Let\/ $\subspace$ be a linear subspace of\/ $\gambles(\pspace)$ and let $\cone\subset\subspace$ be a convex cone containing the zero gamble $0$.
  Consider any subset $\assessment$ of $\subspace$, then $\assessment$ avoids non-positivity relative to $\spacecone$ iff there is some maximal $\desirs\in\maxdesirs[\spacecone](\pspace)$ such that $\assessment\subseteq\desirs$.
\end{theorem}

\begin{corollary}\label{cor:maximality}
  Let\/ $\subspace$ be a linear subspace of\/ $\gambles(\pspace)$ and let $\cone\subset\subspace$ be a convex cone containing the zero gamble $0$.
  Consider any subset $\assessment$ of $\subspace$, then
  \begin{equation}\label{eq:atomicity}
    \natex[\spacecone]{\assessment}
    =\bigcap\set{\desirs\in\maxdesirs[\spacecone](\pspace)}
    {\assessment\subseteq\desirs}.
  \end{equation}
\end{corollary}

\section{An important special case}\label{sec:desir}
We now turn to the important special case, commonly considered in the literature \citep{walley2000,moral2003}, where $\subspace=\gambles(\pspace)$ and $\cone=\gambles^+_0(\pspace)$ is the cone of all non-negative gambles, for which the associated partial order $\succeq$ is the point-wise ordering $\geq$.
\par
If $\desirs$ avoids non-positivity relative to $\smash[b]{\bigl(\gambles(\pspace),\gambles^+_0(\pspace)\bigr)}$, we simply say that $\desirs$ \emph{avoids non-positivity}: $\gambles^-(\pspace)\cap\posi(\desirs)=\emptyset$.\footnote{A related, but weaker condition, is that $\desirs$ \emph{avoids partial loss}, meaning that $f\not<0$ for all gambles~$f$ in~$\posi(\desirs)$. We need the stronger condition because we have excluded the zero gamble from being desirable.}
\par
Similarly, if $\desirs$ is coherent relative to $\bigl(\gambles(\pspace),\gambles^+_0(\pspace)\bigr)$, we simply say that $\desirs$ is \emph{coherent}, and we denote the set of coherent sets of desirable gambles by $\alldesirs(\pspace)$.
In this case, the coherence conditions D\ref{item:Dnonzero}--D\ref{item:Dnonzeroapl} are to be seen as \emph{rationality criteria}.
In particular, D\ref{item:Dapg} is now an `accepting partial gain' condition: $\gambles^+_0(\pspace)\subseteq\desirs$.
D\ref{item:Dnonzeroapl} is an `avoiding partial loss' condition, together with the convention that the zero gamble is never considered desirable.
We give two illustrations, the first is a general one and the second models certainty about~$\omega$ happening.
The dashed line indicates a non-included border.
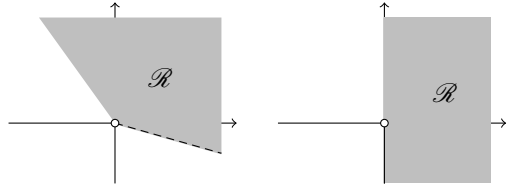
\begin{center}
  \begin{tikzpicture}[scale=2]
    \draw[->] (-.7,0) -- (.8,0);
    \draw[->] (0,-.4) -- (0,.8);
    \fill[lightgray] (-.5,.7) -- (0,0) -- (.7,-.2) -- (.7,.7) --cycle;
    \draw[lightgray,thick] (0,0) -- (.7,-.2);
    \draw[densely dashed] (0,0) -- (.7,-.2);
    \node[fill=white,draw,circle,inner sep=1pt] (zero) at (0,0) {};
    \node at (.3,.3) {$\desirs$};
  \end{tikzpicture}
  \hspace{1em}
  \begin{tikzpicture}[scale=2]
    \draw[->] (-.7,0) -- (.8,0);
    \draw[->] (0,-.4) -- (0,.8);
    \fill[lightgray] (0,-.4) rectangle (.7,.7);
    \draw (0,-.4) -- (0,0);
    \draw[lightgray,thick] (0,0) -- (0,.7);
    \node[fill=white,circle,draw,inner sep=1pt] (zero) at (0,0) {};
    \node at (.4,.2) {$\desirs$};
  \end{tikzpicture}
\end{center}
\par
The $\bigl(\gambles(\pspace),\gambles^+_0(\pspace)\bigr)$-natural extension of an assessment $\assessment\subseteq\gambles(\pspace)$ is simply denoted by $\natex{\assessment}$, and is called the \emph{natural extension} of $\assessment$.
In that case we can visualise Eq.~\eqref{eq:posi-natex}, i.e., natural extension as a positive hull operation, with a small illustration:
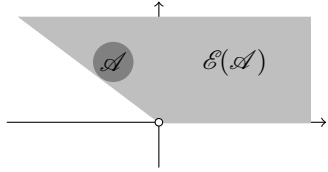
\begin{center}
  \begin{tikzpicture}[scale=2]
    \draw[->] (-1,0) -- (1.1,0);
    \draw[->] (0,-.3) -- (0,.8);
    \node[circle,fill=gray,inner sep=1pt] (A) at (-.3,.4) {$\assessment$};
    \coordinate (Atan) at (tangent cs:node=A,point={(0,0)},solution=2);
    \draw[lightgray] (0,0) -- (intersection of 0,0--Atan and -1,.7--1,.7) coordinate (Atanaway);
    \fill[lightgray] (Atanaway) -- (0,0) -- (1,0) -- (1,.7) -| (Atanaway);
    \draw[lightgray,thick] (0,0) -- (1,0);
    \node[fill=white,draw,circle,inner sep=1pt] (zero) at (0,0) {};
    \node at (.5,.4) {$\natex{\assessment}$};
    \node[circle,fill=gray,inner sep=1pt] (A) at (-.3,.4) {$\assessment$};
  \end{tikzpicture}
\end{center}
\par
More generally, consider a linear subspace $\subspace$ of $\gambles(\pspace)$, and let $\cone=\set{f\in\subspace}{f\geq0}=\gambles^+_0(\pspace)\cap\subspace$ be the cone of all non-negative gambles in that subspace.
If a set $\desirs\subseteq\subspace$ is coherent relative to $(\subspace,\gambles^+_0(\pspace)\cap\subspace)$, we also say that is \emph{coherent relative to $\subspace$}.
\par
In Secs.~\ref{sec:finite-polynomials} and~\ref{sec:infinite-xch}, we shall come across other useful types of coherence, relative to more involved structures $\subspace$ and $\cone$.

\subsection{Weakly desirable gambles, previsions \& marginally desirable gambles}
We now define \emph{weak desirability}: a useful modification of \citeauthor{walley1991}'s \citeyearpar[Section~3.7]{walley1991} notion of \emph{almost-desirability}.
Our conditions for a gamble~$f$ to be weakly desirable are more stringent than \citeauthor{walley1991}'s for almost-desirability: he only requires that adding any constant strictly positive amount of utility to~$f$ should make the resulting gamble desirable.
We require that adding anything desirable (be it constant or not) to~$f$ should make the resulting gamble desirable.
Weak desirability is better behaved under updating: we shall see in Proposition~\ref{prop:updating-preserves-exchangeability} that it makes sure that the exchangeability of a set of desirable gambles, whose definition hinges on the notion of weak desirability, is preserved under updating after observing a sample.
This is not necessarily true if weak desirability is replaced by almost-desirability in the definition of exchangeability, as was for instance done in our earlier work \citep{cooman2005c}.

\begin{definition}[Weak desirability]\label{def:weak-desirability}
  Consider a coherent set $\desirs$ of desirable gambles.
  Then a gamble $f$ is called \emph{weakly desirable} if $f+f'$ is desirable for all desirable $f'$, i.e., if
  $f+f'\in\desirs$ for all $f'$ in $\desirs$.
  We denote the set of weakly desirable gambles by $\weakly{\desirs}$:
  \begin{equation}\label{eq:weak-desirability}
    \weakly{\desirs}
    =\set{f\in\gambles(\pspace)}{f+\desirs\subseteq\desirs}.
  \end{equation}
\end{definition}
In particular, every desirable gamble is also weakly desirable, so $\desirs\subseteq\weakly{\desirs}$.

\begin{proposition}\label{prop:weak-desirability}
  Let $\desirs$ be a coherent set of desirable gambles, and let $\weakly{\desirs}$ be the associated set of weakly desirable gambles.
  Then $\weakly{\desirs}$ has the following properties, for all gambles $f$, $f_1$, and~$f_2$ in $\gambles(\pspace)$ and all real $\lambda\geq0$:
  \begin{compactenum}[\upshape {WD}1.]
    \item if $f<0$ then $f\notin\weakly{\desirs}$, or equivalently $\gambles^-_0(\pspace)\cap\weakly{\desirs}=\emptyset$ [avoiding partial loss];\footnote{Compare this to the less stringent requirement for almost-desirability \cite[Section~3.7.3]{walley1991}: if $f\in\weakly{\desirs}$ then $\sup f\geq 0$ [avoiding sure loss].}\label{item:WDapl}
    \item if $f\geq0$ then $f\in\weakly{\desirs}$, or equivalently $\gambles^+_0(\pspace)\subseteq\weakly{\desirs}$ [accepting partial gain];\label{item:WDapg}
    \item if $f\in\weakly{\desirs}$ then $\lambda f\in\weakly{\desirs}$ [scaling];\label{item:WDscl}
    \item if $f_1,f_2\in\weakly{\desirs}$ then $f_1+f_2\in\weakly{\desirs}$ [combination].\label{item:WDcmb}
  \end{compactenum}
\end{proposition}
Like $\desirs$, $\weakly{\desirs}$ is a convex cone.
\par
With a set of gambles $\assessment$, we associate a \emph{lower prevision} $\lpr_\assessment$ and an \emph{upper prevision} $\upr_\assessment$ by letting
\begin{align}
  \lpr_\assessment(f)
  &\coloneqq\sup\set{\mu\in\reals}{f-\mu\in\assessment} \label{eq:lpr}\\
  \upr_\assessment(f)
  &\coloneqq\inf\set{\mu\in\reals}{\mu-f\in\assessment}
\end{align}
for all gambles $f$.
Observe that $\lpr_\assessment$ and $\upr_\assessment$ always satisfy the \emph{conjugacy relation} $\lpr_\assessment(-f)=-\upr_\assessment(f)$.
We call a real functional $\lpr$ on $\gambles(\pspace)$ a \emph{coherent lower prevision} if  there is some coherent set of desirable gambles $\desirs$ on $\gambles(\pspace)$ such that $\lpr=\lpr_\desirs$.

\begin{theorem}\label{theo:lpr-desirs}
  Let $\desirs$ be a coherent set of desirable gambles.
  Then $\lpr_\desirs$ is real-valued, ${\lpr_\desirs=\lpr_{\weakly{\desirs}}}$, and $\lpr_\desirs(f)\geq0$ for all $f\in\weakly{\desirs}$.
  Moreover, a real functional $\lpr$ is a coherent lower prevision iff it satisfies the following properties, for all gambles $f_1$ and $f_2$ in $\gambles(\pspace)$ and all real $\lambda\geq0$:
  \begin{compactenum}[\upshape P1.]
    \item $\lpr(f)\geq\inf f$ [accepting sure gain];\label{item:Pasg}
    \item $\lpr(f_1+f_2)\geq\lpr(f_1)+\lpr(f_2)$ [super-additivity];\label{item:Psad}
    \item $\lpr(\lambda f)=\lambda\lpr(f)$ [non-negative homogeneity].\label{item:Phom}
  \end{compactenum}
\end{theorem}
\par
A coherent lower prevision $\lpr$ is called a \emph{linear prevision} if it is self-conjugate, in the sense that $\lpr=\upr$.
Such a linear prevision can be seen as an expectation operator associated with a (finitely additive) probability.
Using Proposition~\ref{prop:maximality}, it is not difficult to prove that the lower prevision $\lpr_\desirs$ associated with a maximal coherent set $\desirs$ is a linear prevision.\footnote{The proof for finite $\pspace$, given in \cite{couso2009isipta}, can be trivially extended to the infinite case. See also the proof of Proposition~\ref{prop:maximal-iid} in Appendix~\ref{app:proofs}.}

Finally, we turn to marginal desirability.
Given a coherent set of desirable gambles $\desirs$, we define the associated set of \emph{marginally desirable} gambles as
\begin{equation}
  \marginally{\desirs}\coloneqq\set{f-\lpr_\desirs(f)}{f\in\gambles(\pspace)}.
\end{equation}
The set of marginally desirable gambles $\marginally{\desirs}$ is completely determined by the lower prevision $\lpr_\desirs$.
The converse is also true:

\begin{proposition}\label{prop:marginal-desirability}
  Let $\desirs$ be a coherent set of desirable gambles.
  Then $\lpr_{\marginally{\desirs}}=\lpr_\desirs$ and
  \begin{equation}
    \marginally{\desirs}
    =\marginally{\lpr_\desirs}
    \coloneqq\set{f\in\gambles(\pspace)}{\lpr_\desirs(f)=0}.
  \end{equation}
\end{proposition}
The set of marginally desirable gambles $\marginally{\desirs}$ is the entire cone surface of $\desirs$ and $\weakly{\desirs}$, possibly including gambles that incur a partial (but not a sure) loss.
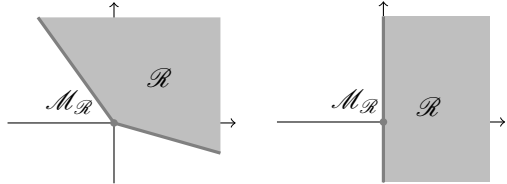
\begin{center}
  \begin{tikzpicture}[scale=2]
    \draw[->] (-.7,0) -- (.8,0);
    \draw[->] (0,-.4) -- (0,.8);
    \fill[lightgray] (0,0) -- (-.5,.7) -| (.7,-.2) -- cycle;
    \node at (.3,.3) {$\desirs$};
    \draw[gray,very thick] (-.5,.7) -- (0,0) -- (.7,-.2);
    \node[fill=gray,circle,inner sep=1pt] (zero) at (0,0) {};
    \node[above left,xshift=-4pt] at (0,0) {$\marginally{\desirs}$};
  \end{tikzpicture}
  \hspace{1em}
  \begin{tikzpicture}[scale=2]
    \draw[->] (-.7,0) -- (.8,0);
    \draw[->] (0,-.4) -- (0,.8);
    \fill[lightgray] (0,-.4) -| (.7,0) |-  (0,.7) -- cycle;
    \node at (.3,.1) {$\desirs$};
    \draw[gray,very thick] (0,-.4) -- (0,.7);
    \node[fill=gray,circle,inner sep=1pt] (zero) at (0,0) {};
    \node[above left,xshift=2pt] at (0,0) {$\marginally{\desirs}$};
  \end{tikzpicture}
\end{center}
When $\desirs$ is maximal, $\marginally{\desirs}$ constitutes a hyperplane.

\subsection{Updating sets of desirable gambles}\label{sec:updating}
Consider a set of desirable gambles~$\desirs$ on~$\pspace$.
With a non-empty subset~$B$ of~$\pspace$, we associate an \emph{updated} set of desirable gambles on $\pspace$, as defined by \citet{walley2000}:
\begin{equation}
  \walleyupdate{\desirs}{B}
  \coloneqq\set{f\in\gambles(\pspace)}{I_Bf\in\desirs}.
\end{equation}
We find it more convenient to work with the following, slightly different but completely equivalent, version:
\begin{align}
  \update{\desirs}{B}
  \coloneqq&\set{f\in\desirs}{I_Bf=f}
  =\desirs\cap\update{\gambles(\pspace)}{B},
\end{align}
which completely determines $\walleyupdate{\desirs}{B}$: for all $f\in\gambles(\pspace)$,
\begin{equation}
  f\in\walleyupdate{\desirs}{B}\asa I_Bf\in\update{\desirs}{B}.
\end{equation}
In our version, updating corresponds to intersecting the convex cone $\desirs$ with the linear subspace $\update{\gambles(\pspace)}{B}$, which results in a convex cone $\update{\desirs}{B}$ of lower dimension.
And since we can uniquely identify a gamble $f=I_Bf$ in $\update{\gambles(\pspace)}{B}$ with a gamble on~$B$, namely its restriction~$f_B$ to~$B$, and {\itshape vice versa}, we can also identify $\update{\desirs}{B}$ with a set of desirable gambles on~$B$:
\begin{equation}
  \restrict{\desirs}{B}
  \coloneqq\set{f_B}{f\in\update{\desirs}{B}}
  =\set{f_B}{f\in\walleyupdate{\desirs}{B}}
  \subseteq\gambles(B).
\end{equation}

\begin{proposition}\label{prop:updating-preserves-coherence}
  If $\desirs$ is a coherent set of desirable gambles on $\pspace$, then $\update{\desirs}{B}$ is coherent relative to $\update{\gambles(\pspace)}{B}$, or equivalently, $\restrict{\desirs}{B}$ is a coherent set of desirable gambles on~$B$.
\end{proposition}

\noindent Our subject takes $\update{\desirs}{B}$ (or $\restrict{\desirs}{B}$) as his set of desirable gambles contingent on observing the event~$B$.

\section{Finite exchangeable sequences}\label{sec:finite-xch}
Now that we have become better versed in the theory of sets of desirable gambles, we are going to focus on the first main topic: reasoning about finite exchangeable sequences.
We first show how they are related to count vectors (Section~\ref{sec:countvec}).
Then we are ready to give a desirability-based definition of exchangeability (Section~\ref{sec:xchdef}) and treat natural extension and updating under exchangeability (Sections~\ref{sec:exchangeable-natural-extension} and \ref{sec:exchangeable-update}).
After presenting our Finite Representation Theorem (Section~\ref{sec:finite-representation}), we can show what natural extension and updating under exchangeability look like in terms of the count vector representation (Sections~\ref{sec:finite-representation-exnatex} and~\ref{sec:finite-representation-update}).
Finally, we take a look at multinomial processes (Section~\ref{sec:multproc}), which will allow us to present a version of the Representation Theorem in terms of frequency vectors (Section~\ref{sec:finite-polynomials}).

Consider random variables $\rv_1$, \dots, $\rv_N$ taking values in a non-empty finite set~$\values$,\footnote{A lot of functions and sets introduced below will depend on the set~$\values$. We do not indicate this explicitly, to not overburden the notation and because we do not consider different sets of values in this paper.} where $N\in\nats_0$, i.e., a positive (non-zero) integer.
The possibility space is $\pspace=\values^N$.

\subsection{Count vectors}\label{sec:countvec}
We denote by $\sample{x}=\tuple{x}{N}$ an arbitrary element of~$\values^N$.
$\permuts_N$~is the set of all permutations~$\pi$ of the index set $\{1,\dots,N\}$.
With any such permutation~$\pi$, we associate a permutation of $\values^N$, also denoted by~$\pi$, and defined by $(\pi\sample{x})_k=\smash[b]{\sample{x}_{\pi(k)}}$, or in other words, $\pi\tuple{x}{N}=(x_{\pi(1)},\dots,x_{\pi(N)})$.
Similarly, we lift~$\pi$ to a permutation $\pi^t$ of $\gambles(\values^N)$ by letting $\pi^tf=f\circ\pi$, so $(\pi^tf)(\sample{x})=f(\pi\sample{x})$.
\par
The permutation invariant atoms $\batom{x}\coloneqq\set{\pi\sample{x}}{\pi\in\permuts_N}$, $x\in\values^N$ are the smallest permutation invariant subsets of $\values^N$.
We introduce the \emph{counting map}
\begin{equation}
  \vbcntf{N}\colon\values^N\to\vcnts{N}
  \colon\sample{x}\mapsto\vbcntf{N}(\sample{x}),
\end{equation}
where $\vbcntf{N}(\sample{x})$ is the $\values$-tuple with components
\begin{equation}
 \cntf_z^N(\sample{x})
 \coloneqq\abs{\set{k\in\{1,\dots,N\}}{x_k=z}}
 \;\;\text{ for all $z\in\values$},
\end{equation}
and the set of possible \emph{count vectors} is given by
\begin{equation}
  \vcnts{N}\coloneqq\biggset{\cnt{m}\in\nats^\values}{\smashoperator{\sum_{x\in\values}}m_x=N}.
\end{equation}
If $\cnt{m}=\vbcntf{N}(\sample{x})$, then $\batom{x}=\set{\sample{y}\in\values^N}{\vbcntf{N}(\sample{y})=\cnt{m}}$,
so the atom $\batom{x}$ is completely determined by the count vector $\cnt{m}$ of all its the elements, and is therefore also denoted by $\batom{m}$.

\begin{example}[Running example]
  To familiarise ourselves with some of the concepts introduced, we will use a running example, in which we assume a sample space $\values\coloneqq\{b,w\}$ (for \emph{black} and \emph{white} -- the colours of marbles drawn from an urn containing a mixture of such marbles).
  Consider the situation $N\coloneqq2$, then
  \begin{equation*}
    \values^N = \{b,w\}^2 = \bigl\{(b,b),(b,w),(w,b),(w,w)\bigr\} \quad\text{ and }\quad
    \vcnts{N} = \bigl\{(2,0),(1,1),(0,2)\bigr\}.
  \end{equation*}
  Their correspondence and the non-trivial permutations are implicitly given by
  \begin{align*}
    \batom{2,0} = \bigl\{(b,b)\bigr\}, && \batom{1,1} = \bigl\{(b,w),(w,b)\bigr\}, && \batom{0,2} = \bigl\{(w,w)\bigr\}.\mspace{20mu}\blacklozenge\mspace{-20mu}
  \end{align*}
\end{example}

\subsection{Defining exchangeability}\label{sec:xchdef}
If a subject assesses that $\rv_1$, \dots, $\rv_N$ are exchangeable, this means that for any gamble $f$ and any permutation~$\pi$, he finds exchanging $\pi^tf$ for $f$ weakly desirable, because he is indifferent between them \citep[cf.][Section~4.1.1]{walley1991}.
Let
\begin{equation}
  \weakly{\permuts_N}
  \coloneqq\set{f-\pi^tf}{f\in\gambles(\values^N)\text{ and }\pi\in\permuts_N},
\end{equation}
then we should have that $\weakly{\permuts_N}\subseteq\weakly{\desirs}$.\footnote{Note that the gambles in $\weakly{\permuts_N}$ cannot be assumed to be desirable, because $\weakly{\permuts_N}$ does not avoid non-positivity.}
This is the basis for our definition of exchangeability.
\par
What we would like to do now, is to provide alternative characterisations of exchangeability.
These will be useful for the further development, and provide additional insight into what an assessment of exchangeability amounts to.
\par
We begin by defining a special linear transformation $\vex{N}$ of the linear space of gambles $\gambles(\values^N)$:
\begin{equation}\label{eq:ex}
  \vex{N}\colon\gambles(\values^N)\to\gambles(\values^N)
  \colon f\mapsto\vex{N}(f)\coloneqq\frac{1}{N!}\smashoperator[r]{\sum_{\pi\in\permuts_N}}\pi^tf.
\end{equation}
The idea behind this linear transformation $\vex{N}$ is that it renders a gamble $f$ insensitive to permutation by replacing it with the uniform average $\vex{N}(f)$ of all its permutations $\pi^tf$.
Indeed, observe that for all gambles $f$ and all permutations $\pi$:
\begin{equation}\label{eq:ex-permutation}
  \vex{N}(\pi^tf)=\vex{N}(f)
  \quad\text{ and }\quad
  \pi^t\bigl(\vex{N}(f)\bigr)=\vex{N}(f).
\end{equation}
So $\vex{N}(f)$ is permutation invariant and therefore constant on the permutation invariant atoms $\batom{m}$, and it assumes the same value for all gambles that can be related to each other through some permutation.
But then, what is the value that $\vex{N}(f)$ assumes on each such permutation invariant atom $\batom{m}$?
It is not difficult to see that
\begin{equation}\label{eq:ex-muhy}
  \vex{N}
  =\smashoperator{\sum_{\cnt{m}\in\vcnts{N}}}\vmuhy{N}(\cdot\vert\cnt{m})I_{\batom{m}},
\end{equation}
or in other words, the gamble $\vex{N}(f)$ assumes the constant value $\vmuhy{N}(f\vert\cnt{m})$ on $\batom{m}$,  where we let
\begin{equation}
  \vmuhy{N}(f\vert\cnt{m})
  \coloneqq\frac{1}{\abs{\batom{m}}}\smashoperator[r]{\sum_{\sample{y}\in\batom{m}}}f(\sample{y})
  \quad\text{ and }\quad
  \abs{\batom{m}}
  =\binom{N}{\cnt{m}}
  \coloneqq\frac{N!}{\prod_{z\in\values}m_z!}.
\end{equation}
$\vmuhy{N}(\cdot\vert\cnt{m})$ is the linear expectation operator associated with the uniform distribution on the invariant atom $\batom{m}$.
It characterises a (multivariate) \emph{hyper-geometric distribution} \citep[Section~39.2]{johnson1997}, associated with random sampling without replacement from an urn with $N$ balls of types $\values$, whose composition is characterised by the count vector $\cnt{m}$.

\begin{example}
  To get some feeling for what this means, let us go back to our running example $\values\coloneqq\{b,w\}$ and let $N\coloneqq4$ and $m=(m_b,m_w)\coloneqq(2,2)$.
  Then
  \begin{equation*}
    \batom{m}=[2,2]=\{(w,w,b,b),(w,b,w,b),(w,b,b,w),(b,b,w,w),(b,w,w,b),(b,w,b,w)\}.
  \end{equation*}
  Consider the event $A$ that amongst the first three observations, we see twice $b$ and once $w$:
  \begin{equation*}
    A=\{(b,b,w),(b,w,b),(w,b,b)\}\times\{b,w\}.
  \end{equation*}
  Then
  \begin{equation*}
    \muhy^4(I_A\vert2,2)
    =\frac{1}{\binom{4}{2}}\sum_{y\in[2,2]}I_A(y)=\frac{1}{6}\cdot3=\frac{1}{2}
  \end{equation*}
  is the probability of getting two black balls and one white when drawing three balls without replacement from an urn containing two black and two white balls.~$\blacklozenge$
\end{example}

\noindent So we see that the linear transformation $\vex{N}$ is intricately linked with the $N$-variate hypergeometric distribution.
If we also observe that $\vex{N}\circ\vex{N}=\vex{N}$, we see that \emph{$\vex{N}$ is the linear projection operator} of the linear space $\gambles(\values^N)$ to the linear subspace
\begin{equation}
  \gambles_{\permuts_N}(\values^N)
  \coloneqq\set{f\in\gambles(\values^N)}{(\forall\pi\in\permuts_N)\pi^tf=f}
\end{equation}
of all permutation invariant gambles.
\par
The linear transformation $\vex{N}$ is also tightly connected with the previously defined set $\weakly{\permuts_N}$  of gambles $f-\pi^tf$  that play a role in defining exchangeability.
Indeed, if we look at the linear subspace $\weakly{\average_N}$ that is generated by such gambles $f-\pi^tf$, then it is not hard to see that
\begin{align}
  \weakly{\average_N}
  \coloneqq{}&\opspanning(\weakly{\permuts_N})\label{eq:span}
  =\biggset{\sum_{k=1}^n\lambda_kf_k}{n\geq0,\,\lambda_k\in\reals,\,f_k\in\weakly{\permuts_N}}\\
  ={}&\set{f-\vex{N}(f)}{f\in\gambles(\values^N)} 
  =\set{f\in\gambles(\values^N)}{\vex{N}(f)=0},\label{eq:kernel} 
\end{align}
where `$\opspanning(\cdot)$' denotes \emph{linear span} of its argument set: the set of all linear combinations of elements from that set.
The last equality tells us that the linear subspace $\weakly{\average_N}$ is the \emph{kernel} of the linear projection operator $\vex{N}$: it contains precisely those gambles that are mapped to $0$ by $\vex{N}$.

\begin{example}
  Let us return for a moment to our running example $\values\coloneqq\{b,w\}$ and $N\coloneqq2$, then
  \begin{equation}\label{eq:ex-weakly-permuts}
    \weakly{\permuts_N} = \set{f\in\gambles(\values^N)}{f(b,b)=f(w,w)=0\text{ and }f(b,w)=-f(w,b)},
  \end{equation}
  and $\weakly{\average_N}=\weakly{\permuts_N}$.
  Let $f$ be some gamble on $\values^N$ and $f'\coloneqq\vex{N}(f)$, then
  \begin{equation*}
  \begin{split}
    f'(b,b) &= \vmuhy{N}(f\vert2,0) = f(b,b),\\ f'(w,w) &= \vmuhy{N}(f\vert0,2) = f(w,w),\\ f'(b,w) = f'(w,b) &= \vmuhy{N}(f\vert1,1) = \tfrac{1}{2}\bigl(f(b,w)+f(w,b)\bigr).
  \end{split}
  \end{equation*}
  The plane $\update{\gambles(\values^N)}{\batom{1,1}}=\set{f\in\gambles(\values^N)}{f(b,b)=f(w,w)=0}$, which includes $\weakly{\permuts_N}$ and thus $\weakly{\average_N}$, can be used for a graphical illustration:
  
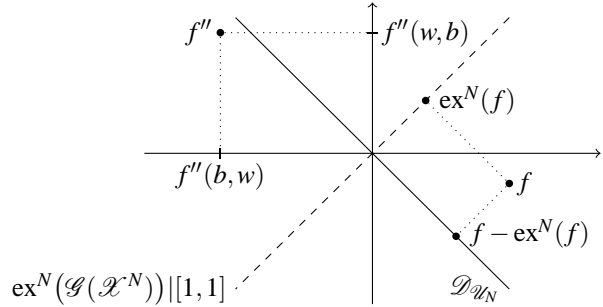
\begin{center}
    \begin{tikzpicture}[scale=2]
      \draw[->] (-1.5,0) -- (1.5,0);
      \draw[->] (0,-1) -- (0,1);
      \draw (-.9,.9) -- (.9,-.9) node[left] {$\weakly{\average_N}$};
      \node[fill,circle,inner sep=1pt,label={[inner sep=1pt]left:$f''$}] (fdp) at (-1,.8) {};
      \draw[dotted] (fdp) -- (fdp |- 1,0) coordinate (fdpbw) \xtick node[below] {$f''(b,w)$};
      \draw[dotted] (fdp) -- (fdp -| 0,1) coordinate (fdpwb) \ytick node[right] {$f''(w,b)$};
      \draw[dashed] (.9,.9) -- (-.9,-.9) node[left] {$\update{\vex{N}\bigl(\gambles(\values^N)\bigr)}{\batom{1,1}}$};
      \node[fill,circle,inner sep=1pt,label={[inner sep=1pt]right:$f$}] (f) at (.9,-.2) {};
      \node[fill,circle,inner sep=1pt,label={right:$\vex{N}(f)$}] (fp) at (.35,.35) {};
      \node[fill,circle,inner sep=1pt,label={[yshift=1pt]right:$f-\vex{N}(f)$}] (fminfp) at (.55,-.55) {};
      \draw[dotted] (f) -- (fp) (f) -- (fminfp);
    \end{tikzpicture}
  \end{center}
  We show the intersection of this plane and the range of the operator $\vex{N}$ as well as its effect on a gamble $f$ in $\update{\gambles(\values^N)}{\batom{1,1}}$.
  The gamble $f''$ is included to make it explicit which components are actually depicted.~$\blacklozenge$
\end{example}

\begin{definition}[Exchangeability]\label{def:exchangeability}
  A coherent set  $\desirs$ of desirable gambles on $\values^N$ is called \emph{exchangeable} if any (and hence all) of the following equivalent conditions is (are) satisfied:
  \begin{compactenum}[\upshape (i)]
    \item all gambles in $\weakly{\permuts_N}$ are weakly desirable: $\weakly{\permuts_N}\subseteq\weakly{\desirs}$;\label{item:xch-permutssubsetdef}
    \item $\weakly{\permuts_N}+\desirs\subseteq\desirs$;\label{item:xch-permutsreldef}
    \item all gambles in $\weakly{\average_N}$ are weakly desirable: $\weakly{\average_N}\subseteq\weakly{\desirs}$;\label{item:xch-averagesubsetdef}
    \item $\weakly{\average_N}+\desirs\subseteq\desirs$;\label{item:xch-averagereldef}
  \end{compactenum}
  We call a coherent lower prevision $\lpr$ on $\gambles(\values^N)$ \emph{exchangeable} if there is some exchangeable coherent set of desirable gambles $\desirs$ such that $\lpr=\lpr_\desirs$.
\end{definition}

\noindent Because they are stated in terms of the kernel $\weakly{\average_N}$ of the projection operator $\vex{N}$, which we have seen is intricately linked with (multivariate) hypergeometric distributions, conditions (\ref{item:xch-averagesubsetdef}) and (\ref{item:xch-averagereldef}) of this definition are quite closely related to the desirability version of a de Finetti-like representation theorem for finite exchangeable sequences in terms of sampling without replacement from an urn.
They allow us to talk about exchangeability without explicitly invoking permutations.
This is what we will address in Section~\ref{sec:finite-representation}.

\begin{example}
  In our running example, if~$f\in\gambles(\values^N)$ is desirable, then because of the definition of exchangeability and Eq.~\eqref{eq:ex-weakly-permuts}, all gambles in the linear subspace
  \begin{multline*}
    \bigl\{f'\in\gambles(\values^N): f'(b,b)=f(b,b), f'(w,w)=f(w,w),\\ f'(b,w)+f'(w,b)=f(b,w)+f(w,b)\bigr\}
  \end{multline*}
  are also desirable.
  So under exchangeability transfers between $(b,w)$ and $(w,b)$-components are irrelevant for desirability.
  This illustrates that, generally and geometrically speaking, any
  exchangeable set of desirable gambles
  $\desirs\subseteq\gambles(\values^N)$ for any $\values$ and any $N>1$ must be cylindrical along the directions in the linear subspace $\weakly{\average_N}$.~$\blacklozenge$
\end{example}

A number of useful results follow from Definition~\ref{def:exchangeability}:

\begin{proposition}\label{prop:permutability}
  Let $\desirs$ be a coherent set of desirable gambles.
  If $\desirs$ is exchangeable then it is also \emph{permutable}: $\pi^tf\in\desirs$ for all $f\in\desirs$ and all $\pi\in\permuts_N$.
\end{proposition}

\noindent We have seen above that the gambles $f-\pi^tf$ in $\weakly{\permuts_N}$ span the kernel $\weakly{\average_N}$ of the linear operator $\vex{N}$ that projects any gamble $f$ on its symmetrised counterpart $\vex{N}(f)$.
It should therefore not come as a surprise that for an exchangeable model, we can determine whether a gamble is desirable by looking at this symmetrised counterpart.

\begin{proposition}\label{prop:exchangeability-ex}
  Let $\desirs$ be a coherent and exchangeable set of desirable gambles.
  For all gambles $f$ and $f'$ on $\values^N$:
  \begin{compactenum}[\upshape (i)]
    \item $f\in\desirs\Leftrightarrow\vex{N}(f)\in\desirs$;
    \item If $\vex{N}(f)=\vex{N}(f')$, then  $f\in\desirs\Leftrightarrow f'\in\desirs$.
  \end{compactenum}
\end{proposition}
\noindent It follows from this last proposition and Eq.~\eqref{eq:kernel} that for any coherent and exchangeable set of desirable gambles $\desirs$:
\begin{equation}\label{eq:no-intersection}
  \desirs\cap\weakly{\average_N}=\emptyset.
\end{equation}

\noindent We can use these ideas to derive a direct characterisation for the exchangeability of a lower prevision, without the intervention of sets of desirable gambles.

\begin{theorem}\label{theo:exchangeability-lpr}
  Let $\lpr$ be a coherent lower prevision on $\gambles(\values^N)$.
  Then the following statements are equivalent:\footnote{This also shows that the exchangeability of a lower prevision can also be expressed using marginally desirable gambles \citep[see][Section~3.1.1]{Quaeghebeur2009phd}.}
  \begin{compactenum}[\upshape (i)]
    \item $\lpr$ is exchangeable;\label{item:xch-prev}
    \item $\lpr(f)=\upr(f)=0$ for all $f\in\weakly{\permuts_N}$;\label{item:xch-prevpermuts}
    \item $\lpr(f)=\upr(f)=0$ for all $f\in\weakly{\average_N}$.\label{item:xch-prevaverage}
  \end{compactenum}
\end{theorem}

\subsection{Exchangeable natural extension}\label{sec:exchangeable-natural-extension}
Let us denote the set of all coherent and exchangeable sets of desirable gambles on~$\values^N$ by
\begin{equation}
  \allexdesirs(\values^N)
  \coloneqq\set{\desirs\in\alldesirs(\values^N)}
  {\weakly{\average_N}+\desirs\subseteq\desirs}.
\end{equation}
This set is closed under arbitrary non-empty intersections.
We shall see further on in Corollary~\ref{cor:smallest-exchangeable} that it is also non-empty, and therefore has a smallest element.
\par
Suppose our subject has an assessment, or in other words, a set $\assessment$ of gambles on $\values^{N}$ that he finds desirable.
Then we can ask if there is some coherent and exchangeable set of desirable gambles~$\desirs$ that includes~$\assessment$.
In other words, we want a set of desirable gambles~$\desirs$ to satisfy the requirements:
\begin{inparaenum}[(i)]
  \item $\desirs$ is coherent;
  \item $\assessment\subseteq\desirs$; and
  \item $\weakly{\average_N}+\desirs\subseteq\desirs$.
\end{inparaenum}
The intersection $\bigcap_{i\in I}\desirs_i$ of an arbitrary non-empty family of sets of desirable gambles~$\desirs_i$, $i\in I$ that satisfy these requirements, will satisfy these requirements as well.
This is the idea behind the following definition and results.

\begin{definition}[Avoiding non-positivity under exchangeability]
  We say that a set  $\assessment$ of gambles on  $\values^N$ \emph{avoids non-positivity under exchangeability} if $[\gambles^+_0(\values^N)\cup\assessment]+\weakly{\average_N}$ avoids non-positivity.
\end{definition}

\begin{proposition}\label{prop:anpune}
  \begin{compactenum}[(i)]
    \item $\emptyset$ avoids non-positivity under exchangeability; 
    \item A non-empty set of gambles $\assessment$ on $\values^N$ avoids non-positivity under exchangeability iff $\assessment+\weakly{\average_N}$ avoids non-positivity.
  \end{compactenum}
\end{proposition}

\begin{example}
  For our running example, avoiding non-positivity under exchangeability is best illustrated graphically---again in the plane $\update{\gambles(\values^N)}{\batom{1,1}}$---for the case that the given assessment $\assessment\subset\update{\gambles(\values^N)}{\batom{1,1}}$ avoids non-positivity, but not so under exchangeability:
  
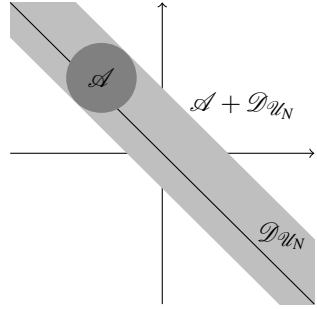
\begin{center}
    \begin{tikzpicture}[scale=2]
      \clip (-2,-1) rectangle (2,1);
      \begin{pgfonlayer}{background}
        \draw[->] (-1,0) -- (1,0);
        \draw[->] (0,-1) -- (0,1);
      \end{pgfonlayer}
      \draw (-1,1) -- (1,-1) node[pos=.85,above,xshift=5pt] {$\weakly{\average_N}$};
      \node[circle,fill=gray,inner sep=5pt] (A) at (-.4,.5) {$\assessment$};
      \coordinate[shift={(-2,2)}] (ce1) at (A.north east);
      \coordinate[shift={(-2,2)}] (ce2) at (A.south west);
      \coordinate[shift={(4,-4)}] (ce3) at (A.north east);
      \coordinate[shift={(4,-4)}] (ce4) at (A.south west);
      \begin{pgfonlayer}{background}
        \clip (-1,-1) rectangle (1,1);
        \fill[lightgray] (ce1) -- (ce3) -- (ce4) -- (ce2);
      \end{pgfonlayer}
      \path (ce1) -- (ce3) node[pos=.45,right] {$\assessment+\weakly{\average_N}$};
    \end{tikzpicture}
  \end{center}
  (As a reminder: $\update{\gambles(\values^N)}{\batom{1,1}}=\set{f\in\gambles(\values^N)}{f(b,b)=f(w,w)=0}$.)~$\blacklozenge$
\end{example}

\begin{theorem}[Exchangeable natural extension]\label{theo:exnatex}
  Consider a set $\assessment$ of gambles on $\values^N$, and define its \emph{exchangeable natural extension} $\exnatex{\assessment}$ by
  \begin{equation}
    \exnatex{\assessment}
    \coloneqq{} \bigcap\set{\desirs\in\allexdesirs(\values^N)}{\assessment\subseteq\desirs}.
  \end{equation}
  Then the following statements are equivalent:
  \begin{compactenum}[\upshape (i)]
  \item\label{item:Aaplunderxch} $\assessment$ \emph{avoids non-positivity under exchangeability};
  \item\label{item:Asubsetcohxch} $\assessment$ is included in some coherent and exchangeable set of desirable gambles;
  \item\label{item:exnatexA-nonempty} $\exnatex{\assessment}\neq\gambles(\values^N)$;
  \item\label{item:exnatexA-coh} $\exnatex{\assessment}$ is a coherent and exchangeable set of desirable gambles;
  \item\label{item:exnatexA-smallestcoh} $\exnatex{\assessment}$ is the smallest coherent and exchangeable set of desirable gambles that includes $\assessment$.
  \end{compactenum}
  When any (and hence all) of these equivalent statements hold, then
  \begin{align}
    \exnatex{\assessment}
    &=\posi\left(\weakly{\average_N}+[\gambles^+_0(\values^N)\cup\assessment]\right)
    \label{eq:posi-exnatex}\\
    &=\weakly{\average_N}+\natex{\assessment}.
    \label{eq:posi-exnatex-better}
  \end{align}
\end{theorem}

\begin{example}
  Exchangeable natural extension is again best illustrated graphically.
  We use an assessment~$\assessment\subset\update{\gambles(\values^N)}{\batom{1,1}}$ that avoids non-positivity under exchangeability and contrast natural extension with exchangeable natural extension:
  
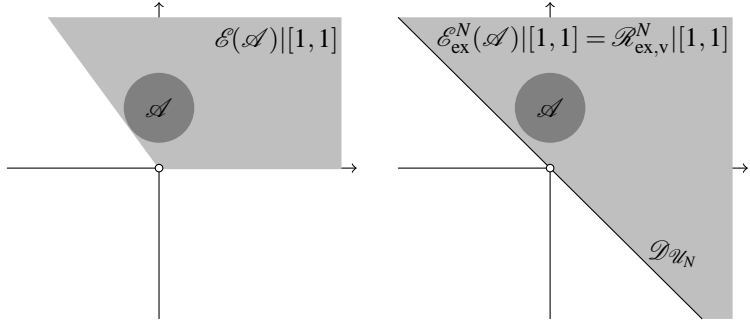
\begin{center}
    \begin{tikzpicture}[scale=2]
      \draw[->] (-1,0) -- (1.3,0);
      \draw[->] (0,-1) -- (0,1.1);
      \node[circle,fill=gray,inner sep=5pt] (A) at (0,.4) {$\assessment$};
      \coordinate (Atan) at (tangent cs:node=A,point={(0,0)},solution=2);
      \draw[lightgray] (0,0) -- (intersection of 0,0--Atan and -1,1--0,1) coordinate (Atanaway);
      \fill[lightgray] (Atanaway) -- (0,0) -- (1.2,0) -- (1.2,.8) |- (Atanaway);
      \draw[lightgray,thick] (0,0) -- (1.2,0);
      \node[fill=white,draw,circle,inner sep=1pt] (zero) at (0,0) {};
      \node[left,inner sep=1pt] at (1.2,.85) {$\update{\natex{\assessment}}{\batom{1,1}}$};
      \node[circle,fill=gray,inner sep=5pt] at (0,.4) {$\assessment$};
    \end{tikzpicture}
    \hspace{1em}
    \begin{tikzpicture}[scale=2]
      \draw[->] (-1,0) -- (1.3,0);
      \draw[->] (0,-1) -- (0,1.1);
      \fill[lightgray] (-1,1) -- (1,-1) -- (1.2,-1) -- (1.2,1) -- cycle;
      \draw (-1,1) -- (1,-1) node[pos=.86,above,xshift=5pt] {$\weakly{\average_N}$};
      \node[circle,fill=gray,inner sep=5pt] (A) at (0,.4) {$\assessment$};
      \node[fill=white,draw,circle,inner sep=1pt] (zero) at (0,0) {};
      \node[left,inner sep=1pt] at (1.2,.85)  {$\update{\exnatex{\assessment}}{\batom{1,1}}=\update{\vacuous}{\batom{1,1}}$};
    \end{tikzpicture}
  \end{center}
  Only the parts of an assessment that fall outside of $\vacuous$ have a nontrivial impact.
  This is something that cannot be illustrated in the context of drawings as these, but think about what the extensions would look like if the assessment would consist of the singleton $\{f\}$ with $f(b,b)=-f(w,w)=1$ and $f(b,w)=f(w,b)=0$: the natural extension $\natex{\{f\}}$ consists of all gambles $\mu\lambda f+(1-\mu)g$, with $\mu\in[0,1]$, $\lambda\in\reals^+_0$, and $g\in\gambles^+_0(\values^N)$, so the exchangeable natural extension $\exnatex{\{f\}}$ consists of all gambles $\rho h+\mu\lambda f+(1-\mu)g$, where additionally $\rho\in\reals$ and $h(b,w)=-h(w,b)=1$ and $h(b,b)=h(w,w)=0$.~$\blacklozenge$
\end{example}

\noindent Eq.~\eqref{eq:posi-exnatex-better} shows that if we have an assessment $\assessment$ with a finite description, it is possible, but not necessarily very efficient, to represent its exchangeable natural extension on a computer: besides the finite description of its extreme rays of $\natex{\assessment}$, we need to account for taking the Minkowski sum with $\weakly{\average_N}$.
We shall see further on in Theorem~\ref{theo:finite-representation-exnatex} that this extra complication can be circumvented by working with so-called count representations.
\par
There is always a most conservative exchangeable belief model, which represents the effects of making only an assessment of exchangeability, and nothing more:

\begin{corollary}\label{cor:smallest-exchangeable}
  The set\/ $\allexdesirs(\values^N)$ is non-empty, and has a smallest element
  \begin{equation}\label{eq:smallest-exchangeable}
    \vacuous
    \coloneqq\exnatex{\emptyset}
    =\weakly{\average_N}+\gambles^+_0(\values^N).
  \end{equation}
\end{corollary}

\subsection{Updating exchangeable models}\label{sec:exchangeable-update}
Consider an exchangeable and coherent set of desirable gambles $\desirs$ on $\values^N$, and assume that we have observed the values $\osample{x}=(\obs{x}_1,\obs{x}_2,\dots,\obs{x}_{\obs{n}})$ of the first $\obs{n}$ variables $\rv_1$, \dots, $\rv_{\obs{n}}$, and that we want to make inferences about the remaining $\rest{n}\coloneqq N-\obs{n}$ variables.
To do this, we simply update the set $\desirs$ with the event $\event_{\osample{x}}\coloneqq\{\osample{x}\}\times\values^{\rest{n}}$, to obtain the set $\update{\desirs}{\event_{\osample{x}}}$, also denoted as $\update{\desirs}{\osample{x}}\coloneqq\set{f\in\desirs}{fI_{\event_{\osample{x}}}=f}$.
As we have seen in Section~\ref{sec:updating}, this set can be identified with a coherent set of desirable gambles on $\values^{\rest{n}}$, which we denote by $\restrict{\desirs}{\osample{x}}$.
With obvious notations:\footnote{Here and further on we silently use cylindrical extension on gambles, i.e., let them `depend' on extra variables whose value does not influence the value they take.}
\begin{equation}
  \restrict{\desirs}{\osample{x}}
  \coloneqq\set{f\in\gambles(\values^{\rest{n}})}{fI_{\event_{\osample{x}}}\in\desirs}.
\end{equation}
We already know that updating preserves coherence.
We now see that this type of updating on an observed sample also preserves exchangeability.

\begin{proposition}\label{prop:updating-preserves-exchangeability}
  Consider $\osample{x}\in\values^{\obs{n}}$ and a coherent and exchangeable set of desirable gambles~$\desirs$ on $\values^N$.
  Then $\restrict{\desirs}{\osample{x}}$ is a coherent and exchangeable set of desirable gambles on~$\values^{\rest{n}}$.
\end{proposition}

We also introduce another type of updating, where we observe a count vector $\ocnt{m}\in\vcnts{\obs{n}}$, and we update the set $\desirs$ with the set $\event_{\ocnt{m}}\coloneqq\oatom{m}\times\values^{\rest{n}}$, to obtain the set $\update{\desirs}{\event_{\ocnt{m}}}$, also denoted as $\update{\desirs}{\ocnt{m}}\coloneqq\set{f\in\desirs}{fI_{\event_{\ocnt{m}}}=f}$.
This set can be identified with a coherent set of desirable gambles on $\values^{\rest{n}}$, which we also denote by $\restrict{\desirs}{\ocnt{m}}$.
With obvious notations:
\begin{equation}
  \restrict{\desirs}{\ocnt{m}}
  \coloneqq\set{f\in\gambles(\values^{\rest{n}})}{fI_{\event_{\ocnt{m}}}\in\desirs}.
\end{equation}
Interestingly, the count vector $\ocnt{m}$ for an observed sample $\osample{x}$ is a \emph{sufficient statistic} in that it extracts from $\osample{x}$ all the information that is needed to characterise the updated model:

\begin{proposition}[Sufficiency of observed count vectors]\label{prop:sufficiency}
  Consider $\osample{x},\osample{y}\in\values^{\obs{n}}$ and a coherent and exchangeable set of desirable gambles~$\desirs$ on $\values^N$.
  If\/ $\osample{y}\in\oatom{x}$, or in other words if\/ $\vbcntf{\obs{n}}(\osample{x})=\vbcntf{\obs{n}}(\osample{y})\eqqcolon\ocnt{m}$, then $\restrict{\desirs}{\osample{x}}=\restrict{\desirs}{\osample{y}}=\restrict{\desirs}{\ocnt{m}}$.
\end{proposition}

\subsection{Finite representation}\label{sec:finite-representation}
We can use the symmetry that an assessment of exchangeability generates to represent an exchangeable coherent set of desirable gambles in a much more economical, or condensed, fashion.
This has already been made apparent in Proposition~\ref{prop:exchangeability-ex}, where we saw that the desirability of any gamble $f$ can be determined by looking at the desirability of its symmetrised counterpart $\vex{N}(f)$.
We have seen that this projection $\vex{N}(f)$ of $f$ onto the linear subspace $\gambles_{\permuts_N}(\values^N)$ of permutation invariant gambles assumes the constant value $\vmuhy{N}(f\vert m)$ on the permutation invariant atoms $\batom{m}$, $m\in\cnts^N$.
\par
Now, since a gamble is permutation invariant if and only if it is constant on these permutation invariant atoms, we can \emph{identify} permutation invariant gambles on $\values^N$ with gambles on $\cnts^N$.
This identification is made more formal using the following linear isomorphism $\vocntf{N}$ between the linear spaces $\gambles(\vcnts{N})$ and $\gambles_{\permuts_N}(\values^N)$:
\begin{equation}
  \vocntf{N}\colon\gambles(\vcnts{N})\to\gambles_{\permuts_N}(\values^N)
  \colon g\mapsto\vocntf{N}(g)\coloneqq g\circ\vbcntf{N},
\end{equation}
so $\vocntf{N}(g)$ is the permutation invariant gamble on $\values^N$ that assumes the constant value $g(\cnt{m})$ on the invariant atom~$\batom{m}$.
\par
Through the mediation of this identification $\vocntf{N}$, we can use the projection operator $\vex{N}$ to turn a gamble $f$ on $\values^N$ into a gamble on $\cnts^N$, as follows:
\begin{equation}
  \vmuhy{N}\colon\gambles(\values^N)\to\gambles(\vcnts{N})\colon
  f\mapsto\vmuhy{N}(f)\coloneqq\vmuhy{N}(f\vert\cdot),
\end{equation}
so $\vmuhy{N}(f)$ is the gamble on $\vcnts{N}$ that assumes the value $\vmuhy{N}(f\vert\cnt{m})$ in the count vector ${\cnt{m}\in\vcnts{N}}$.
By definition, $\vex{N}(f)=\vocntf{N}\bigl(\vmuhy{N}(f)\bigr)$ for all $f\in\gambles(\values^n)$, and similarly $\vmuhy{N}\bigl(\vocntf{N}(g)\bigr)=g$, for all $g\in\gambles(\vcnts{N})$.
Hence:
\begin{equation}\label{eq:muhy-identities}
  \vex{N}=\vocntf{N}\circ\vmuhy{N}
  \quad\text{ and }\quad
  \vmuhy{N}\circ\vocntf{N}=\identity{\gambles(\vcnts{N})}.
\end{equation}
Since $\vex{N}$ is a projection operator, its restriction to $\gambles_{\permuts_N}(\values^N)$ is the identity map, and therefore we infer from Eq.~\eqref{eq:muhy-identities} that the restriction of $\vmuhy{N}$ to $\gambles_{\permuts_N}(\values^N)$ is the inverse of $\vocntf{N}$, and therefore also a linear isomorphism between $\gambles_{\permuts_N}(\values^N)$ and $\gambles(\vcnts{N})$.
\par
If we invoke Eq.~\eqref{eq:ex-permutation} we find that
\begin{equation}\label{eq:muhy-permutation}
  \vmuhy{N}(\pi^tf)=\vmuhy{N}(f).
\end{equation}
Also taking into account the linearity of $\vmuhy{N}$ and Eq.~\eqref{eq:ex}, this leads to
\begin{equation}\label{eq:ex-muhy-permutation}
  \vmuhy{N}\bigl(\vex{N}(f)\bigr)=\vmuhy{N}(f).
\end{equation}
\par
\begin{figure}[htb]
  \centering
  \begin{tikzpicture}
    \matrix[matrix of math nodes,column sep=3em,row sep=10ex] {
      |(X)| \gambles(\values^N)\vphantom{\gambles_{\permuts_N}} && |(P)| \gambles_{\permuts_N}(\values^N)\\
      & |(N)| \gambles(\vcnts{N}) & \\
      & |(V)| \poly_N(\vsimplex) & \\
    };
    \path (X) edge[->] node[fill=white] {$\vex{N}$} (P);
    \path (X) edge[->] node[fill=white] {$\vmuhy{N}$} (N);
    \path (X) edge[->,out=-90,in=180] node[fill=white] {$\vmult{N}$} (V);
    \path (N) edge[->] ([yshift=-1pt]V.north);
    \path (N) edge[-,double,double distance=1pt] node[fill=white] {$\vcmult{N}$} (V);
    \path (P) edge[->,out=-90,in=0] ([xshift=-1pt]V.east);
    \path (P) edge[-,double,double distance=1pt,out=-90,in=0] node[fill=white] {$\vmult{N}$} (V);
    \path (N) edge[->] (P);
    \path (N) edge[-,double,double distance=1pt] node[fill=white] {$\vocntf{N}$} ($(N)!.84!(P)$);
  \end{tikzpicture}
  \caption{Single sequence length commutative diagram. Single arrows indicate linear monomorphisms (injective). Double arrows indicate linear isomorphisms (bijective).}
  \label{dia:same}
\end{figure}
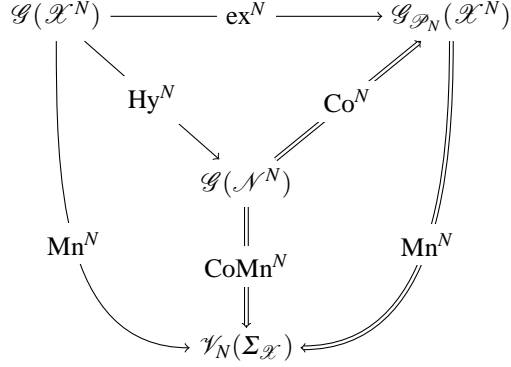

The relationships between the three important linear maps $\vex{N}$, $\vmuhy{N}$ and $\vocntf{N}$ we have introduced above are clarified by the commutative diagram in Fig.~\ref{dia:same}.
(The bottom part of the diagram can be safely ignored for now.)

\begin{example}
  In the context of our running example, we have the following:
  Take some gamble~$f$ on $\values^N$ and let $g\coloneqq\vmuhy{N}(f)$, then
  \begin{align*}
    g(2,0) = f(b,b), && g(0,2) = f(w,w), && g(1,1)=\tfrac{1}{2}\bigl(f(b,w)+f(w,b)\bigr).
  \end{align*}
  Conversely, take some gamble $g$ on~$\vcnts{N}$ and let $f\coloneqq\vocntf{N}(g)$, then
  \begin{align*}
    f(b,b)=g(2,0), && f(w,w)=g(0,2), && f(b,w)=f(w,b)=g(1,1). \mspace{20mu}\blacklozenge\mspace{-20mu}
  \end{align*}
\end{example}

For every gamble $f$ on $\values^N$, $f=\vex{N}(f)+[f-\vex{N}(f)]$, so it can be decomposed as a sum of a permutation invariant gamble $\vex{N}(f)$ and an element $f-\vex{N}(f)$ of the kernel $\weakly{\average_N}$ of the linear projection operator $\vex{N}$.
Elements of this kernel are, by definition, irrelevant as far as desirability under exchangeability is concerned, so the only part of this decomposition that matters is the element $\vex{N}(f)$ of $\gambles_{\permuts_N}(\values^N)$.
Since we have seen that $\vmuhy{N}$ acts as a linear isomorphism between the linear spaces $\gambles_{\permuts_N}(\values^N)$ and $\gambles(\vcnts{N})$, we now investigate whether we can use $\vmuhy{N}$ to represent a coherent and exchangeable $\desirs$ by some set of desirable \emph{count gambles} on $\vcnts{N}$.

\begin{theorem}[Finite Representation]\label{theo:finite-representation}
  A set of desirable gambles $\desirs$ on $\values^N$ is coherent and exchangeable iff there is some coherent set $\cdesirs$ of desirable gambles on $\vcnts{N}$ such that
  \begin{equation}\label{eq:representation-in-between}
    \desirs=(\vmuhy{N})^{-1}(\cdesirs),
  \end{equation}
  and in that case this $\cdesirs$ is uniquely determined by
  \begin{equation}
    \cdesirs=\set{g\in\gambles(\vcnts{N})}{\vocntf{N}(g)\in\desirs}=\vmuhy{N}(\desirs).
  \end{equation}
\end{theorem}

\noindent This leads to the following representation result for lower previsions, formulated without the mediation of coherent sets of desirable gambles.

\begin{corollary}\label{cor:finite-representation}
  A lower prevision $\lpr$ on $\gambles(\values^N)$ is coherent and exchangeable iff there is some coherent lower prevision $\clpr$ on $\gambles(\vcnts{N})$ such that $\lpr=\clpr\circ\vmuhy{N}$.
  In that case $\clpr$ is uniquely determined by $\clpr=\lpr\circ\vocntf{N}$.
\end{corollary}
\noindent
We call the set~$\cdesirs$ and the lower prevision~$\clpr$ the \emph{count representations} of the exchangeable set~$\desirs$ and the exchangeable lower prevision $\lpr$, respectively.
Our Finite Representation Theorem allows us to give an appealing geometrical interpretation to the notions of exchangeability and representation.
The exchangeability of~$\desirs$ means that it is completely determined by its count representation $\vmuhy{N}(\desirs)$, or what amounts to the same thing since~$\vocntf{N}$ is a linear isomorphism: by its projection $\vex{N}(\desirs)$ on the linear space of all permutation invariant gambles.
This turns count vectors into useful sufficient statistics (compare with Proposition~\ref{prop:sufficiency}), because the dimension of $\gambles(\vcnts{N})$ is typically much smaller than that of~$\gambles(\values^N)$.
To give an easy example: when $\values$ has two elements, $\gambles(\values^N)$ has dimension $2^N$, whereas the dimension of $\gambles(\vcnts{N})$ is only $N+1$.

\subsection{Exchangeable natural extension and representation}\label{sec:finite-representation-exnatex}
The exchangeable natural extension is easy to calculate using natural extension in terms of count representations, and the following simple result therefore has important consequences for practical implementations of reasoning and inference under exchangeability.

\begin{theorem}\label{theo:finite-representation-exnatex}
  Let $\assessment$ be a set of gambles on $\values^N$, then
  \begin{compactenum}[\upshape (i)]
    \item $\assessment$ avoids non-positivity under exchangeability iff\/ $\vmuhy{N}(\assessment)$ avoids non-positivity.
    \item $\vmuhy{N}\bigl(\exnatex{\assessment}\bigr)=\natex*{\vmuhy{N}(\assessment)}$.
  \end{compactenum}
\end{theorem}

\noindent This result gives us an extra approach to calculating the exchangeable natural extension of an assessment.
It reduces calculating the exchangeable natural extension to calculating a natural extension in the lower dimensional  space of count gambles.
The commutative diagram that corresponds to it is given in Fig.~\ref{dia:exnatex-natex}.

\begin{figure}[htb]
  \centering
  \begin{tikzpicture}
    \matrix[matrix of math nodes,column sep=8em,row sep=10ex] {
      |(Xa)| 2^{\gambles(\values^N)} & |(Na)| 2^{\gambles(\vcnts{N})}\\
      |(Xn)| 2^{\gambles(\values^N)} & |(Nn)| 2^{\gambles(\vcnts{N})}\\
    };
    \draw[->] (Xa) -- (Na) node[pos=.5,fill=white] {$\vmuhy{N}$};
    \draw[->] (Xn) -- (Nn) node[pos=.5,fill=white] {$\vmuhy{N}$};
    \draw[->] (Xa) -- (Xn) node[pos=.5,fill=white] {$\mathcal{E}_{\mathrm{ex}}^{N}$};
    \draw[->] (Na) -- (Nn) node[pos=.5,fill=white] {$\mathcal{E}$};
  \end{tikzpicture}
  \caption{The relationship between exchangeable natural extension and count representation natural extension. The arrows indicate monomorphisms (injective).}
  \label{dia:exnatex-natex}
\end{figure}

\subsection{Updating and representation}\label{sec:finite-representation-update}
Suppose, as in Section~\ref{sec:exchangeable-update}, that we update a coherent and exchangeable set of desirable gambles $\desirs$ after observing a sample $\osample{x}$ with count vector $\ocnt{m}$.
This leads to an updated coherent and exchangeable set of desirable gambles $\restrict{\desirs}{\osample{x}}=\restrict{\desirs}{\ocnt{m}}$ on $\values^{\rest{n}}$.
Here, we take a closer look at the corresponding set of desirable gambles on $\vcnts{\rest{n}}$, which we denote (symbolically) by $\restrict{\cdesirs}{\ocnt{m}}$.
(But we do not want to suggest with this notation that this is in some way an updated set of gambles!)
The Finite Representation Theorem~\ref{theo:finite-representation} tells us that
$\restrict{\cdesirs}{\ocnt{m}}=\vmuhy{\rest{n}}(\restrict{\desirs}{\ocnt{m}})$, but is there a direct way to infer the count representation $\restrict{\cdesirs}{\ocnt{m}}$ of $\restrict{\desirs}{\ocnt{m}}$ from the count representation $\cdesirs=\vmuhy{N}(\desirs)$ of $\desirs$?
\par
To show that there is, we need to introduce two new notions: the \emph{likelihood function}
  \begin{equation}\label{eq:likelihood}
    L_{\ocnt{m}}\colon\vcnts{\rest{n}}\to\reals
    \colon\rcnt{m}\mapsto L_{\ocnt{m}}(\rcnt{m})
    \coloneqq\frac{\abs{\oatom{m}}\,\abs{\ratom{m}}}{\abs{\atom{\ocnt{m}+\rcnt{m}}}},
  \end{equation}
  associated with sampling without replacement, and the linear map
  $+_{\ocnt{m}}$ from the linear space $\gambles(\vcnts{\rest{n}})$ to
  the linear space $\gambles(\vcnts{N})$ given by
  \begin{equation}
    +_{\ocnt{m}}\colon\gambles(\vcnts{\rest{n}})\to\gambles(\vcnts{N})
    \colon g\mapsto+_{\ocnt{m}}g
  \end{equation}
  where
  \begin{equation}\label{eq:count-reduction}
    +_{\ocnt{m}}g(\cnt{M})=
    \begin{cases}
      g(\cnt{M}-\ocnt{m})&\text{ if $\cnt{M}\geq\ocnt{m}$}\\
      0&\text{ otherwise}.
    \end{cases}
  \end{equation}

\begin{proposition}\label{prop:finite-representation-update}
  Consider a coherent and exchangeable set of desirable gambles~$\desirs$ on~$\values^N$, with count representation~$\cdesirs$.
  Let  $\restrict{\cdesirs}{\ocnt{m}}$ be the count representation of the coherent and exchangeable set of desirable gambles $\restrict{\desirs}{\ocnt{m}}$, obtained after updating~$\desirs$ with a sample~$\osample{x}$ with count vector~$\ocnt{m}$.
  Then
  \begin{equation}\label{eq:updating-counts}
    \restrict{\cdesirs}{\ocnt{m}}
    =\set{g\in\gambles(\vcnts{\rest{n}})}
    {+_{\ocnt{m}}(L_{\ocnt{m}}g)\in\cdesirs}.
  \end{equation}
\end{proposition}

\begin{example}
  In the context of our running example, where, recall, $\values\coloneqq\{b,w\}$ and $N\coloneqq2$, let $\rest{n}\coloneqq1$ and $\osample{x}\coloneqq w$ with count vector $\ocnt{m}=(1,0)$.
  Take $g$ to be some gamble on $\vcnts{\rest{n}}$ and let $g'\coloneqq+_{\ocnt{m}}(L_{\ocnt{m}}g)$, then
  \begin{align*}
    g'(2,0)=g(1,0), && g'(0,2)=0, && g'(1,1)=\tfrac{1}{2}g(0,1),
  \end{align*}
  So $g\in\restrict{\cdesirs}{\ocnt{m}}$ if $g'\in\cdesirs$.
  Contrast this with updating $\cdesirs$ with the information that one of the two observations is $w$, i.e., with conditioning event $\{\ocnt{m}\}+\vcnts{\rest{n}}=\{(2,0),(1,1)\}$.
  In that case $g\in\restrict{\cdesirs}{\bigl(\{\ocnt{m}\}+\vcnts{\rest{n}}\bigr)}$ if $g''\in\cdesirs$, where
  \begin{align*}
    g''(2,0)=g(1,0), && g''(0,2)=0, && g''(1,1)=g(0,1).~\mspace{20mu}\blacklozenge\mspace{-20mu}
  \end{align*}
\end{example}

\subsection{Multinomial processes}\label{sec:multproc}
Next, we turn to  a number of important ideas related to multinomial processes. They are at the same time useful for comparisons with the existing literature, and essential for our treatment of countable exchangeable sequences in Section~\ref{sec:infinite-xch}.
\par
Consider the \emph{$\values$-simplex}
\begin{equation}
  \vsimplex
  \coloneqq\biggset{\btheta\in\reals^\values}
  {\theta\geq0
    \text{ and }
  \smashoperator{\sum_{x\in\values}}\theta_x=1}.
\end{equation}
and, for $N\in\nats_0$, the linear map $\vcmult{N}$ from $\gambles(\vcnts{N})$ to $\gambles(\vsimplex)$ defined by
\begin{equation}\label{eq:cmultf}
  \vcmult{N}\colon\gambles(\vcnts{N})\to\gambles(\vsimplex)\colon
  g\mapsto\vcmult{N}(g)=\vcmult{N}(g\vert\cdot),
\end{equation}
where for all $\btheta\in\vsimplex$,
\begin{equation}\label{eq:cmult}
  \vcmult{N}(g\vert\btheta)
  \coloneqq\smashoperator{\sum_{\cnt{m}\in\vcnts{N}}}g(\cnt{m})\bern{m}(\btheta)
\end{equation}
is the expectation associated with the \emph{count multinomial} distribution with parameters~$N$ and~$\btheta$, and $\bern{m}$ is the  multivariate \emph{Bernstein (basis) polynomial of degree $N$} given by
\begin{equation}\label{eq:bernstein}
  \bern{m}(\btheta)
  \coloneqq\binom{N}{\cnt{m}}\smashoperator{\prod_{z\in\values}}\theta_z^{m_z}
  =\abs{\batom{m}}\smashoperator{\prod_{z\in\values}}\theta_z^{m_z}.
\end{equation}
$\vcmult{N}(\{m\}\vert\btheta)=\bern{m}(\btheta)$ is the probability of observing a count vector $\cnt{m}$ in a multinomial process where the possible outcomes $z\in\values$ have probability $\theta_z$.
\par
We also consider the related linear map $\vmult{N}$ from $\gambles(\values^N)$ to $\gambles(\vsimplex)$ defined by
\begin{equation}
  \vmult{N}\colon\gambles(\values^N)\to\gambles(\vsimplex)
  \colon f\mapsto\vmult{N}(f)=\vmult{N}(f\vert\cdot),
\end{equation}
where for all $\btheta\in\vsimplex$,
\begin{equation}\label{eq:mult}
  \vmult{N}(f\vert\btheta)
  \coloneqq\smashoperator{\sum_{\cnt{m}\in\vcnts{N}}}
  \vmuhy{N}(f\vert\cnt{m})\bern{m}(\btheta)
\end{equation}
is the expectation associated with the \emph{multinomial} distribution with parameters $N$ and $\btheta$.
We then have that
\begin{equation}\label{eq:mult-identities}
  \vcmult{N}=\vmult{N}\circ\vocntf{N}
  \quad\text{ and }\quad
  \vmult{N}=\vcmult{N}\circ\vmuhy{N}.
\end{equation}
If we consider a sequence of observations $\sample{x}$ with count vector $\cnt{m}$, then  $\vmult{N}(\{x\}\vert\btheta)=\bern{m}(\btheta)/\abs{\batom{m}}=\prod_{z\in\values}\theta_z^{m_z}$ is the probability of observing this sequence in a multinomial process where the possible outcomes $z\in\values$ have probability $\theta_z$.
\par
The Bernstein basis polynomials~$\bern{m}$, $\cnt{m}\in\vcnts{N}$ form a basis for the linear space $\vpoly[N]$ of all polynomials on $\vsimplex$ of degree up to $N$.
This means that for each polynomial $p$ whose degree $\deg(p)$ does not exceed $N$, there is a unique gamble $\bexp{p}{N}$ on $\vcnts{N}$ such that $p=\vcmult{N}(\bexp{p}{N})$.
We denote by $\vpoly$ the linear space of all polynomials on $\vsimplex$.
More details on Bernstein basis polynomials can be found in Appendix~\ref{app:bernstein}.

\begin{example}\label{ex:bernstein}
  For our running example, the unit simplex $\vsimplex=\simplex{\{b,w\}}$ is a line of unit length; frequency vectors $(\theta_b,\theta_w)$ can be parametrised by $\theta_b\in[0,1]$, as $\theta_w=1-\theta_b$.
  One of the line's extreme points corresponds to the frequency vector $(1,0)$, `$b$', the other to the frequency vector $(0,1)$, `$w$'.
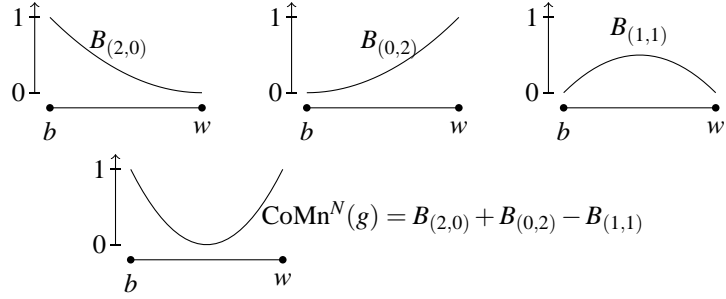
  \begin{center}
    \begin{tikzpicture}[xscale=2]
      \draw (0,-.2) node[fill,circle,inner sep=1pt,label=-90:$b$] {} -- (1,-.2) node[fill,circle,inner sep=1pt,label=-90:$w$] {};
      \draw[->] (-.1,0) -- (-.1,1.2);
      \path (-.1,0) \ytick node[left] {$0$} (-.1,1) \ytick node[left] {$1$};
      \draw plot[domain=0:1] (1-\x,\x*\x);
      \node[right] at (.2,.64) {$\bern{(2,0)}$};
    \end{tikzpicture}
    \hspace{1em}
    \begin{tikzpicture}[xscale=2]
      \draw (0,-.2) node[fill,circle,inner sep=1pt,label=-90:$b$] {} -- (1,-.2) node[fill,circle,inner sep=1pt,label=-90:$w$] {};
      \draw[->] (-.1,0) -- (-.1,1.2);
      \path (-.1,0) \ytick node[left] {$0$} (-.1,1) \ytick node[left] {$1$};
      \draw plot[domain=0:1] (1-\x,{(1-\x)*(1-\x)});
      \node[left] at (.8,.64) {$\bern{(0,2)}$};
    \end{tikzpicture}
    \hspace{1em}
    \begin{tikzpicture}[xscale=2]
      \draw (0,-.2) node[fill,circle,inner sep=1pt,label=-90:$b$] {} -- (1,-.2) node[fill,circle,inner sep=1pt,label=-90:$w$] {};
      \draw[->] (-.1,0) -- (-.1,1.2);
      \path (-.1,0) \ytick node[left] {$0$} (-.1,1) \ytick node[left] {$1$};
      \draw plot[domain=0:1] (1-\x,{2*\x*(1-\x)});
      \node[above] at (.5,.5) {$\bern{(1,1)}$};
    \end{tikzpicture}
    \begin{tikzpicture}[xscale=2]
      \draw (0,-.2) node[fill,circle,inner sep=1pt,label=-90:$b$] {} -- (1,-.2) node[fill,circle,inner sep=1pt,label=-90:$w$] {};
      \draw[->] (-.1,0) -- (-.1,1.2);
      \path (-.1,0) \ytick node[left] {$0$} (-.1,1) \ytick node[left] {$1$};
      \draw plot[domain=0:1] (1-\x,{\x*\x + (1-\x)*(1-\x) - 2*\x*(1-\x)});
      \node[right] at (.8,.36) {$\vcmult{N}(g)=\bern{(2,0)}+\bern{(0,2)}-\bern{(1,1)}$};
    \end{tikzpicture}
  \end{center}
  The Bernstein basis polynomials for the case $N=2$ are given on the top row; on the bottom row, we give the polynomial $\vcmult{N}(g)$ corresponding to the gamble $g$ on $\vcnts{N}$ defined by $g(2,0)=g(0,2)=-g(1,1)=1$.~$\blacklozenge$
\end{example}

\subsection{Finite representation in terms of polynomials}\label{sec:finite-polynomials}
We see that the range of the linear maps $\vcmult{N}$ and $\vmult{N}$ is the linear space $\vpoly[N]$.
Moreover, since for every polynomial $p$ of degree up to~$N$, i.e., for every $p\in\vpoly[N]$, there is a \emph{unique} count gamble $\smash[b]{\bexp{p}{N}}\in\gambles(\vcnts{N})$ such that $p=\vcmult{N}(\smash[b]{\bexp{p}{N}})$, $\vcmult{N}$ is a \emph{linear isomorphism} between the linear spaces $\gambles(\vcnts{N})$ and $\vpoly[N]$.
The relationships between the five important linear maps we have introduced so far are clarified by the commutative diagram in Fig.~\ref{dia:same}.
\par
In summary, everything that can be expressed using the language of gambles on $\vcnts{N}$, can also be expressed using the language of polynomial gambles on $\vsimplex$ of degree up to $N$, and {\itshape vice versa}.
Again, as explained above, the fundamental reason why this is possible, is that the Bernstein basis polynomials of degree $N$ constitute a basis for the linear space of all polynomials of degree up to $N$, where a count gamble $g$ plays the r\^ole of a coordinate representation for a polynomial $p$ in this basis.
The map $\vcmult{N}$ and its inverse are the tools that take care of the translation between the two languages.
This is essentially what is behind the representation theorem for countable exchangeable sequences that we will turn to in Section~\ref{sec:infinite-xch}.
\par
In order to lay the proper foundations for this work, we now prove a version of the finite representation theorem in terms of polynomial gambles of degree $N$ on $\vsimplex$, rather than count gambles on $\vcnts{N}$.

\begin{definition}[Bernstein coherence]\label{def:Bern-coh-fixed-degree}
  We call a set $\tdesirs$ of polynomials in $\vpoly[N]$ \emph{Bernstein coherent at degree $N$} if it satisfies the following properties: for all $p,p_1,p_2\in\vpoly[N]$ and all real $\lambda>0$,
  \begin{compactenum}[\upshape {B}$_N$1.]
  \item if $p=0$ then $p\notin\tdesirs$;\label{item:BNnonzero}
  \item if $p$ is such that $\bexp{p}{N}>0$ then $p\in\tdesirs$;\label{item:BNapg}
  \item if $p\in\tdesirs$ then $\lambda p\in\tdesirs$;\label{item:BNscl}
  \item if $p_1,p_2\in\tdesirs$ then $p_1+p_2\in\tdesirs$.\label{item:BNcmb}
  \end{compactenum}
\end{definition}

\noindent Bernstein coherence at degree $N$ is very closely related to coherence, the only difference being that we do not consider whether a polynomial $p$ is positive, but whether its Bernstein expansion $\bexp{p}{N}$ is.
This means that models in terms of sequences or count vectors are authoritative over those in terms of frequency vectors in the sense that polynomials are not directly behaviourally interpreted as gambles.
This is related to the fact that not all possible frequency vectors can practically be observed.

\begin{example}
  Any polynomial with a positive expansion in terms of Bernstein basis polynomials is positive.
  But the pair $\bigl(g,\vcmult{N}(g)\bigr)$ of Example~\ref{ex:bernstein} shows that a polynomial can be positive, while its Bernstein expansion is not.
  So the smallest set of polynomials Bernstein coherent at degree~$2$ is $\posi\bigl(\{\bern{(2,0)},\bern{(0,2)},\bern{(1,1)\}}\bigr)$.
  Taking $\vcmult{N}(g)$ to be desirable corresponds to the assessment that observing differing colors is less likely than observing identical ones.~$\blacklozenge$
\end{example}

\noindent Bernstein coherence at degree $N$ is a special case of the general concept of coherence relative to $\spacecone$, discussed in Section~\ref{sec:gendesir}, where $\subspace\coloneqq\vpoly[N]$ and $\cone$ is the convex cone of all polynomials of degree at most $N$ with a non-negative expansion $\bexp{p}{N}$ in the Bernstein basis of degree $N$:
\begin{equation}
  \cone\coloneqq\set{p\in\vpoly[N]}{\bexp{p}{N}\geq0}.
\end{equation}

\begin{theorem}[Finite Representation]\label{theo:finite-representation-polynomials}
  A set of desirable gambles $\desirs$ on $\values^N$, with count representation $\cdesirs\coloneqq\vmuhy{N}(\desirs)$, is coherent and exchangeable iff there is some subset\/ $\tdesirs$ of\/ $\vpoly[N]$,  Bernstein coherent at degree $N$, such that
  \begin{equation}
    \desirs=(\vmult{N})^{-1}(\tdesirs)
    \quad\text{ or equivalently }\quad
    \cdesirs=(\vcmult{N})^{-1}(\tdesirs),
  \end{equation}
  and in that case this $\tdesirs$ is uniquely determined by
  \begin{equation}
    \tdesirs
    =\vmult{N}(\desirs)
    =\vcmult{N}(\cdesirs).
  \end{equation}
\end{theorem}
\noindent
We call the set $\tdesirs =\vmult{N}(\desirs)$ the \emph{frequency representation} of the coherent and exchangeable set of desirable gambles $\desirs$.

\section{Countable exchangeable sequences}\label{sec:infinite-xch}
With the experience gained in investigating finite exchangeable sequences, we are now ready to tackle the problem of reasoning about countably infinite exchangeable sequences, our second main topic.
The first step is to use the finite frequency representation results of Sections~\ref{sec:multproc} and~\ref{sec:finite-polynomials} to find a Representation Theorem for infinite exchangeable sequences  (Section~\ref{sec:infinite-representation}).
We can then show what updating and natural extension look like in terms of this frequency representation (respectively Section~\ref{sec:infrep-updating} and Sections~\ref{sec:bernstein-natural-extension} and~\ref{sec:exchangeable-natural-extension-infinite}).

\subsection{Infinite representation}\label{sec:infinite-representation}
We consider a countable sequence $\rv_1$, \dots, $\rv_N$, \dots\ of random variables assuming values in the same finite set~$\values$.
We call this sequence \emph{exchangeable} if each of its finite subsequences is, or equivalently, if for all $n\in\nats_0$, the random variables $\rv_1$,~\dots,~$\rv_n$ are exchangeable.
\par
How can we model this? First of all, this means that for each $n\in\nats_0$, there is a coherent and exchangeable set of desirable gambles $\vdesirs{n}$ on $\values^n$.
Equivalently, we have a coherent set of desirable gambles (count representation) $\vcdesirs{n}\coloneqq\vmuhy{n}(\vdesirs{n})$ on $\vcnts{n}$, or a set (frequency representation) $\vtdesirs[n]\coloneqq\vmult{n}(\vdesirs{n})=\vcmult{n}(\vcdesirs{n})$ of polynomials in $\vpoly[n]$, Bernstein coherent at degree $n$.
\par
In addition, there is a time-consistency constraint.
Consider the following linear \emph{projection operators}, with $n_1\leq n_2$:
\begin{equation}\label{eq:projection}
  \vproj{n_2}{n_1}\colon\values^{n_2}\to\values^{n_1}\colon
  \tuple{x}{n_2}\mapsto
  \vproj{n_2}{n_1}\tuple{x}{n_2}\coloneqq\tuple{x}{n_1}.
\end{equation}
With each such operator there corresponds a linear map $\vexten{n_1}{n_2}$ between the linear spaces $\gambles(\values^{n_1})$ and $\gambles(\values^{n_2})$, defined as follows:
\begin{equation}\label{eq:extension}
  \vexten{n_1}{n_2}\colon\gambles(\values^{n_1})\to\gambles(\values^{n_2})\colon
  f\mapsto\vexten{n_1}{n_2}(f)=f\circ\vproj{n_2}{n_1}.
\end{equation}
In other words, $\vexten{n_1}{n_2}(f)$ is the \emph{cylindrical extension} of the gamble $f$ on $\values^{n_1}$ to a gamble on $\values^{n_2}$.
\par
\emph{Time-consistency} now means that if we consider a gamble on $\values^{n_2}$ that really only depends on the first $n_1$ variables, it should not matter, as far as its desirability is concerned, whether we consider it to be a gamble on $\values^{n_1}$ or a gamble on $\values^{n_2}$.
More formally:
\begin{equation}\label{eq:time-consistency}
  (\forall n_1\leq n_2)
  \vexten{n_1}{n_2}(\vdesirs{n_1})
  =\vdesirs{n_2}\cap\vexten{n_1}{n_2}\bigl(\gambles(\values^{n_1})\bigr).
\end{equation}
How can we translate this constraint in terms of the count representations $\vcdesirs{n}$ or the frequency representations $\vtdesirs[n]$? Using the Finite Representation Theorem~\ref{theo:finite-representation}, we see that
$f\in\vdesirs{n_k}\Leftrightarrow\vmuhy{n_k}(f)\in\vcdesirs{n_k}$.
It follows from a few algebraic manipulations that for any gamble $f$ on $\values^{n_1}$ and all $\cnt{M}\in\vcnts{n_2}$:
\begin{equation}
  \vmuhy{n_2}(\vexten{n_1}{n_2}(f)\vert\cnt{M})
  =\smashoperator{\sum_{\cnt{m}\in\vcnts{n_1}}}
  \frac{\abs{\batom{M-m}}\,\abs{\batom{m}}}{\abs{\batom{M}}}
  \vmuhy{n_1}(f\vert\cnt{m}).
\end{equation}
So if we introduce the linear \emph{extension map} $\venl{n_1}{n_2}$ from the linear space $\gambles(\vcnts{n_1})$ to the linear space $\gambles(\vcnts{n_2})$ as follows:
\begin{equation}
  \venl{n_1}{n_2}\colon\gambles(\vcnts{n_1})\to\gambles(\vcnts{n_2})\colon
  g\mapsto\venl{n_1}{n_2}(g)\coloneqq\smashoperator{\sum_{\cnt{m}\in\vcnts{n_1}}}
  \frac{\abs{\batom{\cdot-m}}\,\abs{\batom{m}}}{\abs{\batom{\cdot}}}
  g(\cnt{m}),
\end{equation}
this can be summarised succinctly as:
\begin{equation}\label{eq:enl-identity}
  \vmuhy{n_2}\circ\vexten{n_1}{n_2}
  =\venl{n_1}{n_2}\circ\vmuhy{n_1},
\end{equation}
and we see that the time-consistency requirement~\eqref{eq:time-consistency} is then equivalent to [see Appendix~\ref{app:proofs} for a detailed proof]:
\begin{equation}\label{eq:time-consistency-counts}
  (\forall n_1\leq n_2)
  \venl{n_1}{n_2}(\vcdesirs{n_1})
  =\vcdesirs{n_2}\cap\venl{n_1}{n_2}\bigl(\gambles(\vcnts{n_1})\bigr),
\end{equation}
which is in turn equivalent to [see Appendix~\ref{app:proofs} for a detailed proof]:
\begin{equation}\label{eq:time-consistency-frequencies}
  (\forall n_1\leq n_2)
  \vtdesirs[n_1]
  =\vtdesirs[n_2]\cap\vpoly[n_1].
\end{equation}
We see that the time consistency condition can be most elegantly expressed in terms of the frequency representations.
\begin{example}
  Let us illustrate the newly introduced operators in the context of our running example.
  Take $n_1\coloneqq1$, $n_1\coloneqq2$, and $f$ a gamble on $\values^{n_1}$; let $f'\coloneqq\vexten{n_1}{n_2}(f)$, then
  \begin{equation*}
    f'(b,b) = f'(b,w) = f(b) \quad\text{ and }\quad f'(w,w) = f'(w,b) = f(w).
  \end{equation*}
  Now take a gamble $g$ on $\vcnts{n_1}$; and let $g'\coloneqq\venl{n_1}{n_2}(g)$, then
  \begin{align*}
    g'(2,0) = g(1,0), && g'(0,2) = g(0,1), && g'(1,1) = \tfrac{1}{2}\bigl(g(1,0)+g(0,1)\bigr). \mspace{20mu}\blacklozenge\mspace{-20mu}
  \end{align*}
\end{example}
\par
We call the family $\vdesirs{n}$, $n\in\nats_0$ \emph{time-consistent, coherent and exchangeable} when each member $\vdesirs{n}$ is coherent and exchangeable, and when the family $\vdesirs{n}$, $n\in\nats_0$ satisfies Eq.~\eqref{eq:time-consistency}.
\par
The (count) multinomial expectations introduced in the previous section also satisfy a nice time consistency property.
If we consider a gamble $f_1$ on $\values^{n_1}$, then we can also consider it as a gamble $\vexten{n_1}{n_2}(f_1)$ on $\values^{n_2}$, and of course both versions of this gamble should have the same multinomial expectation.
This leads to the following identities:
\begin{equation}\label{eq:mult-extensions}
  \vmult{n_2}\circ\vexten{n_1}{n_2}=\vmult{n_1}
  \quad\text{ and }\quad
  \vcmult{n_2}\circ\venl{n_1}{n_2}=\vcmult{n_1},
\end{equation}
where the second identity follows from combining the first with Eqs.~\eqref{eq:enl-identity} and~\eqref{eq:mult-identities}.
\par
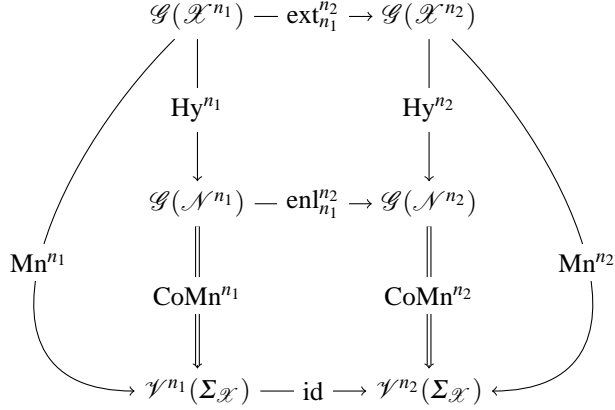
\begin{figure}[htb]
  \centering
  \begin{tikzpicture}
    \matrix[matrix of math nodes,column sep=2em,row sep=12ex] {
      |(X1)| \gambles(\values^{n_1}) && |(X2)| \gambles(\values^{n_2})\\
      |(N1)| \gambles(\vcnts{n_1}) && |(N2)| \gambles(\vcnts{n_2})\\
      |(S1)| \vpoly[n_1] && |(S2)| \vpoly[n_2] \\
    };
    \path[->] (X1) edge node[fill=white] {$\vexten{n_1}{n_2}$} (X2);
    \path[->] (N1) edge node[fill=white] {$\venl{n_1}{n_2}$} (N2);
    \path[->] (S1) edge node[fill=white] {$\iden$} (S2);
    \path[->] (X1) edge node[fill=white] {$\vmuhy{n_1}$} (N1);
    \path[->] (X2) edge node[fill=white] {$\vmuhy{n_2}$} (N2);
    \path (N1) edge[->] (S1);
    \path (N1) edge[-,double,double distance=1pt] node[fill=white] {$\vcmult{n_1}$} ($(N1)!.98!(S1.north)$);
    \path (N2) edge[->] (S2);
    \path (N2) edge[-,double,double distance=1pt] node[fill=white] {$\vcmult{n_2}$} ($(N2)!.98!(S2.north)$);
    \path[->] (X1) edge[out=-135,in=180,min distance=15ex] node[fill=white] {$\vmult{n_1}$} (S1);
    \path[->] (X2) edge[out=-45,in=0,min distance=15ex] node[fill=white] {$\vmult{n_2}$} (S2);
  \end{tikzpicture}
  \caption{Different sequence length commutative diagram. Single arrows indicate linear monomorphisms (injective). Double arrows indicate linear isomorphisms (bijective).}
  \label{dia:different}
\end{figure}

The relationships between three of the linear maps we encountered earlier and the maps related to time-consistency introduced here are clarified by the commutative diagrams in Fig.~\ref{dia:different}.

We can generalise the concept of Bernstein coherence given in Definition~\ref{def:Bern-coh-fixed-degree} to sets of polynomials of arbitrary degree:
\begin{definition}[Bernstein coherence]\label{def:Bern-coh}
  We call a set $\tdesirs$ of polynomials in $\vpoly$ \emph{Bernstein coherent} if it satisfies the following properties: for all $p,p_1,p_2\in\vpoly$ and all real $\lambda>0$,
  \begin{compactenum}[\upshape {B}1.]
  \item if $p=0$ then $p\notin\tdesirs$;\label{item:Bnonzero}
  \item if $p$ is such that $\bexp{p}{n}>0$ for some $n\geq\deg(p)$, then $p\in\tdesirs$;\label{item:Bapg}
  \item if $p\in\tdesirs$ then $\lambda p\in\tdesirs$;\label{item:Bscl}
  \item if $p_1,p_2\in\tdesirs$ then $p_1+p_2\in\tdesirs$.\label{item:Bcmb}
  \end{compactenum}
\end{definition}

It is clear that we can replace B\ref{item:Bnonzero} by the following requirement, because it is equivalent to it under B\ref{item:Bapg}--B\ref{item:Bcmb} [see Appendix~\ref{app:proofs} for a proof]:
\begin{compactenum}[\upshape {B}1.]\addtocounter{enumi}{4}
  \item If $p$ is such that $\bexp{p}{n}\leq0$ for some $n\geq\deg(p)$, then $p\notin\tdesirs$.
\end{compactenum}

This type of Bernstein coherence is again very closely related to coherence, the only difference being that not all positive polynomials, but rather all polynomials with some positive Bernstein expansion are required to belong to a Bernstein coherent set.

\begin{example}
  The parabola $\vcmult{N}(g)$ of Example~\ref{ex:bernstein} also shows that a polynomial can be positive, while no Bernstein expansion of any order is.
  This follows from the fact that all Bernstein basis polynomials are strictly positive on the interior of the unit simplex and that this parabola has a minimum of\/ $0$ within this interior.~$\blacklozenge$
\end{example}

\noindent Bernstein coherence is a special case of the general concept of coherence relative to $\spacecone$, discussed in Section~\ref{sec:gendesir}, where $\subspace\coloneqq\vpoly$ and $\cone$ is the convex cone of all polynomials with some non-negative Bernstein expansion:
\begin{equation}\label{eq:vnnegpoly}
  \cone
  \coloneqq\vnnegpoly
  \coloneqq\set{p\in\vpoly}{(\exists n\geq0)\bexp{p}{n}\geq0}.
\end{equation}
We also denote the set $\alldesirs[(\vpoly,\vnnegpoly)](\vsimplex)$ of all Bernstein coherent subsets of $\vpoly$ by $\allberndesirs(\vsimplex)$.

We are now ready to formulate our Infinite Representation Theorem~\ref{theo:infinite-representation}, which is a significant generalisation of \citeauthor{finetti1937}'s representation result for countable sequences \citeyearpar{finetti1937}.
A similar result can also be proved for coherent lower previsions \cite{cooman2006d}.

\begin{theorem}[Infinite Representation]\label{theo:infinite-representation}
  A family $\vdesirs{n}$, $n\in\nats_0$ of sets of desirable gambles on~$\values^n$, with associated count representations $\vcdesirs{n}\coloneqq\vmuhy{n}(\vdesirs{n})$ and frequency representations $\vtdesirs[n]\coloneqq\vmult{n}(\vdesirs{n})=\vcmult{n}(\vcdesirs{n})$, is time-consistent, coherent and exchangeable iff there is some Bernstein coherent set $\vtdesirs$ of polynomials in $\vpoly$ such that, for all $n\in\nats_0$,
  \begin{equation}\label{eq:infinite-representation-backtooriginal}
    \vcdesirs{n}=(\vcmult{n})^{-1}(\vtdesirs) \quad\text{ and }\quad \vdesirs{n}=(\vmult{n})^{-1}(\vtdesirs),
  \end{equation}
  and in that case this $\vtdesirs$ is uniquely given by
  \begin{equation}
    \vtdesirs=\bigcup_{n\in\nats_0}\vtdesirs[n].
  \end{equation}
\end{theorem}
\noindent
We call $\vtdesirs$ the \emph{frequency representation} of the coherent, exchangeable and time-consistent family of sets of desirable gambles $\vdesirs{n}$, $n\in\nats_0$.

\subsection{Updating and infinite representation}\label{sec:infrep-updating}
Suppose we have a coherent, exchangeable and time-consistent family of sets of desirable gambles $\vdesirs{n}$, $n\in\nats_0$, with associated count representations $\vcdesirs{n}\coloneqq\vmuhy{n}(\vdesirs{n})$ and associated frequency representation $\vtdesirs\coloneqq\bigcup_{n\in\nats}\vtdesirs[n]$ with $\vtdesirs[n]\coloneqq\vmult{n}(\vdesirs{n})$.
\par
Now suppose we observe the values $\osample{x}$ of the first $\obs{n}$ variables, with associated count vector $\ocnt{m}\coloneqq\vbcntf{\obs{n}}(\osample{x})$, then we have seen in Section~\ref{sec:finite-representation-update} that these models $\vdesirs{n}$ and $\vcdesirs{n}$ (for $n>\obs{n}$) get updated to coherent and exchangeable models
$\vrestrict{n}{m}$ with count representations $\vcrestrict{n}{m}=\vmuhy{\rest{n}}(\vrestrict{n}{m})$ for $\rest{n}\coloneqq n-\obs{n}=1,2,\dots$.
It turns out that updating becomes especially easy in terms of the frequency representation.

\begin{theorem}\label{theo:updated-simplex-representation}
  Consider a coherent, exchangeable and time-consistent family of sets of desirable gambles $\vdesirs{n}$, $n\in\nats_0$, with associated frequency representation $\vtdesirs$.
  After updating with a sample with count vector $\ocnt{m}\in\cnts_\values^{\obs{n}}$, the family $\vrestrict{n}{m}$, $\rest{n}\in\nats_0$ is still coherent, exchangeable and time-consistent, and has frequency representation
  \begin{equation}\label{eq:updated-simplex-representation}
    \vtrestrict{m}\coloneqq\set{p\in\vpoly}{\bern{\ocnt{m}}p\in\vtdesirs}.
  \end{equation}
\end{theorem}

\subsection{Independence: iid sequences}
Theorem~\ref{theo:updated-simplex-representation} can be used to find an easy and quite intriguing characterisation of a sequence of \emph{independent and identically distributed} (iid) random variables $\rv_1$, \dots, $\rv_N$, \dots\ assuming values in a finite set $\values$. This is an exchangeable sequence where learning the value of any finite number of variables does not change our beliefs about the remaining, unobserved ones. We infer from Theorem~\ref{theo:updated-simplex-representation} that such will be the case iff the frequency representation $\vtdesirs$ of the sequence satisfies
\begin{equation}\label{eq:iid-one}
  (\forall\obs{n}\in\nats_0)
  (\forall\ocnt{m}\in\cnts_\values^{\obs{n}})
  \vtrestrict{m}=\vtdesirs,
\end{equation}
which is equivalent to
\begin{equation}\label{eq:iid-two}
  (\forall\obs{n}\in\nats_0)
  (\forall\ocnt{m}\in\cnts_\values^{\obs{n}})
  \bigl(\forall p\in\vpoly\bigr)
 (p\in\vtdesirs\asa\bern{\ocnt{m}}p\in\vtdesirs).
\end{equation}
Any Bernstein coherent set of polynomials that satisfies one of the equivalent conditions~\eqref{eq:iid-one} or~\eqref{eq:iid-two} is an imprecise-probabilistic model for a (discrete-time) iid-process, or equivalently, a multinomial process, assuming values in a set $\values$.
\par
Let us define  $e_z$ as the special count vector corresponding to a single observation of $z\in\values$: the $z$-component of $e_z$ is one, and all other components are zero.
Observe that $\bern{e_z}(\theta)=\theta_z$.
The precise-probabilistic iid-processes, or in other words, the multinomial processes, correspond to the \emph{maximal} coherent sets of polynomials that satisfy the iid condition:

\begin{proposition}\label{prop:maximal-iid}
  Consider any maximal element $\vtdesirs$ of\/ $\allberndesirs(\vsimplex)$ that satisfies either of the equivalent conditions~\eqref{eq:iid-one} or~\eqref{eq:iid-two}.
  Let $\lpr_{\vtdesirs}$ be the lower prevision defined on $\vpoly$ in the usual way by letting $\lpr_{\vtdesirs}(p)\coloneqq\sup\set{\alpha}{p-\alpha\in\vtdesirs}$ for all $p\in\vpoly$.
  Then $\lpr_{\vtdesirs}$ is a linear functional that dominates the $\min$ functional, and is completely determined by $\lpr_{\vtdesirs}(p)=p(\vartheta)$ for all $p\in\vpoly$, where $\vartheta_z\coloneqq\lpr_{\vtdesirs}(\bern{e_z})$ for all $z\in\values$.
  In addition, consider $n\in\nats$ and let $\vcdesirs{n}\coloneqq(\vcmult{n})^{-1}(\vtdesirs)$ be the corresponding set of desirable count gambles on $\vcnts{n}$, with associated lower prevision $\lpr_{\vcdesirs{n}}$.
  Then $\lpr_{\vcdesirs{n}}(g)=\lpr_{\vtdesirs}(\vcmult{n}(g))=\vcmult{n}(g\vert\vartheta)=\sum_{\cnt{m}\in\vcnts{n}}g(\cnt{m})\bern{m}(\vartheta)$ for all gambles $g$ on $\vcnts{n}$, and in particular the probability of observing a count vector $\cnt{m}\in\vcnts{n}$ is given by $\lpr_{\vcdesirs{n}}(\{\cnt{m}\})=\bern{m}(\vartheta)$.
 \end{proposition}
\noindent While it appears that such imprecise iid-processes are interesting, much more work needs to be done before we can get a complete picture of their structural properties and practical relevance.
We leave this as a potential avenue for further research.

\subsection{Bernstein natural extension}\label{sec:bernstein-natural-extension}
The intersection of an arbitrary non-empty family of Bernstein coherent sets of polynomials is still Bernstein coherent.
This is the idea behind the following theorem, which is a special instance of Theorem~\ref{theo:natex} with $\subspace\coloneqq\vpoly$ and $\cone\coloneqq\vnnegpoly$.
\par
We denote by $\vpospoly$ the set of all polynomials on $\vsimplex$ with some positive Bernstein expansion:
\begin{equation}
  \vpospoly=\set{p\in\vpoly}{\bigl(\exists n\geq\deg(p)\bigr)\bexp{p}{n}>0}.
\end{equation}
and by $\vnpospoly$ the set of all polynomials on $\vsimplex$ with some non-positive Bernstein expansion:
\begin{equation}
  \vnpospoly=\set{p\in\vpoly}{\bigl(\exists n\geq\deg(p)\bigr)\bexp{p}{n}\leq0}.
\end{equation}
Moreover, we say that a set $\passessment$ of polynomials \emph{avoids Bernstein non-positivity} if no polynomial in its positive hull $\posi(\passessment)$ has any non-positive Bernstein expansion, i.e.
\begin{equation}
  \posi(\passessment)\cap\vnpospoly=\emptyset;
\end{equation}
clearly, this is the case iff $\passessment$ avoids non-positivity relative to $\bigl(\vpoly,\vnnegpoly\bigr)$.
We also call  the $\bigl(\vpoly,\vnnegpoly\bigr)$-natural extension $\natex[(\vpoly,\vnnegpoly)]{\passessment}$ of $\passessment$ its Bernstein natural extension, and denote it by $\bernnatex{\passessment}$.

\begin{theorem}[Bernstein natural extension]\label{theo:bernstein-natex}
  Consider a set of polynomials $\passessment\subseteq\vpoly$, and define its \emph{Bernstein natural extension}
  \begin{equation}\label{eq:bernstein-natex}
    \bernnatex{\passessment}
    \coloneqq\bigcap\set{\vtdesirs\in\allberndesirs(\vsimplex)}
    {\passessment\subseteq\vtdesirs}.
  \end{equation}
  Then the following statements are equivalent:
  \begin{compactenum}[\upshape (i)]
  \item $\passessment$ avoids Bernstein non-positivity;\label{item:natextBapl}
  \item $\passessment$ is included in some Bernstein coherent set of polynomials;\label{item:natextBsubsetcoh}
  \item $\bernnatex{\passessment}\neq\vpoly$;\label{item:natextBnonempty}
  \item $\bernnatex{\passessment}$ is a Bernstein coherent set of polynomials;\label{item:natextBiscoh}
  \item $\bernnatex{\passessment}$ is the smallest Bernstein coherent set of polynomials that includes $\passessment$.\label{item:natextBsmallestcoh}
  \end{compactenum}
  When any (and hence all) of these equivalent statements hold, then
  \begin{equation}\label{eq:posi-bernstein-natex}
    \bernnatex{\passessment}
    =\posi\bigl(\vpospoly\cup\assessment\bigr).
  \end{equation}
\end{theorem}

\begin{example}
  Recall that, for our running example, the unit simplex $\vsimplex=\simplex{\{b,w\}}$ is a line of unit length parametrised by $\theta_b\in[0,1]$, with $\theta_w=1-\theta_b$.
  Consider the polynomial $p\in\vpoly$ such that $p(\theta_b)=-1+3\theta_b-\theta_b^2$.
  Because $p(1)=1$, we have that $\max p>0$, and thus via Eq.~\eqref{eq:bernstein-coefficients-4} that $\max\bexp{p}{n}>0$ for all $n\geq2$.
  So the assessment $\{p\}$ avoids Bernstein non-positivity, because its Bernstein expansions of degree 2 and up are not non-positive.
\end{example}

\subsection{Exchangeable natural extension for infinite sequences}\label{sec:exchangeable-natural-extension-infinite}
To finish this discussion of exchangeability for infinite sequences of random variables, we take up the issue of inference, and extend the notion of exchangeable natural extension, discussed in Section~\ref{sec:exchangeable-natural-extension}, from finite to infinite sequences.
\par
This extension is fairly straightforward.
Suppose that for each~$i$ in the non-empty index set~$I$, we consider a coherent, exchangeable and time-consistent family $\fdesirs{n}{i}$, $n\in\nats_0$ of sets of desirable gambles.
As we know from our Infinite Representation Theorem~\ref{theo:infinite-representation}, each such family is represented by a Bernstein coherent set of polynomials on $\vsimplex$:
\begin{equation}
  \ftdesirs{i}=\bigcup_{n\in\nats_0}\vmult{n}(\fdesirs{n}{i})
\end{equation}
in the sense that, for all $n\in\nats_0$,
\begin{equation}
  \fdesirs{n}{i}=(\vmult{n})^{-1}(\ftdesirs{i}).
\end{equation}
We know from the previous section that the intersection of a non-empty family of Bernstein coherent sets of polynomials is still Bernstein coherent.
This implies that $\vtdesirs\coloneqq\bigcap_{i\in I}\ftdesirs{i}$ is a Bernstein coherent set of polynomials such that,  for all $n\in\nats_0$,
\begin{equation}
  \vdesirs{n}
  =\bigcap_{i\in I}\fdesirs{n}{i}
  =\bigcap_{i\in I}(\vmult{n})^{-1}(\ftdesirs{i})
  =(\vmult{n})^{-1}\Bigl(\bigcap_{i\in I}\ftdesirs{i}\Bigr)
  =(\vmult{n})^{-1}(\vtdesirs),
\end{equation}
implying that \emph{the (element-wise) intersection $\vdesirs{n}$, $n\in\nats_0$ of the coherent, exchangeable and time-consistent families $\fdesirs{n}{i}$, $n\in\nats_0$ is still a coherent, exchangeable and time-consistent family, whose frequency representation $\vtdesirs$ is the intersection of the frequency representations~$\ftdesirs{i}$.}
\par
Now suppose we have an assessment in the form of a set $\vassessment{n}$ of desirable gambles on~$\values^n$ for each $n\in\nats_0$.
We are looking for the (element-wise) smallest coherent, exchangeable and time-consistent family $\vdesirs{n}$, $n\in\nats_0$ that includes this assessment in the sense that $\vassessment{n}\subseteq\vdesirs{n}$  for all $n\in\nats_0$, which is equivalent to $\vmult{n}(\vassessment{n})\subseteq\vmult{n}(\vdesirs{n})$ for all $n\in\nats_0$, which is in turn---because of Eq.~\eqref{eq:infinite-representation-backtooriginal}---equivalent to
\begin{equation}
  \bigcup_{n\in\nats_0}\vmult{n}(\vassessment{n})
  \subseteq\bigcup_{n\in\nats_0}\vmult{n}(\vdesirs{n})
  \eqqcolon\vtdesirs,
\end{equation}
a condition formulated in terms of the frequency representation $\vtdesirs$ of the family $\vdesirs{n}$, $n\in\nats_0$.
The existence of this smallest family is implied by what we found in the previous paragraph.
If we combine all this with the arguments in the previous section, we are led to the following theorem.

\begin{theorem}\label{theo:exchangeable-natural-extension-infinite}
  Suppose we have an assessment in the form of a set $\vassessment{n}$ of desirable gambles on $\values^n$ for each $n\in\nats_0$, and consider the corresponding set of polynomials:
  \begin{equation}
     \passessment\coloneqq\bigcup_{n\in\nats_0}\vmult{n}(\vassessment{n}).
  \end{equation}
  Then there is a coherent, exchangeable and time-consistent family $\vdesirs{n}$, $n\in\nats_0$ that includes this assessment iff $\passessment$ avoids Bernstein non-positivity, and in that case $\bernnatex\passessment$ is the frequency representation of the (element-wise)  smallest  coherent, exchangeable and time-consistent family that includes this assessment.
\end{theorem}

\section{Extending finite exchangeable sequences}\label{sec:seq-ext-finite}
Suppose we have $n$ random variables $\rv_1$, \dots, $\rv_n$, that a subject judges to be exchangeable, and for which he has an assessment $\vassessment{n}$ of desirable gambles on $\values^n$, with corresponding count representation $\vcdesirs{n}=\vmuhy{n}(\vassessment{n})$.
We here answer the question of when it is possible and how, if so, to extend such a sequence to a longer, finite or infinite sequence that is still exchangeable.

\subsection{Extension to a longer, finite exchangeable sequence}\label{sec:seq-ext-to-finite}
In this section we ask: \emph{Can the assessment~$\vassessment{n}$ be extended to a coherent exchangeable model for $n+k$ variables? And if so, what is the most conservative such extended model?}
\par
It is well-known \citep{diaconis1980} that when the subject's assessment is an exchangeable linear prevision, such an extension is not generally possible.
In the much more general case that we are considering here, we now look at our Theorems~\ref{theo:exnatex} and~\ref{theo:finite-representation-exnatex}  to provide us with an elegant answer: the extension problem considered here is a special case of the one studied in Section~\ref{sec:exchangeable-natural-extension}.
\par
Indeed, since any gamble $f$ on the first $n$ variables $\rv_1$, \dots, $\rv_n$ corresponds to the gamble $\vexten{n}{n+k}(f)$ on the $n+k$ variables $\rv_1$, \dots, $\rv_n$, \dots, $\rv_{n+k}$, we see that the assessment $\vassessment{n}$ corresponds to an assessment
\begin{equation}
  \vassessment{n+k}
  \coloneqq\vexten{n}{n+k}(\vassessment{n})
  =\set{\vexten{n}{n+k}(f)}{f\in\vassessment{n}}
\end{equation}
of desirable gambles on $\values^{n+k}$.
It is then clear from  Theorem~\ref{theo:exnatex} that
\begin{inparaenum}[(i)]
\item $\vassessment{n}$ can  be extended to a coherent exchangeable model for $n+k$ variables iff  this $\vassessment{n+k}$ avoids non-positivity under exchangeability; and if such is the case, that
\item the smallest such coherent exchangeable extension is given by $\exnatex{\vassessment{n+k}}$.
\end{inparaenum}
\par
But we know from Theorem~\ref{theo:finite-representation-exnatex}  that it is easier to express this in terms if the count representations.
Since moreover, by Eq.~\eqref{eq:enl-identity},
\begin{equation}
  \vmuhy{n+k}(\vassessment{n+k})
  =\vmuhy{n+k}\bigl(\vexten{n}{n+k}(\vassessment{n})\bigr)
  =\venl{n}{n+k}\bigl(\vmuhy{n}(\vassessment{n})\bigr)
  =\venl{n}{n+k}(\vcdesirs{n}),
\end{equation}
we are led to the following simple solution to the extension problem.

\begin{theorem}
  Consider $n$ and $k$ in $\nats_0$.
  An assessment $\vassessment{n}$ of desirable gambles on $\values^n$, with corresponding count representation $\vcdesirs{n}\coloneqq\vmuhy{n}(\vassessment{n})$, can be extended to a coherent exchangeable model for $n+k$ variables iff\/ $\venl{n}{n+k}(\vcdesirs{n})$ avoids non-positivity.
  In that case the most conservative such coherent exchangeable model has count representation $\natex[][\big]{\venl{n}{n+k}(\vcdesirs{n})}$.
\end{theorem}

\begin{example}\label{ex:seq-ext-to-finite}
  In the context of our running example, take $n=2$ and consider the assessment $\vassessment{n}\coloneqq\{f\}$ where $f$ is the gamble on $\values^n$ given by
  \begin{equation*}
    f(b,b)=f(w,w)\coloneqq-3\text{ and }f(b,w)=f(w,b)\coloneqq1,
  \end{equation*}
  whence
  \begin{equation*}
    \vcdesirs{n}
    \coloneqq\vmuhy{n}(\vassessment{n})
    \coloneqq\set{g\in\gambles(\vcnts{n})}{g(2,0)=g(0,2)=-3\text{ and }g(1,1)=1}.
  \end{equation*}
  This singleton assessment avoids non-positivity under exchangeability and could be interpreted to express a strong belief that both colours will appear on the upcoming two draws, which could, e.g., be based on an observation of one black and one white marble being put in a seemingly empty urn.
  Now, let us see if this assessment can be extended to an exchangeable model for $n+k\coloneqq3$ variables:
  Let $g'\coloneqq\venl{n}{n+k}(g)$, then
  \begin{equation*}
  \begin{aligned}
    g'(3,0)&=g(2,0)=-3, &\quad g'(2,1)&=\tfrac{1}{3}g(2,0)+\tfrac{2}{3}g(1,1)=-\tfrac{1}{3},\\
    g'(0,3)&=g(0,2)=-3, & g'(1,2)&=\tfrac{1}{3}g(0,2)+\tfrac{2}{3}g(1,1)=-\tfrac{1}{3}.
  \end{aligned}
  \end{equation*}
  The gamble $g'$ is clearly non-positive, so the assessment cannot be extended to a coherent exchangeable model.
  Learning that there are more marbles in the urn would force us to revise the initial assessment: making this assessment when there are (at least) three balls in the urn leads to a sure loss.~$\blacklozenge$
\end{example}

\subsection{Extension to an infinite exchangeable sequence}\label{sec:seq-ext-to-infinite}
Let us now extend the course of reasoning in the previous section to make it deal with infinite sequences of random variables.
So in this section we ask: \emph{Can the assessment~$\vassessment{n}$ be extended to a coherent exchangeable model for an infinite sequence of variables? And if so, what is the most conservative such extended model?}
\par
Here, we look at Theorem~\ref{theo:exchangeable-natural-extension-infinite}  to provide us with an elegant answer: the present extension problem is a special case of that studied in Section~\ref{sec:exchangeable-natural-extension-infinite}.
Indeed, the set of desirable gambles $\vassessment{n}$ corresponds to an assessment of polynomials  $\vmult{n}(\vassessment{n})$, leading to the following simple solution to the extension problem.

\begin{theorem}\label{theo:extendability:infinite}
  Consider $n\in\nats_0$.
  An assessment $\vassessment{n}$ of desirable gambles on $\values^n$ can be extended to a coherent, exchangeable and time-consistent family iff\/ $\vmult{n}(\vassessment{n})$ avoids Bernstein non-positivity.
  In that case the most conservative such family has frequency representation $\bernnatex[\big]{\vmult{n}(\vassessment{n})}$.
\end{theorem}

\begin{example}
  The singleton assessment~$\vassessment{n}$ of Example~\ref{ex:seq-ext-to-finite} can very quickly be seen to not be extendable to a coherent, exchangeable and time-consistent family, because the single polynomial in $\vmult{n}(\vassessment{n})$---depicted below---is strictly negative, which by Proposition~\ref{prop:ranges} assures us it incurs Bernstein non-positivity.
  
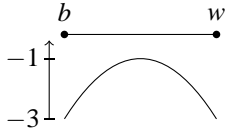
\begin{center}
    \begin{tikzpicture}[xscale=2,yscale=.4]
      \draw (0,-.2) node[fill,circle,inner sep=1pt,label=90:$b$] {} -- (1,-.2) node[fill,circle,inner sep=1pt,label=90:$w$] {};
      \draw[->] (-.1,-3) -- (-.1,-.4);
      \path (-.1,-3) \ytick node[left] {$-3$} (-.1,-1) \ytick node[left] {$-1$};
      \draw plot[domain=0:1] (1-\x,{-3*\x*\x + -3*(1-\x)*(1-\x) + 2*\x*(1-\x)});
    \end{tikzpicture}
  \end{center}
  Had the parabola's top value been~$0$, this would not have been the case due to BP\ref{item:bern-positive}.~$\blacklozenge$
\end{example}

\section{Conclusions}
We have shown that modelling a finite or infinite exchangeability assessment using sets of desirable gambles is not only possible, but also quite elegant.
Our results indicate that, using sets of desirable gambles, it is conceptually easy to reason about exchangeable sequences.
\par
Calculating the natural extension and updating are but simple geometrical operations: taking unions, sums and positive hulls and taking intersections, respectively.
This approach has the added advantage that the exchangeability assessment is preserved under updating, also when the conditioning event has lower probability zero, which does not hold when using (lower) previsions (although this might be remedied by using full conditional measures, for which \citet{Cozman-Seidenfeld-2009} give a good number of references).
\par
Using our Finite Representation Theorem, reasoning about finite exchangeable sequences can be reduced to reasoning about count vectors or (polynomials of) frequency vectors.
Working with these representations automatically guarantees that exchangeability is satisfied.
The representation for the natural extension and for updated models can be derived directly from the representation of the original model, without having to go back to the (more complex) world of sequences.
\par
Moreover, using our Infinite Representation Theorem, reasoning about infinite exchangeable sequences is reduced to reasoning about (polynomials of) frequency vectors.
Doing this automatically guarantees that, next to the exchangeability of finite subsequences, time consistency of these subsequences is satisfied.
Again, the representation for the natural extension and for updated models can be derived directly from the representation of the original model.
\par
Additionally, using our results about representation and natural extension, we have shown when and how finite exchangeable sequences can be extended to longer, finite or infinite exchangeable sequences.
However, we suspect there may be a more elegant characterisation of $\vnpospoly$ than the one given above, which might make the characterisation in Theorem~\ref{theo:extendability:infinite} more efficient to implement in terms of computer algorithms.
\par
What are the advantages of our approach?
It makes it easy for us to represent and reason with a finite number of expert assessments, and to see what its consequences are under exchangeability.
Also, we have seen that there are simple geometrical representations and interpretations of coherence and exchangeability: due to the symmetry, the assessments can be represented in simpler, lower dimensional spaces, and there are linear maps effecting that representation.
\par
The conceptual techniques employed in this paper are not restricted in use to a treatment of exchangeability.
They could be applied to other structural assessments, e.g., invariance assessments, as long as this assessment allows us to identify a characterising set of weakly desirable gambles that is sufficiently well-behaved (cf.~the first paragraph of Section~\ref{sec:xchdef}).
This idea was briefly taken up by one of us in another paper \cite{cooman2005c}, but clearly merits further attention.
\par
Thinking in even broader terms, we feel that using sets of desirable gambles can provide a refreshing and fruitful approach to many problems in uncertainty modelling, not only those related to structural assessments.
\par
While writing this paper, we regularly wondered what Henry Kyburg would have thought about it.
The topic surely has connections with his interests: exchangeability is an important basic assumption used in many models for statistical inference and our use of a model for uncertainty that is not just a precise probability, sets of desirable gambles.
What we tried to do in this paper is in some sense clarify, in a very general setting, what the consequences are of an assessment of exchangeability.
We know from his work that he thought it important for people to realise they are quite strong \citep[p.~111ff., p.~122ff]{kyburg1974}.
The preceding pages elaborately underline this point, and we agree it is an important one.
\par
Although this paper sprouted from minds mildly seduced by subjectivist betting frameworks, nothing in it precludes using it objectively.
We say this with a slightly mischievous smile, mirroring a similar twinkle in Henry's eyes when we met him last, at a 2005 conference in Pittsburgh.

\section*{Acknowledgements}
The authors wish to thank Teddy Seidenfeld for sharing some of the things he knows about Henry Kyburg and his work; they were both interesting and relevant.
Erik Quaeghebeur was supported by a Fellowship of the Belgian American Educational Foundation and wishes to thank Carnegie Mellon University's Department of Philosophy for its hospitality.

\appendix

\section{Proofs}\label{app:proofs}

\begin{proof}[\bfseries Proof of Theorem~\ref{theo:natex}]
  It follows from the fact that $\alldesirs[\spacecone](\pspace)$ is closed under arbitrary non-empty intersections, the definition of $\natex[\spacecone]{\assessment}$, and the fact that $\subspace$ is not coherent relative to $\spacecone$ [because $\cone\subset\subspace$], that the last four statements are equivalent.
  \par
  Next, we prove that (\ref{item:Aapl})$\Leftrightarrow$(\ref{item:Asubsetcoh}):
  \begin{compactitem}
    \item[$\Leftarrow$] Assume that~$\assessment$ is included in some set of desirable gambles $\desirs$ that is coherent relative to $\spacecone$.
    Since $\desirs=\posi(\desirs)$, $\desirs$~avoids non-positivity relative to $\spacecone$ by~D\ref{item:Dnonzeroapl}, and therefore so do all its subsets, including~$\assessment$.
    \item[$\Rightarrow$] Conversely, assume that $\assessment$ avoids non-positivity relative to $\spacecone$.
      For notational convenience, let  $\desirs^*\coloneqq\posi(\positive\cup\assessment)$.
      It is clear that $\desirs^*$ satisfies D\ref{item:Dapg}, D\ref{item:Dscl} and~D\ref{item:Dcmb}.
      Consider any $f\in\desirs^*$, so there are ${n\geq1}$, real $\lambda_k>0$, $f_k\in\positive\cup\assessment$ such that $f=\sum_{k=1}^n\lambda_kf_k$.
      Let $I\coloneqq\set{k\in\{1,\dots,n\}}{f_k\succ0}$, then $f_\ell\in\assessment$ for all $\ell\notin I$ and  ${f=f_0+\sum_{\ell\notin I}\lambda_kf_k}$ with $f_0\coloneqq\sum_{\ell\in I}\lambda_kf_k\succ0$.
      It therefore follows from the assumption that $\sum_{\ell\not\in I}\lambda_kf_k\not\preceq0$ and therefore {\itshape a fortiori} $f\not\preceq0$, so $\desirs^*$ also satisfies D\ref{item:Dnonzero} (or D\ref{item:Dnonzeroapl}), and is therefore coherent relative to $\spacecone$.
  \end{compactitem}
  Finally, we prove that  $\natex{\assessment}=\desirs^*$ whenever any (and hence all) of the equivalent statements (\ref{item:Aapl})--(\ref{item:natextAsmallestcoh}) hold.
  Any coherent set of desirable gambles that includes~$\assessment$, must also include $\desirs^*$, by the axioms D\ref{item:Dapg}, D\ref{item:Dscl}, and D\ref{item:Dcmb}.
  Since we have proved above that $\desirs^*$  also satisfies D\ref{item:Dnonzero} and is therefore coherent relative to $\spacecone$, it is the smallest  set of desirable gambles that is coherent relative to $\spacecone$ and includes~$\assessment$. Hence it is equal to $\natex[\spacecone]{\assessment}$, by (\ref{item:natextAsmallestcoh}).
\end{proof}

\begin{proof}[\bfseries Proof of Proposition~\ref{prop:maximality}]
  We first prove sufficiency.
  Assume that Eq.~\eqref{eq:maximality} holds.
  Consider any $\desirs'$ in $\alldesirs[\spacecone](\pspace)$ such that $\desirs\subseteq\desirs'$, then we prove that also $\desirs'\subseteq\desirs$.
  To this effect, consider any $f\in\desirs'$, so $-f\notin\desirs'$ by coherence, and therefore also $-f\notin\desirs$.
  Now invoke Eq.~\eqref{eq:maximality} to find that $f\in\desirs$.
  \par
  Next, we turn to necessity.
  Assume that $\desirs$ is maximal, consider any $f\in\subspace_0$, and assume that $f\notin\desirs$.
  We have to prove that $-f\in\desirs$.
  By Lemma~\ref{lem:maximality}, we get that $\posi(\desirs\cup\{-f\})\in\alldesirs[\spacecone](\pspace)$, but since $\desirs\subseteq\posi(\desirs\cup\{-f\})$ and $\desirs$ is maximal, we conclude that $\desirs=\posi(\desirs\cup\{-f\})$ and therefore indeed $-f\in\desirs$.
\end{proof}

\begin{lemma}\label{lem:maximality}
   Let\/ $\subspace$ be a linear subspace of\/ $\gambles(\pspace)$ and let $\cone\subset\subspace$ be a convex cone containing the zero gamble $0$.
   Let $\desirs\in\alldesirs[\spacecone](\pspace)$ and let $f$ be any non-zero gamble in $\subspace$.
   Then $f\notin\desirs$ implies that $\desirs\cup\{-f\}$ avoids non-positivity relative to $\spacecone$, and therefore $\natex[\spacecone]{\desirs\cup\{-f\}}=\posi(\desirs\cup\{-f\})$ is coherent relative to $\spacecone$:
   \begin{equation}
     (\forall f\in\subspace_0)
     \bigl(f\notin\desirs
     \then\posi(\desirs\cup\{-f\})
     \in\alldesirs[\spacecone](\pspace)\bigr).
   \end{equation}
\end{lemma}

\begin{proof}[\bfseries Proof of Lemma~\ref{lem:maximality}]
  We give a proof by contradiction.
  Let $f\in\subspace_0\setminus\desirs$ and assume that $\desirs\cup\{-f\}$ does not avoid non-positivity relative to $\spacecone$.
  This means that ${\posi(\desirs\cup\{-f\})}\cap\nonpositive\neq\emptyset$, and since $\desirs$ does avoid non-positivity  relative to $\spacecone$, this tells us that there are $n\in\nats_0$, $f_1$, \dots, $f_n$ in $\desirs$, $\lambda$ in $\reals^+_0$, and $\lambda_1$, \dots, $\lambda_n$ in $\reals^+$ such that
  \begin{equation}
    \sum_{k=1}^n\lambda_kf_k+\lambda(-f)\preceq0
    \text{ and therefore }
    f\succeq\sum_{k=1}^n\frac{\lambda_k}{\lambda}f_k.
  \end{equation}
  Then obviously $f\in\desirs$ since $f\neq0$, a contradiction.
  The rest of the proof now follows from Theorem~\ref{theo:natex} and $\positive\subseteq\desirs$.
\end{proof}

\begin{proof}[\bfseries Proof of Theorem~\ref{theo:maximality}]
  Sufficiency follows readily from Theorem~\ref{theo:natex}.
  \par
  For necessity, assume that  $\assessment$ avoids non-positivity relative to $\spacecone$, and consider the set ${\uparrow}\assessment\coloneqq\set{\desirs\in\alldesirs[\spacecone](\pspace)}{\assessment\subseteq\desirs}$. This set is non-empty by Theorem~\ref{theo:natex}, and partially ordered by set inclusion.
  We show that this poset has a maximal element, which is then automatically also a maximal element of $\alldesirs[\spacecone](\pspace)$.
  \par
  Consider any chain $\chain\subseteq{\uparrow}\assessment$.
  We show that $\bigcup\chain$ avoids non-positivity relative to $\spacecone$.
  Consider arbitrary $n\in\nats_0$, $f_1$, \dots, $f_n$ in $\bigcup\chain$.
  $f_k\in\bigcup\chain$ means that there is some $\desirs_k\in\chain$ such that $f_k\in\desirs_k$, and therefore $\{f_1,\dots,f_n\}\subseteq\bigcup_{k=1}^n\desirs_n\eqqcolon\tilde{\desirs}$.
  But $\tilde{\desirs}\in\chain$ because $\chain$ is a chain, and therefore $\tilde{\desirs}$ is coherent relative to $\spacecone$.
  This implies that $\nonpositive\cap\tilde{\desirs}=\emptyset$, and therefore {\itshape a fortiori}  $\nonpositive\cap\posi(\{f_1,\dots,f_n\})=\emptyset$.
  So we find that $\bigcup\chain$ indeed avoids non-positivity relative to $\spacecone$.
  \par
  By Theorem~\ref{theo:natex}, $\natex[\spacecone]{\bigcup\chain}=\posi(\bigcup\chain)$ is coherent relative to $\spacecone$ and includes $\bigcup\chain$, so $\desirs\subseteq\posi(\bigcup\chain)$ for all $\desirs\in\chain$.
  Because also $\assessment\subseteq\posi(\bigcup\chain)$, we have just shown that every chain $\chain$ in  the poset ${\uparrow}\assessment$ has an upper bound $\posi(\bigcup\chain)$ in ${\uparrow}\assessment$.
  By Zorn's Lemma, ${\uparrow}\assessment$ has a maximal element.
\end{proof}

\begin{proof}[\bfseries Proof of Corollary~\ref{cor:maximality}]
   We use the notation $\mathbb{M}^*\coloneqq\set{\desirs\in\maxdesirs[\spacecone](\pspace)}{\assessment\subseteq\desirs}$ for the sake of brevity.
   \par
   If $\assessment$ does not avoid non-positivity relative to $\spacecone$, then $\mathbb{M}^*=\emptyset$, by Theorem~\ref{theo:natex}, so $\bigcap\mathbb{M}^*=\subspace$.
   Again by Theorem~\ref{theo:natex}, also  $ \natex[\spacecone]{\assessment}=\subspace$.
   \par
   Assume, therefore, that $\assessment$ avoids non-positivity relative to $\spacecone$.
   Then $\mathbb{M}^*\neq\emptyset$ by Theorem~\ref{theo:maximality}.
   Since $\bigcap\mathbb{M}^*$ is coherent relative to $\spacecone$ and $\assessment\subseteq\bigcap\mathbb{M}^*$, we infer from Theorem~\ref{theo:natex} that $\natex[\spacecone]{\assessment}\subseteq\bigcap\mathbb{M}^*$.
   Assume {\itshape ex absurdo} that  $\natex[\spacecone]{\assessment}\subset\bigcap\mathbb{M}^*$, so there is some (non-zero) $f\in\bigcap\mathbb{M}^*$ such that $f\notin\natex[\spacecone]{\assessment}$, and therefore $\natex[\spacecone]{\assessment}\cup\{-f\}$ avoids non-positivity relative to $\spacecone$, by Lemma~\ref{lem:maximality}.
   By Theorem~\ref{theo:maximality}, there is some $\desirs^*$ in $\maxdesirs[\spacecone](\pspace)$ such that $\natex[\spacecone]{\assessment}\cup\{-f\}\subseteq\desirs^*$.
   On the one hand, we infer that $\assessment\subseteq\desirs^*$, so $\desirs^*\in\mathbb{M}^*$, and therefore $f\in\desirs^*$.
   On the other hand, we infer that $-f\in\desirs^*$, which contradicts ${f\in\desirs^*}$, since $\desirs^*$ is coherent.
\end{proof}

\begin{proof}[\bfseries Proof of Proposition~\ref{prop:weak-desirability}]
  The defining property of any gamble~$f$ in~$\weakly{\desirs}$ is that ${f+f'\in\desirs}$ for all gambles~$f'$ in~$\desirs$.
  \begin{compactitem}
    \item[WD\ref{item:WDapl}.] Let $f<0$; let $f'=-f/2$ then $f'>0$ and therefore $f'\in\desirs$, by~D\ref{item:Dnonzero}.
      But $f''=f+f'=f/2<0$ and thus, by~D\ref{item:Dnonzeroapl}, $f''\not\in\desirs$.
      Hence indeed $f\notin\weakly{\desirs}$.
    \item[WD\ref{item:WDapg}.] Since clearly $0+\desirs=\desirs$, we see that $0\in\weakly{\desirs}$.
      If ${f>0}$ then $f\in\desirs$ by D\ref{item:Dapg}, and therefore $f\in\weakly{\desirs}$, because $\desirs\subseteq\weakly{\desirs}$.
    \item[WD\ref{item:WDscl}.] Consider $f\in\weakly{\desirs}$.
      If $\lambda=0$ then $\lambda f=0\in\weakly{\desirs}$ by~WD\ref{item:WDapg}.
      Assume therefore that $\lambda>0$.
      Consider any $f'\in\desirs$.
      Then $f'/\lambda\in\desirs$ by~D\ref{item:Dscl}, so $f+f'/\lambda\in\desirs$, and therefore $\lambda f+f'\in\desirs$, again by~D\ref{item:Dscl}.
      Hence indeed $\lambda f\in\weakly{\desirs}$.
    \item[WD\ref{item:WDcmb}.] Consider $f_1,f_2\in\weakly{\desirs}$, and any $f'\in\desirs$, so ${f'/2\in\desirs}$ by~D\ref{item:Dscl}.
      Then $f_1+f'/2\in\desirs$ and $f_2+f'/2\in\desirs$, and therefore $f_1+f_2+f'\in\desirs$, by~D\ref{item:Dcmb}.
      Hence indeed ${f_1+f_2\in\weakly{\desirs}}$. \qedhere
  \end{compactitem}
\end{proof}

\begin{proof}[\bfseries Proof of Theorem~\ref{theo:lpr-desirs}]
  To prove that $\lpr_\desirs$ is real-valued, we prove that $\lpr_\desirs(f)$ is bounded for all gambles $f\in\gambles(\pspace)$---which are bounded and real-valued by definition.
  It follows from~D\ref{item:Dnonzeroapl} that if $f-\mu\in\desirs$, then $f\not\leq\mu$, so $\sup{f}>\mu$, whence $\lpr_\desirs(f)\leq\sup f<+\infty$.
  It follows from D\ref{item:Dapg} that $f-\mu\in\desirs$ if $f-\mu>0$; let $\mu$ be any real number such that $\mu<\inf f$, then $f-\mu>f-\inf f\geq0$, so ${f-\mu\in\desirs}$, whence $\lpr_\desirs(f)\geq\inf f>\mu>-\infty$.
  \par
  To prove the equality of $\lpr_\desirs$ and $\lpr_{\weakly{\desirs}}$, consider any gamble $f\in\gambles(\pspace)$.
  Since $\desirs\subseteq\weakly{\desirs}$, we immediately get that
  \begin{equation}
    \set{\mu\in\reals}{f-\mu\in\desirs}
    \subseteq\set{\mu\in\reals}{f-\mu\in\weakly{\desirs}}
  \end{equation}
  and therefore $\lpr_\desirs(f)\leq\lpr_{\weakly{\desirs}}(f)$.
  Conversely, consider any $\alpha>0$, then $\alpha\in\desirs$ by coherence [D\ref{item:Dapg}], and therefore
  \begin{align}
    \set{\mu\in\reals}{f-\mu\in\weakly{\desirs}}
    &\subseteq\set{\mu}{f-\mu+\alpha\in\desirs}\\
    &=\alpha+\set{\mu}{f-\mu\in\desirs},
  \end{align}
  whence $\lpr_{\weakly{\desirs}}(f)\leq\alpha+\lpr_\desirs(f)$.
  Since this holds for all ${\alpha>0}$, we also have $\lpr_{\weakly{\desirs}}(f)\leq\lpr_\desirs(f)$.
  \par
  Next, consider any $f\in\weakly{\desirs}$.
  Because ${f=f-0}$ this tells us that $\lpr_\desirs(f)=\lpr_{\weakly{\desirs}}(f)\geq0$.
  \par
  The rest of the proof is now standard, see for instance \citep[Section~6]{walley2000}.
\end{proof}

\begin{proof}[\bfseries Proof of Proposition~\ref{prop:marginal-desirability}]
  Since it follows from Theorem~\ref{theo:lpr-desirs} that $\lpr_\desirs(f-\lpr_\desirs(f))=\lpr_\desirs(f)-\lpr_\desirs(f)=0$   for all gambles~$f$, it follows that $\marginally{\desirs}\subseteq\set{f\in\gambles(\pspace)}{\lpr_\desirs(f)=0}$.
  For the converse inequality, assume that $\lpr_\desirs(f)=0$ holds; then $f=f-\lpr_\desirs(f)\in\marginally{\desirs}$.
  \par
  This also means that $\lpr_\desirs(g)=0$ iff $g\in\marginally{\desirs}$, so for every gamble~$f$ we can write:
  \begin{align}
    \lpr_{\marginally{\desirs}}(f)
    &=\sup\set{\mu\in\reals}{f-\mu\in\marginally{\desirs}}\\
    &=\sup\set{\mu\in\reals}{\lpr_\desirs(f-\mu)=0}\\
    &=\sup\set{\mu\in\reals}{\mu=\lpr_\desirs(f)}
    =\lpr_\desirs(f),
  \end{align}
  which proves the equality of $\lpr_{\marginally{\desirs}}$ and $\lpr_\desirs$.
\end{proof}

\begin{proof}[\bfseries Proof of Proposition~\ref{prop:updating-preserves-coherence}]
  We need to prove that the appropriate versions of D\ref{item:Dnonzero}--D\ref{item:Dcmb} hold for $\update{\desirs}{B}$, with $\subspace=\update{\gambles(\pspace)}{B}$ and $\cone=\update{\gambles(\pspace)}{B}\cap\gambles^+(\pspace)$.
  For D\ref{item:Dnonzero}, consider $f\in\update{\gambles(\pspace)}{B}$ and assume that $f=0$.
  Then by coherence $f\not\in\desirs$ and hence $f\not\in\update{\desirs}{B}$.
  For D\ref{item:Dapg}, consider $f\in\update{\gambles(\pspace)}{B}$ and assume that $f>0$.
  Then by coherence $f\in\desirs$ and hence $f\in\update{\desirs}{B}$.
  The proof for D\ref{item:Dscl} is similar to the one for~D\ref{item:Dcmb}.
  For D\ref{item:Dcmb}, consider $f_1,f_2\in\update{\desirs}{B}$, then on the one hand ${f_1,f_2\in\desirs}$ and therefore $f_1+f_2\in\desirs$ by coherence; and on the other hand ${f_1,f_2\in\update{\gambles(\pspace)}{B}}$ and therefore $f_1+f_2=I_Bf_1+I_Bf_2=I_B(f_1+f_2)$, so ${f_1+f_2\in\update{\gambles(\pspace)}{B}}$ and hence $f_1+f_2\in\update{\desirs}{B}$.
\end{proof}

\begin{proof}[\bfseries Proof of the equivalences in Definition~\ref{def:exchangeability}]
  That (\ref{item:xch-permutssubsetdef})$\Leftrightarrow$(\ref{item:xch-permutsreldef}) and (\ref{item:xch-averagesubsetdef})$\Leftrightarrow$(\ref{item:xch-averagereldef}) is an immediate consequence of the definition of weak desirability.
  We show that (\ref{item:xch-permutssubsetdef})$\Leftrightarrow$(\ref{item:xch-averagesubsetdef}).
  For the `$\Rightarrow$' part, observe that $f-\vex{N}(f)=\frac{1}{N!}\sum_{\pi\in\permuts_N}[f-\pi^tf]\in\weakly{\desirs}$, since $\weakly{\desirs}$ is a convex cone by Proposition~\ref{prop:weak-desirability}.
  The `$\Leftarrow$' part follows from $\weakly{\permuts_N}\subseteq\weakly{\average_N}$, i.e., from Eq.~\eqref{eq:span}.
\end{proof}

\begin{proof}[\bfseries Proof of Proposition~\ref{prop:permutability}]
  Consider $f\in\desirs$.
  Since $\pi^tf-f=(-f)-\pi^t(-f)\in\weakly{\permuts_N}$, we see that $\pi^tf=f+\pi^tf-f\in\desirs+\weakly{\permuts_N}\subseteq\desirs$, using the exchangeability condition of Definition~\ref{def:exchangeability}(\ref{item:xch-permutsreldef}).
\end{proof}

\begin{proof}[\bfseries Proof of Proposition~\ref{prop:exchangeability-ex}]
  The first statement is a consequence of the second, with $f'=\vex{N}(f)$, because then $\vex{N}(f')=\vex{N}(\vex{N}(f))=\vex{N}(f)$.
  For the second statement, consider arbitrary gambles $f$ and $f'$ on $\values^N$ such that $\vex{N}(f)=\vex{N}(f')$, and assume that $f\in\desirs$.
  We prove that then also ${f'\in\desirs}$.
  Since $\vex{N}(f)-f=(-f)-\vex{N}(-f)\in\weakly{\desirs}$ and ${f'-\vex{N}(f')\in\weakly{\desirs}}$, we see that $f'-f\in\weakly{\desirs}$ by WD\ref{item:WDcmb}, and therefore $f'=f+f'-f\in\desirs+\weakly{\desirs}\subseteq\desirs$.
\end{proof}

\begin{proof}[\bfseries Proof of Theorem~\ref{theo:exchangeability-lpr}]
  We give a circular proof.
  We first show that (\ref{item:xch-prevpermuts}) holds if $\lpr$ is exchangeable, i.e., if there is some coherent and exchangeable $\desirs$ such that $\lpr=\lpr_\desirs$.
  We already know from Theorem~\ref{theo:lpr-desirs} that $\lpr=\lpr_\desirs$ satisfies P\ref{item:Pasg}--P\ref{item:Phom}, because $\desirs$ is coherent.
  Consider any $f\in\weakly{\permuts_N}$.
  Since $\weakly{\permuts_N}\subseteq\weakly{\desirs}$, it also follows from Theorem~\ref{theo:lpr-desirs} that $\lpr_\desirs(f)\geq0$ and similarly $-\upr_\desirs(f)=\lpr_\desirs(-f)\geq0$ because also $-f\in\weakly{\permuts_N}$.
  Hence indeed $0\leq\lpr_\desirs(f)\leq\upr_\desirs(f)\leq0$, where the second inequality is a consequence of~P\ref{item:Pasg} and~P\ref{item:Psad}.
  \par
  That (\ref{item:xch-prevpermuts}) implies (\ref{item:xch-prevaverage}) follows the super-additivity of $\lpr$ and the sub-additivity of $\upr$.
  \par
  Finally, we show that (\ref{item:xch-prevaverage}) implies that $\lpr$ is exchangeable.
  The standard argument in \cite[Section~6]{walley2000} tells us that $\desirs'\coloneqq\set{f\in\gambles(\values^N)}{f>0\text{ or }\lpr(f)>0}$ is a coherent set of desirable gambles such that $\lpr_{\desirs'}=\lpr$.
  Now consider the set $\desirs\coloneqq\desirs'+\weakly{\average_N}$.
  We show that this $\desirs$ is a coherent and exchangeable set of desirable gambles, and that $\lpr_\desirs=\lpr$.
  It is clear from its definition that $\desirs$ satisfies D\ref{item:Dapg}, D\ref{item:Dscl} and D\ref{item:Dcmb}, so let us assume {\itshape ex absurdo} that $0\in\desirs$, meaning that there is some $f\in\desirs'$ such that $f'\coloneqq-f\in\weakly{\average_N}$.
  There are two possibilities.
  Either $f>0$, so $f'<0$, which contradicts Lemma~\ref{lem:permutation-consistency}.
  Or $\lpr(f)>0$.
  But it follows from~(\ref{item:xch-prevaverage}) and the coherence of the lower prevision $\lpr$ that ${0=\lpr(f+f')=\lpr(f)>0}$, a contradiction too.
  So $\desirs$ satisfies D\ref{item:Dnonzero} as well, and is therefore coherent. It is obvious that $\desirs$ is exchangeable:
  $
  \desirs+\weakly{\average_N}
  =\desirs'+\weakly{\average_N}+\weakly{\average_N}
  =\desirs'+\weakly{\average_N}
  =\desirs
  $.
  The proof is complete if we can show that $\lpr=\lpr_\desirs$.
  Fix any gamble $f$.
  Observe that $f-\alpha\in\desirs$ iff there are $f'\in\desirs'$ and $f''\in\weakly{\average_N}$ such that $f-\alpha={f'+f''}$.
  But then it follows from the coherence of $\lpr$ and the assumption that $\lpr(f)=\alpha+\lpr(f'+f'')=\alpha+\lpr(f')\geq\alpha$, and therefore $\lpr_\desirs(f)\leq\lpr(f)=\lpr_{\desirs'}(f)$.
  For the converse inequality, we infer from $0\in\weakly{\average_N}$ that $\desirs'\subseteq\desirs$, and therefore $\lpr_{\desirs'}\leq\lpr_{\desirs}$.
\end{proof}

\begin{lemma}\label{lem:permutation-consistency}
  For all $f$ in $\weakly{\average_N}$, $f\not<0$.
\end{lemma}

\begin{proof}
  First of all, observe that for any gamble $f'$ on $\values^N$, if $f'>0$ then also $\vex{N}(f')>0$.
  Now consider $f\in\weakly{\average_N}$ and assume {\itshape ex absurdo} that $f<0$.
  Then $-f>0$ and therefore $-\vex{N}(f)=\vex{N}(-f)>0$, whence $\vex{N}(f)<0$.
  But since $f\in\weakly{\average_N}$ we also have that $\vex{N}(f)=0$, a contradiction.
\end{proof}

\begin{proof}[\bfseries Proof of Proposition~\ref{prop:anpune}]
  For the first statement, we have to prove that $\gambles^+_0(\values^N)+\weakly{\average_N}$ avoids non-positivity.
  Consider any $f'\in\weakly{\average_N}$ and any $f''\in\gambles^+_0(\values^N)$, then we have to prove that $f\coloneqq f'+f''\not\leq0$.
  There are two possibilities.
  Either $f'=0$ and then $f=f''>0$.
  Or $f'\neq0$, and then Lemma~\ref{lem:permutation-consistency} tells us that $f'\not\leq0$ and therefore {\itshape a fortiori} $f\not\leq0$.
  \par
  For the second statement, it clearly suffices to prove the `if' part.
  Assume therefore that $\assessment+\weakly{\average_N}$ avoids non-positivity.
  Consider any $f$ in ${\posi([\gambles^+_0(\values^N)\cup\assessment]+\weakly{\average_N})}$, so there are $n\geq1$, $\lambda_k\in\reals^+_0$, $f'\in\weakly{\average_N}$, $f_k\in\gambles^+_0(\values^N)\cup\assessment$ such that $f=f'+\sum_{k=1}^n\lambda_kf_k$.
  Let $I\coloneqq\set{k\in\{1,\dots,n\}}{f_k>0}$ then $f_\ell\in\assessment$ for all $\ell\notin I$, and $f=f_0+f'+\sum_{\ell\notin I}\lambda_\ell f_\ell$ with $f_0>0$.
  By assumption $f'+\sum_{\ell\notin I}\lambda_\ell f_\ell\not\leq0$, and therefore {\itshape a fortiori} $f\not\leq0$.
\end{proof}

\begin{proof}[\bfseries Proof of Theorem~\ref{theo:exnatex}]
  It is immediately clear from the fact that $\allexdesirs(\values^N)$ is closed under arbitrary non-empty intersections, the definition of $\exnatex{\assessment}$, and the fact that $\gambles(\values^N)$ is not a coherent set of desirable gambles, that the last four statements are equivalent.
  \par
  Next, we prove that (\ref{item:Aaplunderxch})$\Leftrightarrow$(\ref{item:Asubsetcohxch}).
  \begin{compactitem}
  \item[$\Leftarrow$] Assume that $\assessment$, and therefore also $\gambles^+_0(\values^N)\cup\assessment$, is included in some coherent and exchangeable set of desirable gambles $\desirs$.
    By exchangeability, we know $[\gambles^+_0(\values^N)\cup\assessment]+\weakly{\average_N}\subseteq\desirs+\weakly{\average_N}\subseteq\desirs$.
    Since $\posi(\desirs)=\desirs$ avoids non-positivity, so does any of its subsets, and therefore in particular
    $[\gambles^+_0(\values^N)\cup\assessment]+\weakly{\average_N}$.
    This means that $\assessment$ indeed avoids non-positivity under exchangeability.
  \item[$\Rightarrow$] Conversely, assume that $\assessment$ avoids non-positivity under exchangeability.
    For the sake of convenience, denote the set on the right-hand side of Eq.~\eqref{eq:posi-exnatex} by $\desirs^*$.
    It is clear that $\desirs^*$ satisfies D\ref{item:Dapg}, D\ref{item:Dscl} and~D\ref{item:Dcmb}.
    Consider any $f\in\desirs^*$, then $f\not\leq0$, precisely because $\assessment$ avoids non-positivity under exchangeability.
    Hence $\desirs^*$ also satisfies D\ref{item:Dnonzero}, and is therefore coherent.
    The exchangeability of $\desirs^*$ immediately follows from the fact that $\weakly{\average_N}+\natex{\assessment}+\weakly{\average_N}=\weakly{\average_N}+\natex{\assessment}$.  \end{compactitem}
  \par
  Finally, we prove Eqs.~\eqref{eq:posi-exnatex} and~\eqref{eq:posi-exnatex-better} whenever any (and hence all) of the equivalent statements (\ref{item:Aaplunderxch})--(\ref{item:exnatexA-smallestcoh}) holds.
  Eq.~\eqref{eq:posi-exnatex-better} follows from Eq.~\eqref{eq:posi-exnatex} and Theorem~\ref{theo:natex}, since $\weakly{\average_N}$ is a convex cone.
  Let us prove that $\exnatex{\assessment}=\desirs^*$.
  It is clear that any coherent and exchangeable set of desirable gambles that includes~$\assessment$, must also include $\desirs^*$, by the axioms D\ref{item:Dapg}, D\ref{item:Dscl}, and D\ref{item:Dcmb}.
  Since we have just proved above that $\desirs^*$ is coherent and exchangeable, it is the smallest coherent and exchangeable set of desirable gambles that includes~$\assessment$, and for this reason it is equal to $\exnatex{\assessment}$, by (\ref{item:exnatexA-smallestcoh}).
\end{proof}

\begin{proof}[\bfseries Proof of Corollary~\ref{cor:smallest-exchangeable}]
 This is an immediate consequence of Proposition~\ref{prop:anpune}(i) and Theorem~\ref{theo:exnatex}.
\end{proof}

\begin{proof}[\bfseries Proof of Proposition~\ref{prop:updating-preserves-exchangeability}]
  The coherence of $\restrict{\desirs}{\osample{x}}$ is guaranteed by Proposition~\ref{prop:updating-preserves-coherence}.
  We show that $\restrict{\desirs}{\osample{x}}$ is exchangeable.
  Consider arbitrary $f\in\gambles(\values^{\rest{n}})$, $\rest{\pi}\in\permuts_{\rest{n}}$ and $f_1\in\restrict{\desirs}{\osample{x}}$.
  Then we must show that $f_1+f-\rest{\pi}^tf\in\restrict{\desirs}{\osample{x}}$, or in other words that $I_{\event_{\osample{x}}}[f_1+f-\rest{\pi}^tf]\in\desirs$.
  But since $f_1\in\restrict{\desirs}{\osample{x}}$, we know that $I_{\event_{\osample{x}}}f_1\in\desirs$.
  And if we consider the permutation $\pi\in\permuts_N$ defined by
  \begin{equation}
    \pi(k)\coloneqq
    \begin{cases}
      k&\text{ $1\leq k\leq\obs{n}$}\\
      \obs{n}+\rest{\pi}(k-\obs{n})&\text{ $\obs{n}+1\leq k\leq N$},
    \end{cases}
  \end{equation}
  then clearly $I_{\event_{\osample{x}}}\rest{\pi}^tf=\pi^t(I_{\event_{\osample{x}}}f)$ and therefore ${I_{\event_{\osample{x}}}[f_1+f-\rest{\pi}^tf]}=I_{\event_{\osample{x}}}f_1+I_{\event_{\osample{x}}}f-\pi^t(I_{\event_{\osample{x}}}f)$ and this gamble belongs to~$\desirs$ because~$\desirs$ is exchangeable.
\end{proof}

\begin{proof}[\bfseries Proof of Proposition~\ref{prop:sufficiency}]
  Consider $\obs{\pi}\in\permuts_{\obs{n}}$ and any gamble $f$ on $\values^{\rest{n}}$.
  Assume that ${I_{\event_{\osample{x}}}f\in\desirs}$.
  \par
  We first prove that $I_{\event_{\obs{\pi}{\osample{x}}}}f\in\desirs$.
  Consider the permutation $\pi\in\permuts_N$ defined by
  \begin{equation}
    \pi(k)\coloneqq
    \begin{cases}
      \obs{\pi}^{-1}(k)&\text{ $1\leq k\leq\obs{n}$}\\
      k&\text{ $\obs{n}+1\leq k\leq N$},
    \end{cases}
  \end{equation}
  then clearly
  $
    \pi^t(I_{\event_{\osample{x}}}f)
    =(I_{\event_{\osample{x}}}f)\circ\pi
    =(I_{\event_{\osample{x}}}\circ\obs{\pi}^{-1})f
    =I_{\event_{\obs{\pi}\osample{x}}}f,
  $
  so it follows from Proposition~\ref{prop:permutability} that indeed $I_{\event_{\obs{\pi}{\osample{x}}}}f\in\desirs$.
  This already implies that $\restrict{\desirs}{\osample{x}}=\restrict{\desirs}{\obs{\pi}\osample{x}}$, and therefore also that $\restrict{\desirs}{\osample{x}}=\restrict{\desirs}{\osample{y}}$.
  \par
  Since $\desirs$ is coherent, it also follows from $I_{\event_{\osample{x}}}f\in\desirs$ and the reasoning above that
  $
    I_{\event_{\ocnt{m}}}f
    =\sum_{\osample{y}\in\atom{\ocnt{m}}}I_{\event_{\osample{y}}}f\in\desirs,
  $
  whence $\restrict{\desirs}{\osample{x}}\subseteq\restrict{\desirs}{\ocnt{m}}$.
  To prove the converse inequality, assume that ${I_{\event_{\ocnt{m}}}f\in\desirs}$.
  We know that $\oatom{m}=\set{\obs{\pi}\osample{x}}{\obs{\pi}\in\permuts_{\obs{n}}}$, and therefore for any $\osample{y}\in\oatom{m}$ we can pick a $\obs{\pi}_{\osample{y}}\in\permuts_{\obs{n}}$ such that  $\obs{\pi}_{\osample{y}}\osample{x}=\osample{y}$.
  With this~$\obs{\pi}_{\osample{y}}$ we construct a permutation $\pi_{\osample{y}}\in\permuts_N$ in the manner described above, which satisfies $\pi_{\osample{y}}^t(I_{\event_{\osample{x}}}f)=I_{\event_{\osample{y}}}f$.
  But then the exchangeability and coherence of $\desirs$ tell us that
  \begin{equation}
    I_{\event_{\ocnt{m}}}f
    +\sum_{\osample{y}\in\oatom{m}}
    [(I_{\event_{\osample{x}}}f)-\pi_{\osample{y}}^t(I_{\event_{\osample{x}}}f)]
    =I_{\event_{\ocnt{m}}}f
    +f\sum_{\osample{y}\in\oatom{m}}
    [I_{\event_{\osample{x}}}-I_{\event_{\osample{y}}}]
    =\abs{\oatom{m}}fI_{\event_{\osample{x}}}
  \end{equation}
  belongs to $\desirs$, whence also $I_{\event_{\osample{x}}}f\in\desirs$, by coherence.
\end{proof}

\begin{proof}[\bfseries Proof of Theorem~\ref{theo:finite-representation}]
  We begin with the sufficiency part.
  Assume that there is some  coherent set $\cdesirs$ of desirable gambles on $\vcnts{N}$ such that $\desirs=(\vmuhy{N})^{-1}(\cdesirs)$.
  We show that~$\desirs$ is coherent and exchangeable, and that $\cdesirs=\vmuhy{N}(\desirs)$.
  \par
  We first show that $\desirs$ is coherent.
  For D\ref{item:Dnonzero}, consider $f\in\gambles(\values^N)$ with $f=0$.
  Then obviously also $\vmuhy{N}(f)=0$ and therefore $\vmuhy{N}(f)\not\in\cdesirs$.
  Hence $f\notin\desirs$.
  For D\ref{item:Dapg}, let ${f>0}$.
  Then obviously also ${\vmuhy{N}(f)>0}$, and therefore $\vmuhy{N}(f)\in\cdesirs$.
  Hence $f\in\desirs$.
  The proof for D\ref{item:Dscl} is similar to the one for D\ref{item:Dcmb}.
  For D\ref{item:Dcmb}, let $f_1,f_2\in\desirs$.
  Then ${g_1\coloneqq\vmuhy{N}(f_1)\in\cdesirs}$ and $g_2\coloneqq\vmuhy{N}(f_2)\in\cdesirs$.
  This implies that $\vmuhy{N}(f_1+f_2)=g_1+g_2\in\cdesirs$, so again $f_1+f_2\in\desirs$.
  \par
  To show that $\desirs$ is exchangeable, consider any $f\in\desirs$ and $f'\in\weakly{\average_N}$.
  We have to show that $f+f'\in\desirs$.
  It is clear that $\vmuhy{N}(f+f')=\vmuhy{N}(f)+0=\vmuhy{N}(f)\in\cdesirs$.
  Hence $f+f'\in(\vmuhy{N})^{-1}(\cdesirs)$, so indeed $f+f'\in\desirs$.
  \par
  We show that $\cdesirs=\vmuhy{N}(\desirs)$.
  Consider any gamble $g\in\gambles(\vcnts{N})$, then using Eq.~\eqref{eq:muhy-identities}, $\vmuhy{N}(\vocntf{N}(g))=g$.
  Since by assumption $\desirs=(\vmuhy{N})^{-1}(\cdesirs)$, we see that
  \begin{equation}
    g\in\cdesirs
    \Leftrightarrow\vmuhy{N}(\vocntf{N}(g))\in\cdesirs
    \Leftrightarrow\vocntf{N}(g)\in\desirs.
  \end{equation}
  This shows that $\cdesirs=\set{g\in\gambles(\vcnts{N})}{\vocntf{N}(g)\in\desirs}$.
  We show that also $\cdesirs=\vmuhy{N}(\desirs)$.
  Let $g\in\cdesirs$, then we have just proved that $\vocntf{N}(g)\in\desirs$, and therefore, using Eq.~\eqref{eq:muhy-identities}, $g=\vmuhy{N}(\vocntf{N}(g))\in\vmuhy{N}(\desirs)$.
  Conversely, let $g\in\vmuhy{N}(\desirs)$.
  Then there is some $f\in\desirs$ such that $g=\vmuhy{N}(f)$ and therefore $\vocntf{N}(g)=\vocntf{N}(\vmuhy{N}(f))=\vex{N}(f)$, where the last equality follows from Eq.~\eqref{eq:muhy-identities}.
  Now Proposition~\ref{prop:exchangeability-ex} tells us that $\vex{N}(f)\in\desirs$, because $f\in\desirs$ and $\desirs$ is exchangeable.
  Hence $\vocntf{N}(g)\in\desirs$ and therefore $g\in\cdesirs$.
  \par
  Next, we turn to the necessity part.
  Suppose that $\desirs$ is coherent and exchangeable.
  It suffices to prove that $\cdesirs\coloneqq\vmuhy{N}(\desirs)$ is a coherent set of desirable gambles on $\vcnts{N}$, and that~Eq.~\eqref{eq:representation-in-between} is satisfied for this choice of $\cdesirs$.
  \par
  We begin with the coherence of $\vmuhy{N}(\desirs)$.
  For D\ref{item:Dnonzero}, consider  $g\in\gambles(\vcnts{N})$ with $g=0$.
  Assume {\itshape ex absurdo} that $g\in\vmuhy{N}(\desirs)$, meaning that there is some $f\in\desirs$ such that $0=g=\vmuhy{N}(f)$, or in other words $f\in\weakly{\average_N}$.
  This is impossible, due to Eq.~\eqref{eq:no-intersection}.
  For D\ref{item:Dapg}, let $g>0$.
  Then obviously also $f\coloneqq\vocntf{N}(g)>0$.
  Therefore $f\in\desirs$ and, because of Eq.~\eqref{eq:muhy-identities}, $g=\vmuhy{N}(\vocntf{N}(g))=\vmuhy{N}(f)\in\vmuhy{N}(\desirs)$.
  The proof for D\ref{item:Dscl} is similar to the one for D\ref{item:Dcmb}.
  For D\ref{item:Dcmb}, let $g_1,g_2\in\vmuhy{N}(\desirs)$, so there are $f_1,f_2\in\desirs$ such that $g_1=\vmuhy{N}(f_1)$ and $g_2=\vmuhy{N}(f_2)$.
  Then by coherence of $\desirs$, $f_1+f_2\in\desirs$, and therefore, by linearity of $\vmuhy{N}$,
    \begin{equation}
      g_1+g_2
      =\vmuhy{N}(f_1)+\vmuhy{N}(f_2)
      =\vmuhy{N}(f_1+f_2)\in\vmuhy{N}(\desirs).
    \end{equation}
  Finally, we show that $\desirs=(\vmuhy{N})^{-1}(\vmuhy{N}(\desirs))$.
  Consider $f\in\desirs$, then $\vmuhy{N}(f)\in\vmuhy{N}(\desirs)$ and therefore $f\in(\vmuhy{N})^{-1}(\vmuhy{N}(\desirs))$.
  Conversely, consider a gamble~$f$ in $(\vmuhy{N})^{-1}(\vmuhy{N}(\desirs))$.
  Then $g\coloneqq\vmuhy{N}(f)\in\vmuhy{N}(\desirs)$, so we infer that there is some ${f'\in\desirs}$ such that $g=\vmuhy{N}(f)=\vmuhy{N}(f')$.
  Hence ${\vmuhy{N}(f-f')=0}$, so $f-f'\in\weakly{\average_N}$ by Eqs.~\eqref{eq:muhy-identities} and~\eqref{eq:kernel}, and therefore $f=f'+f-f'\in\desirs+\weakly{\average_N}$.
  This implies that $f\in\desirs$, since $\desirs$ is exchangeable.
\end{proof}

\begin{proof}[\bfseries Proof of Corollary~\ref{cor:finite-representation}]
  This result can be easily proved as an immediate consequence of Theorem~\ref{theo:finite-representation} and Eq.~\eqref{eq:lpr}.
  As an illustration, we give a more direct proof of the necessity part, based on Theorem~\ref{theo:exchangeability-lpr}.
  This theorem, together with Eq.~\eqref{eq:muhy-identities}, tells us that for any gamble~$f$ on~$\values^N$,
  $
  \lpr(f)
  =\lpr\bigl(\vex{N}(f)\bigr)
  =\lpr\bigl(\vocntf{N}(\vmuhy{N}(f))\bigr)
  =\clpr\bigl(\vmuhy{N}(f)\bigr).
  $
\end{proof}

\begin{proof}[\bfseries Proof of Theorem~\ref{theo:finite-representation-exnatex}]
  We begin with the second statement.
  Recall that $\exnatex\assessment=\weakly{\average_N}+\natex{\assessment}$ from Theorem~\ref{theo:exnatex}.
  Since $\vmuhy{N}$ is a linear operator, it commutes with the $\posi$ operator, and therefore:
  \begin{align*}
    \vmuhy{N}(\exnatex\assessment)
      &=\vmuhy{N}(\weakly{\average_N})+\vmuhy{N}(\natex\assessment)\\
      &=\vmuhy{N}(\natex\assessment)\\
      &=\posi\bigl(\vmuhy{N}(\gambles^+_0(\values^N)\cup\assessment)\bigr)\\
      &=\posi\bigl(\vmuhy{N}(\gambles^+_0(\values^N))\cup\vmuhy{N}(\assessment)\bigr)\\
      &=\posi\bigl(\gambles^+_0(\vcnts{N})\cup\vmuhy{N}(\assessment)\bigr)\\
      &=\natex{\vmuhy{N}(\assessment)},
  \end{align*}
  where the second equality follows from $\vmuhy{N}(\weakly{\average_N})=\{0\}$, the third from Theorem~\ref{theo:exnatex}, and the last from Theorem~\ref{theo:natex}.
  The first statement is an immediate consequence of the second and Theorems~\ref{theo:natex}, \ref{theo:exnatex} and~\ref{theo:finite-representation}.
\end{proof}

\begin{proof}[\bfseries Proof of Proposition~\ref{prop:finite-representation-update}]
  Recall that $g\in\restrict{\cdesirs}{\ocnt{m}}$ iff there is some $f\in\gambles(\values^{\rest{n}})$ such that at the same time $g=\vmuhy{\rest{n}}(f)$ and $I_{\event_{\oatom{m}}}f\in\desirs$, or in other words $\vmuhy{N}(I_{\smash[b]{\event_{\oatom{m}}}}f)\in\cdesirs$.
  We therefore consider $\cnt{M}\in\vcnts{N}$ and observe that
  \begin{equation}
    \vmuhy{N}(I_{\event_{\oatom{m}}}f\vert\cnt{M})
    =\frac{1}{\abs{\batom{M}}}\smashoperator[r]{\sum_{\sample{x}\in\batom{M}}}
    (I_{\event_{\oatom{m}}}f)(\sample{x})
    =\frac{1}{\abs{\batom{M}}}
    \smashoperator[r]{\sum_{\substack{\osample{x}\in\oatom{m},
        \rsample{x}\in\values^{\rest{n}}\\
        (\osample{x},\rsample{x})\in\batom{M}}}}
    f(\rsample{x}),
  \end{equation}
  so this value is zero unless $\cnt{M}\geq\ocnt{m}$.
  In that case we can write $\cnt{M}=\ocnt{m}+\rcnt{m}$, where $\rcnt{m}\coloneqq\cnt{M}-\ocnt{m}$ is a count vector in $\vcnts{\rest{n}}$; so we find that
  \begin{equation}
    \vmuhy{N}(I_{\event_{\oatom{m}}}f\vert\ocnt{m}+\rcnt{m})
    =\frac{1}{\abs{\atom{\ocnt{m}+\rcnt{m}}}}
    \smashoperator[r]{\sum_{\osample{x}\in\oatom{m},\rsample{x}\in\ratom{m}}}
    f(\rsample{x})
    =\frac{\abs{\oatom{m}}\,\abs{\ratom{m}}}{\abs{\atom{\ocnt{m}+\rcnt{m}}}}
    \vmuhy{\rest{n}}(f\vert\rcnt{m}).
  \end{equation}
  Hence indeed $g\in\restrict{\cdesirs}{\ocnt{m}}$ iff $+_{\ocnt{m}}(L_{\ocnt{m}}g)\in\cdesirs$.
\end{proof}

\begin{proof}[\bfseries Proof of Theorem~\ref{theo:finite-representation-polynomials}]
  It clearly suffices to give the proof in terms of count gambles.
  Because we have seen that $\vcmult{N}$ is a linear isomorphism between the linear spaces $\gambles(\vcnts{N})$ and $\vpoly[N]$, it is clear that $\cdesirs=(\vcmult{N})^{-1}(\tdesirs)$ iff $\tdesirs=\vcmult{N}(\cdesirs)$.
  \par
  Suppose that $\cdesirs$ is coherent, then we have to prove that $\tdesirs=\vcmult{N}(\cdesirs)$ is Bernstein coherent at degree $N$.
  Since $\vcmult{N}$ is a linear isomorphism, it is clear that $\tdesirs$ satisfies~B$_N$\ref{item:BNnonzero}, B$_N$\ref{item:BNscl} and~B$_N$\ref{item:BNcmb}, because $\cdesirs$ satisfies~D\ref{item:Dnonzero}, D\ref{item:Dscl} and~D\ref{item:Dcmb}.
  To show that $\tdesirs$ satisfies~B$_N$\ref{item:BNapg}, consider $p$ such that $\bexp{p}{N}>0$ and therefore $\bexp{p}{N}\in\cdesirs$ by~D\ref{item:Dapg}.
  Hence indeed $p=\vcmult{N}(\bexp{p}{N})\in\vcmult{N}(\cdesirs)=\tdesirs$.
  \par
  Suppose that $\tdesirs$ is Bernstein coherent at degree $N$, then we have to prove that $\cdesirs=(\vcmult{N})^{-1}(\tdesirs)$ is coherent.
  Since $(\vcmult{N})^{-1}$ is a linear isomorphism, it is clear that~$\cdesirs$ satisfies~D\ref{item:Dnonzero}, D\ref{item:Dscl} and~D\ref{item:Dcmb}, because $\tdesirs$ satisfies~B$_N$\ref{item:BNnonzero}, B$_N$\ref{item:BNscl} and~B$_N$\ref{item:BNcmb}.
  To show that~$\cdesirs$ satisfies~D\ref{item:Dapg}, consider $g>0$.
  Then $p=\vcmult{N}(g)$ is such that $\bexp{p}{N}=g>0$ and therefore $p\in\tdesirs$ by~B$_N$\ref{item:BNapg}.
  Hence indeed $g=(\vcmult{N})^{-1}(p)\in(\vcmult{N})^{-1}(\tdesirs)=\cdesirs$.
\end{proof}

\begin{proof}[\bfseries Proof of the equivalence of Eqs.~\eqref{eq:time-consistency} and~\eqref{eq:time-consistency-counts}]
  We begin by proving that Eq.~\eqref{eq:time-consistency} implies  Eq.~\eqref{eq:time-consistency-counts}.
  Consider any $n_1\leq n_2$.
  \begin{compactitem}
  \item[$\subseteq$] Consider any $g_2\in\venl{n_1}{n_2}(\vcdesirs{n_1})$, so there is some $g_1\in\vcdesirs{n_1}$ such that $g_2=\venl{n_1}{n_2}(g_1)$.
    Then it remains to prove that $g_2\in\vcdesirs{n_2}$.
    But  $g_1\in\vcdesirs{n_1}$ means that there is some $f_1\in\vdesirs{n_1}$ such that $g_1=\vmuhy{n_1}(f_1)$.
    It then follows from Eq.~\eqref{eq:time-consistency} that $f_2\coloneqq\vexten{n_1}{n_2}(f_1)\in\vdesirs{n_2}$, and therefore $\vmuhy{n_2}(f_2)\in\vcdesirs{n_2}$.
    But Eq.~\eqref{eq:enl-identity} tells us that
    \begin{equation}
      \vmuhy{n_2}(f_2)
      =\vmuhy{n_2}\bigl(\vexten{n_1}{n_2}(f_1)\bigr)
      =\venl{n_1}{n_2}\bigl(\vmuhy{n_1}(f_1)\bigr)
      =\venl{n_1}{n_2}(g_1)=g_2.
    \end{equation}
  \item[$\supseteq$] Consider any $g_2\in\vcdesirs{n_2}\cap\venl{n_1}{n_2}(\gambles(\vcnts{n_1}))$.
    We have to show that $g_2\in\venl{n_1}{n_2}(\vcdesirs{n_1})$.
    On the one hand, $g_2\in\vcdesirs{n_2}$ implies that there is some $f_2\in\vdesirs{n_2}$ such that $g_2=\vmuhy{n_2}(f_2)$.
    On the other hand, $g_2\in\venl{n_1}{n_2}(\gambles(\vcnts{n_1}))$ means that there is some gamble $g_1$ on $\vcnts{n_1}$ such that $g_2=\venl{n_1}{n_2}(g_1)$, and therefore also some gamble $f_1$ on~$\values^{n_1}$ such that $g_1=\vmuhy{n_1}(f_1)$ and therefore
    \begin{equation}
      g_2
      =\venl{n_1}{n_2}(g_1)
      =\venl{n_1}{n_2}\bigl(\vmuhy{n_1}(f_1)\bigr)
      =\vmuhy{n_2}\bigl(\vexten{n_1}{n_2}(f_1)\bigr),
    \end{equation}
    if we also  consider Eq.~\eqref{eq:enl-identity}.
    Hence $\vmuhy{n_2}(f_2)=\vmuhy{n_2}(\vexten{n_1}{n_2}(f_1))$, and therefore also $\vex{n_2}(f_2)=\vex{n_2}(\vexten{n_1}{n_2}(f_1))$, by Eq.~\eqref{eq:muhy-identities}.
    Since $f_2\in\vdesirs{n_2}$ we conclude from Proposition~\ref{prop:exchangeability-ex} that also $\vexten{n_1}{n_2}(f_1)\in\vdesirs{n_2}$.
    Now we invoke Eq.~\eqref{eq:time-consistency} to find that $\vexten{n_1}{n_2}(f_1)\in\vexten{n_1}{n_2}(\vdesirs{n_1})$, and therefore $f_1\in\vdesirs{n_1}$.
    But this implies that $g_1\in\vcdesirs{n_1}$ and consequently $g_2\in\venl{n_1}{n_2}(\vcdesirs{n_1})$.
  \end{compactitem}
  Next, we prove that Eq.~\eqref{eq:time-consistency-counts} implies  Eq.~\eqref{eq:time-consistency}.
  Consider any $n_1\leq n_2$.
 \begin{compactitem}
  \item[$\subseteq$] Consider any $f_2\in\vexten{n_1}{n_2}(\vdesirs{n_1})$, so there is some $f_1\in\vdesirs{n_1}$ such that $f_2=\vexten{n_1}{n_2}(f_1)$.
    Then $\vmuhy{n_1}(f_1)\in\vcdesirs{n_1}$, and therefore Eq.~\eqref{eq:enl-identity} tells us that
    \begin{equation}
      \vmuhy{n_2}(f_2)
      =\vmuhy{n_2}\bigl(\vexten{n_1}{n_2}(f_1)\bigr)
      =\venl{n_1}{n_2}\bigl(\vmuhy{n_1}(f_1)\bigr)
      \in\venl{n_1}{n_2}(\vcdesirs{n_1}).
    \end{equation}
    We then deduce from  Eq.~\eqref{eq:time-consistency-counts} that $\vmuhy{n_2}(f_2)\in\vcdesirs{n_2}$, whence indeed $f_2\in\vdesirs{n_2}$.
  \item[$\supseteq$] Consider any $f_2\in\vdesirs{n_2}\cap\vexten{n_1}{n_2}\bigl(\gambles(\values^{n_1})\bigr)$.
    Then $\vmuhy{n_2}(f_2)\in\vcdesirs{n_2}$ and there is some gamble $f_1$ on $\values^{n_1}$ such that $f_2=\vexten{n_1}{n_2}(f_1)$.
    So we deduce from  Eq.~\eqref{eq:enl-identity} that
    \begin{equation}
      \vmuhy{n_2}(f_2)
      =\vmuhy{n_2}(\vexten{n_1}{n_2}(f_1))
      =\venl{n_1}{n_2}(\vmuhy{n_1}(f_1))
      \in\venl{n_1}{n_2}(\gambles(\vcnts{n_1}))
    \end{equation}
    as well.
    Therefore  Eq.~\eqref{eq:time-consistency-counts} tells us that $\vmuhy{n_2}(f_2)\in\venl{n_1}{n_2}(\vcdesirs{n_1})$, so there is some $g_1\in\vcdesirs{n_1}$ such that $\vmuhy{n_2}(f_2)=\venl{n_1}{n_2}(g_1)$.
    Hence $\venl{n_1}{n_2}(\vmuhy{n_1}(f_1))=\venl{n_1}{n_2}(g_1)$, and we infer from Lemma~\ref{lem:enl-one-to-one} that therefore $\vmuhy{n_1}(f_1)=g_1$, whence $f_1\in\vdesirs{n_1}$.
    This implies that indeed $f_2=\vexten{n_1}{n_2}(f_1)\in\vexten{n_1}{n_2}(\vdesirs{n_1})$.
  \end{compactitem}
  This completes the proof.
\end{proof}

\begin{lemma}\label{lem:enl-one-to-one}
  Consider any $n_1\leq n_2$ in $\nats_0$.
  Then the extension map $\venl{n_1}{n_2}$ is one-to-one.
\end{lemma}

\begin{proof}
  Consider any gambles $g_1$ and $g_2$ on $\vcnts{n_1}$ and assume that $\venl{n_1}{n_2}(g_1)=\venl{n_1}{n_2}(g_2)\eqqcolon g$.
  Then we must prove that $g_1=g_2$.
  Consider the polynomial $p\coloneqq\vcmult{n_2}(g)$, then we infer from Eq.~\eqref{eq:mult-extensions} that
  \begin{equation}
    \vcmult{n_1}(g_1)
    =\vcmult{n_2}(\venl{n_1}{n_2}(g_1))
    =p
    =\vcmult{n_2}(\venl{n_1}{n_2}(g_2))
    =\vcmult{n_1}(g_2),
  \end{equation}
  which means that, with the notations of Appendix~\ref{app:bernstein}, $\bexp{p}{n_1}=g_1=g_2$ is the unique decomposition of the polynomial $p$ in terms of the Bernstein basis polynomials of degree $n_1$.
\end{proof}

\begin{proof}[\bfseries Proof of the equivalence of Eqs.~\eqref{eq:time-consistency-counts} and~\eqref{eq:time-consistency-frequencies}]
  As a first step, we prove that  Eq.~\eqref{eq:time-consistency-counts} implies Eq.~\eqref{eq:time-consistency-frequencies}.
  Consider any $n_1\leq n_2$.
  \begin{compactitem}
  \item[$\subseteq$]  Choose any $p\in\vtdesirs[n_1]$, then we know from Theorem~\ref{theo:finite-representation-polynomials} and the discussion in Appendix~\ref{app:bernstein} that there is a unique $g_1\coloneqq\bexp{p}{n_1}$ in $\vcdesirs{n_1}$ such that $p=\vcmult{n_1}(g_1)$.
    If we let $g_2\coloneqq\venl{n_1}{n_2}(g_1)$ then we infer from Eq.~\eqref{eq:mult-extensions} that $p=\vcmult{n_2}(g_2)$ as well.
    Since we infer from  Eq.~\eqref{eq:time-consistency-counts} that $g_2\in\vcdesirs{n_2}$, we see that indeed $p\in\vcmult{n_2}(\vcdesirs{n_2})=\vtdesirs[n_2]$.
  \item[$\supseteq$] Choose any $p\in\vtdesirs[n_2]\cap\vpoly[n_1]$.
    Since $p\in\vtdesirs[n_2]$ we infer from  Theorem~\ref{theo:finite-representation-polynomials} and the discussion in Appendix~\ref{app:bernstein} that there is a unique $g_2\coloneqq\bexp{p}{n_2}$ in $\vcdesirs{n_2}$ such that $p=\vcmult{n_2}(g_2)$.
    On the other hand,  since $p$ is a polynomial of degree at most~$n_1$, we know from the discussion in Appendix~\ref{app:bernstein} that there is a unique Bernstein expansion $g_1\coloneqq\bexp{p}{n_1}$ in $\gambles(\vcnts{n_1})$ such that $p=\vcmult{n_1}(g_1)$.
    The relation between the unique Bernstein expansions $g_1$ and $g_2$ is given by Zhou's formula: $g_2=\venl{n_1}{n_2}(g_1)$.
    Hence $g_2\in\venl{n_1}{n_2}\bigl(\gambles(\vcnts{n_1})\bigr)$ as well, and we infer from Eq.~\eqref{eq:time-consistency-counts} that there is some $g_3\in\vcdesirs{n_1}$ such that $g_2=\venl{n_1}{n_2}(g_3)$.
    But since we have shown before that $\venl{n_1}{n_2}$ is one-to-one [Lemma~\ref{lem:enl-one-to-one}], we infer that $g_1=g_3$ and therefore $g_1\in\vcdesirs{n_1}$, whence indeed $p\in\vcmult{n_1}(\vcdesirs{n_1})=\vtdesirs[n_1]$.
  \end{compactitem}
  Next, we prove that   Eq.~\eqref{eq:time-consistency-frequencies} implies Eq.~\eqref{eq:time-consistency-counts}.
  Consider any $n_1\leq n_2$.
  \begin{compactitem}
  \item[$\subseteq$] Choose any $g_2\in\venl{n_1}{n_2}(\vcdesirs{n_1})$.
    Then there is some $g_1\in\vcdesirs{n_1}$ such that $g_2=\venl{n_1}{n_2}(g_1)$.
    Let $p=\vcmult{n_1}(g_1)$, then we infer from Theorem~\ref{theo:finite-representation-polynomials} that $p\in\vtdesirs[n_1]$.
    But  Eq.~\eqref{eq:mult-extensions} also tells us that $p=\vcmult{n_1}(g_1)=\vcmult{n_1}\bigl(\smash[b]{\venl{n_1}{n_2}}(g_1)\bigr)=\vcmult{n_2}(g_2)$, and since also $p\in\vtdesirs[n_2]$ by Eq.~\eqref{eq:time-consistency-frequencies}, we see that indeed $g_2\in(\vcmult{n_2})^{-1}(\vtdesirs[n_2])=\vcdesirs{n_2}$.
  \item[$\supseteq$] Choose any $g_2\in\vcdesirs{n_2}\cap\venl{n_1}{n_2}(\gambles(\vcnts{n_1}))$.
    Let $p=\vcmult{n_2}(g_2)$ then it follows from Theorem~\ref{theo:finite-representation-polynomials} that $p\in\vtdesirs[n_2]$.
    But we also know that there is some $g_1\in\gambles(\vcnts{n_1})$ such that $g_2=\venl{n_1}{n_2}(g_1)$ and therefore $p=\vcmult{n_2}(\venl{n_1}{n_2}(g_1))=\vcmult{n_1}(g_1)$, by Eq.~\eqref{eq:mult-extensions}.
    So $p$ is a polynomial of degree at most $n_1$, and we then infer from   Eq.~\eqref{eq:time-consistency-frequencies} that $p\in\vtdesirs[n_1]$, whence $g_1\in(\vcmult{n_1})^{-1}(\vtdesirs[n_1])=\vcdesirs{n_1}$, and therefore indeed $g_2\in\venl{n_1}{n_2}(\vcdesirs{n_1})$.
  \end{compactitem}
  This completes the proof.
\end{proof}

\begin{proof}[\bfseries Proof that  B\ref{item:Bnonzero} is equivalent to B5 under B\ref{item:Bapg}--B\ref{item:Bcmb}]
  It is clear that B5 implies B\ref{item:Bnonzero}, because if a polynomial is zero, then so are all its Bernstein expansions.
  The proof is therefore complete if we can show that B5 follows from B\ref{item:Bnonzero}--B\ref{item:Bcmb}.
  Consider a polynomial $p$ for which there is some $n\geq\deg(p)$ such that $\bexp{p}{n}\leq0$, and assume \textit{ex absurdo} that $p\in\tdesirs$.
  Then clearly $p\neq0$ by B\ref{item:Bnonzero}, and therefore $\bexp{p}{n}<0$.
  But then  $\bexp{-p}{n}=-\bexp{p}{n}>0$, so $-p\in\tdesirs$ by B\ref{item:Bapg}, and then $0=p+(-p)\in\tdesirs$ by B\ref{item:Bcmb}, a contradiction.
\end{proof}

\begin{lemma}\label{lem:sets-of-polynomials}
  Consider a subset $\vatdesirs$ of\/ $\vpoly$, and define the sets $\vatdesirs[n]\coloneqq\vatdesirs\cap\vpoly[n]$ for all $n\in\nats$.
  Then:
  \begin{compactenum}[\upshape (i)]
  \item\label{item:sop1} $\vatdesirs[n_1]=\vatdesirs[n_2]\cap\vpoly[n_1]$ for all $0\leq n_1\leq n_2$;
  \item\label{item:sop2} For all $p\in\vpoly$, if $n\geq\deg(p)$ then $p\in\vatdesirs\asa p\in\vatdesirs[n]$;
  \item\label{item:sop3} For all $k\geq0$, $\vatdesirs=\bigcup_{n\in\nats}\vatdesirs[n]=\bigcup_{n\geq k}\vatdesirs[n]$;
  \item\label{item:sop4} $\vatdesirs$ is Bernstein coherent iff $\vatdesirs[n]$ is Bernstein coherent at degree $n$ for all $n\in\nats_0$.
  \end{compactenum}
\end{lemma}

\begin{proof}[\bfseries Proof of Lemma~\ref{lem:sets-of-polynomials}]
  The proof of the first two statements it trivial.
  \par
  We turn to the proof of~(\ref{item:sop3}).
  Since $\vatdesirs[n]\subseteq\vatdesirs$ for all $n\geq1$, we see at once that $\bigcup_{n\in\nats_0}\vatdesirs[n]\subseteq\vatdesirs$.
  To prove the converse inequality, consider any $p\in\vatdesirs$.
  With $m=\deg(p)$ we infer from~(\ref{item:sop2}) that $p\in\vatdesirs[m]$ and therefore $p\in\bigcup_{n\in\nats}\vatdesirs[n]$.
  The second equality now follows at once from~(\ref{item:sop1}).
  \par
  On to the proof of~(\ref{item:sop4}).
  \begin{compactitem}
    \item[$\Rightarrow$] Assume first of all that $\vatdesirs$ is Bernstein coherent, and consider any $n\in\nats_0$.
      Then we have to prove that $\vatdesirs[n]$ is Bernstein coherent at degree $n$.
      It is obvious that $\vatdesirs[n]$ satisfies~B$_n$\ref{item:BNnonzero}, B$_n$\ref{item:BNscl} and~B$_n$\ref{item:BNcmb} because $\vatdesirs$ satisfies~B\ref{item:Bnonzero}, B\ref{item:Bscl} and~B\ref{item:Bcmb}.
      To prove that $\vatdesirs[n]$ satisfies~B$_n$\ref{item:BNapg}, consider $p\in\vpoly[n]$ with $\bexp{p}{n}>0$.
      Since clearly $n\geq\deg(p)$, we infer from~B\ref{item:Bapg} that $p\in\vatdesirs$ and therefore indeed $p\in\vatdesirs\cap\vpoly[n]=\vatdesirs[n]$.
    \item[$\Rightarrow$] Finally, assume that  $\vatdesirs[n]$ is Bernstein coherent at degree $n$ for all $n\in\nats_0$.
      Then we have to prove that $\vatdesirs$ is Bernstein coherent.
      It follows readily from~(\ref{item:sop3}) that $\vatdesirs$ satisfies~B\ref{item:Bnonzero}, B\ref{item:Bscl} and~B\ref{item:Bcmb}.
      To prove that $\vatdesirs$ satisfies~B\ref{item:Bapg}, consider any polynomial $p$ and assume that $\bexp{p}{n}>0$ for some $n\geq\deg(p)$.
      Then clearly $p\in\vpoly[n]$ and therefore $p\in\vatdesirs[n]$, by~B$_n$\ref{item:BNapg}.
      Hence indeed $p\in\vatdesirs$. \qedhere
  \end{compactitem}
\end{proof}

\begin{lemma}\label{lem:infinite-representation}
  Consider a time-consistent, coherent and exchangeable family $\vdesirs{n}$, $n\in\nats_0$ of sets of desirable gambles on~$\values^n$, and the associated count representations  $\vcdesirs{n}\coloneqq\vmuhy{n}(\vdesirs{n})$ on $\vcnts{n}$ and frequency representations $\vtdesirs[n]\coloneqq\vcmult{n}(\vcdesirs{n})=\vmult{n}(\vdesirs{n})$ on $\vpoly[n]$.
  Let $\vtdesirs\coloneqq\bigcup_{n\in\nats_0}\vtdesirs[n]$.
  Then the sequence $\vtdesirs[n]$ is non-decreasing, and $\vtdesirs[n]=\vtdesirs\cap\vpoly[n]$.
\end{lemma}

\begin{proof}
  Because the family $\vdesirs{n}$, $n\in\nats_0$ is time-consistent, the sets $\vtdesirs[n]$ satisfy the time-consistency property~\eqref{eq:time-consistency-frequencies}.
  This already implies that the sequence $\vtdesirs[n]$ is non-decreasing.
  We now show that $\vtdesirs[n]=\vtdesirs\cap\vpoly[n]$.
  Indeed:
  \begin{multline}
    \vtdesirs\cap\vpoly[n]
    =\bigcup_{k\geq1}\vtdesirs[k]\cap\vpoly[n]
    =\Bigl(\bigcup_{1\leq k\leq n}\vtdesirs[k]\cap\vpoly[n]\Bigr)
    \cup\Bigl(\bigcup_{k>n}\vtdesirs[k]\cap\vpoly[n]\Bigr)\\
    =\Bigl(\bigcup_{1\leq k\leq n}\vtdesirs[k]\Bigr)
    \cup\Bigl(\bigcup_{k>n}\vtdesirs[n]\Bigr)
    =\vtdesirs[n]\cup\vtdesirs[n]
    =\vtdesirs[n],
  \end{multline}
  where the third and fourth equalities follow from the time-consistency condition~\eqref{eq:time-consistency-frequencies}.
\end{proof}

\begin{proof}[\bfseries Proof of Theorem~\ref{theo:infinite-representation}]
  It clearly suffices to give the proof in terms of the count representations.
  \par
  First of all, consider a Bernstein coherent $\vatdesirs\subseteq\vpoly$, then we have to prove that the $\vcdesirs{n}=(\vcmult{n})^{-1}(\vatdesirs)$, $n\in\nats_0$ are coherent and satisfy the time-consistency condition~\eqref{eq:time-consistency-counts}.
  Let $\vatdesirs[n]\coloneqq\vatdesirs\cap\vpoly[n]$ then clearly
  \begin{equation}
    \vcdesirs{n}
    =(\vcmult{n})^{-1}(\vatdesirs)
    =(\vcmult{n})^{-1}\bigl(\vatdesirs\cap\vpoly[n]\bigr)
    =(\vcmult{n})^{-1}(\vatdesirs[n]).
  \end{equation}
  We then infer from Lemma~\ref{lem:sets-of-polynomials}(\ref{item:sop4})\&(\ref{item:sop1}) that $\vatdesirs[n]$ is Bernstein coherent at degree $n$, and that the $\vatdesirs[n]$, $n\in\nats_0$ satisfy the time-consistency condition~\eqref{eq:time-consistency-frequencies}.
  Hence the  $\vcdesirs{n}=(\vcmult{n})^{-1}(\vatdesirs[n])$ satisfy the time consistency condition~\eqref{eq:time-consistency-counts}, and we infer from the Finite Representation Theorem~\ref{theo:finite-representation-polynomials} that all $\vcdesirs{n}$ are coherent.
  \par
  Conversely, suppose that we have a family of coherent $\vcdesirs{n}$ that satisfy the time-consistency condition~\eqref{eq:time-consistency-counts}.
  Let $\vtdesirs[n]=\vcmult{n}(\vcdesirs{n})$ then we know that $\vtdesirs[n]$ is Bernstein coherent at degree $n$ [by Theorem~\ref{theo:finite-representation-polynomials}] and  that the $\vtdesirs[n]$ satisfy the time-consistency condition~\eqref{eq:time-consistency-frequencies}.
  Let $\vtdesirs\coloneqq\bigcup_{n\in\nats_0}\vtdesirs[n]$.
  Then it follows from Lemma~\ref{lem:infinite-representation} that $\vtdesirs[n]=\vtdesirs\cap\vpoly[n]$, and from Lemma~\ref{lem:sets-of-polynomials}(\ref{item:sop4}) that $\vtdesirs$ is Bernstein coherent.
  Moreover, since  $\vtdesirs[n]=\vcmult{n}(\vcdesirs{n})$ and $\vcmult{n}$ is a linear isomorphism,
  \begin{equation}
    \vcdesirs{n}
    =(\vcmult{n})^{-1}(\vtdesirs[n])
    =(\vcmult{n})^{-1}\bigl(\vtdesirs\cap\vpoly[n]\bigr)
    =(\vcmult{n})^{-1}(\vtdesirs).
  \end{equation}
  To prove unicity, consider any $\vatdesirs\subseteq\vpoly$ such that  $\vcdesirs{n}=(\vcmult{n})^{-1}(\vatdesirs)$ and let $\vatdesirs[n]\coloneqq\vatdesirs\cap\vpoly[n]$.
  Then
  \begin{equation}
    \vcdesirs{n}
    =(\vcmult{n})^{-1}(\vatdesirs)
    =(\vcmult{n})^{-1}\bigl(\vatdesirs\cap\vpoly[n]\bigr)
    =(\vcmult{n})^{-1}(\vatdesirs[n])
  \end{equation}
  and therefore $\vatdesirs[n]=\vcmult{n}(\vcdesirs{n})=\vtdesirs[n]$.
  We then infer from Lemma~\ref{lem:sets-of-polynomials}(\ref{item:sop3}) that   $\vatdesirs=\bigcup_{n\in\nats_0}\vatdesirs[n]=\bigcup_{n\in\nats_0}\vtdesirs[n]=\vtdesirs$.
\end{proof}

\begin{proof}[\bfseries Proof of Theorem~\ref{theo:updated-simplex-representation}]
  We already know that the models in the updated family $\vrestrict{n}{m}$, $\rest{n}\in\nats_0$ are coherent and exchangeable, by Propositions~\ref{prop:updating-preserves-exchangeability} and~\ref{prop:sufficiency}.
  To show that this family has a frequency representation, it suffices, by the Infinite Representation Theorem~\ref{theo:infinite-representation}, to show that it is time-consistent (satisfies Eq.~\eqref{eq:time-consistency}).
  Consider any $\rest{r}\leq\rest{s}$ in $\nats_0$, then we have to show that
  \begin{equation}
    \vextenrest{r}{s}(\vrestrict{r}{m})
    =\vrestrict{s}{m}\cap\vextenrest{r}{s}\bigl(\gambles(\values^{\rest{r}})\bigr).
  \end{equation}
  \begin{compactitem}
    \item[$\subseteq$] Let $f'\in\vextenrest{r}{s}(\vrestrict{r}{m})$, so there is some $f\in\vrestrict{r}{m}$ such that $f'=\vextenrest{r}{s}(f)$.
    Now $f\in\vrestrict{r}{m}$ means that $fI_{\event_{\ocnt{m}}}\in\vdesirs{\obs{n}+\rest{r}}$.
    Since the family $\vdesirs{n}$, $n\in\nats_0$ is by assumption time-consistent, we infer that $\vextenrest{r}{s}(f)I_{\event_{\ocnt{m}}}=\vexten{\obs{n}+\rest{r}}{\obs{n}+\rest{s}}(fI_{\event_{\ocnt{m}}})\in\vdesirs{\obs{n}+\rest{s}}$, and therefore $f'=\vextenrest{r}{s}(f)\in\vrestrict{s}{m}$.
    \item[$\supseteq$] To prove the converse inequality, let $f'\in\vrestrict{s}{m}\cap\vextenrest{r}{s}\bigl(\gambles(\values^{\rest{r}})\bigr)$.
      $f'\in\vrestrict{s}{m}$ means that $f'I_{\event_{\ocnt{m}}}\in\vdesirs{\obs{n}+\rest{s}}$.
      On the other hand, $f'\in\vextenrest{r}{s}\bigl(\gambles(\values^{\rest{r}})\bigr)$ means that there is some $f\in\gambles(\values^{\rest{r}})$ such that $f'=\vextenrest{r}{s}(f)$, and therefore $\vextenrest{r}{s}(f)I_{\event_{\ocnt{m}}}\in\vdesirs{\obs{n}+\rest{s}}$.
      So we infer from the time-consistency of the family $\vdesirs{n}$, $n\in\nats_0$ that $\vextenrest{r}{s}(f)I_{\event_{\ocnt{m}}}\in\vexten{\obs{n}+\rest{r}}{\obs{n}+\rest{s}}(\vdesirs{\obs{n}+\rest{r}})$.This means that there is some $f''\in\vdesirs{\obs{n}+\rest{r}}$ such that $\vextenrest{r}{s}(f)I_{\event_{\ocnt{m}}}=\vexten{\obs{n}+\rest{r}}{\obs{n}+\rest{s}}(f'')$, which clearly implies that $f''=fI_{\event_{\ocnt{m}}}$, and therefore indeed $f\in\vrestrict{r}{m}$.
  \end{compactitem}
  \par
  The only thing that remains to be proved is Eq.~\eqref{eq:updated-simplex-representation}.
  We already know from Theorem~\ref{theo:infinite-representation} that $\vtrestrict{m}=\bigcup_{\rest{n}\in\nats_0}\vcmultrest{n}(\vcrestrict{n}{m})$.
  This triggers a concatenation of equivalences:
  \begin{align}
    p\in\vtrestrict{m}
    \asa
    &\,(\exists\rest{n}\in\nats_0)
    \bigl(\exists g\in\gambles(\vcntsrest{n})\bigr)
    (\exists g'\in\vcdesirs{\obs{n}+\rest{n}})
    \bigl(g'=+_{\ocnt{m}}(L_{\ocnt{m}}g)
    \text{ and }
    p=\vcmultrest{n}(g)\bigr)\notag\\
    \asa
    &\,(\exists\rest{n}\in\nats_0)
    (\exists g'\in\vcdesirs{\obs{n}+\rest{n}})
    \vcmult{\obs{n}+\rest{n}}(g')=\bern{\ocnt{m}}p\notag\\
    \asa
    &\,\bern{\ocnt{m}}p\in\bigcup_{\rest{n}\in\nats_0}\vcmult{\obs{n}+\rest{n}}
    (\vcdesirs{\obs{n}+\rest{n}})\notag\\
    \asa
    &\,\bern{\ocnt{m}}p\in\vtdesirs,
  \end{align}
  where the first equivalence follows from Eq.~\eqref{eq:updating-counts} and the second from Lemma~\ref{lem:updated-simplex-representation}.
  For the last equivalence, consider Lemma~\ref{lem:infinite-representation} and the fact that $\bern{\ocnt{m}}p$ is a polynomial of degree at least~$\obs{n}$.
\end{proof}

\begin{lemma}\label{lem:updated-simplex-representation}
  Consider $\obs{n},\rest{n}\in\nats_0$, and $\ocnt{m}\in\vcnts{\obs{n}}$.
  For all gambles $g$ on $\vcntsrest{n}$ and $g'$ on $\vcnts{\obs{n}+\rest{n}}$:
  \begin{equation}
    \vcmult{\obs{n}+\rest{n}}(g')=\bern{\ocnt{m}}\vcmultrest{n}(g)
    \asa
    g'=+_{\ocnt{m}}(L_{\ocnt{m}}g).
\end{equation}
\end{lemma}

\begin{proof}
  We find that
  \begin{align}
    \vcmult{\obs{n}+\rest{n}}\bigl(+_{\ocnt{m}}(L_{\ocnt{m}}g)\bigr)
    &=\smashoperator[r]{\sum_{\cnt{M}\in\vcnts{\obs{n}+\rest{n}}}}
    +_{\ocnt{m}}(L_{\ocnt{m}}g)(\cnt{M})\bern{M}\notag\\
    &=\smashoperator[r]{\sum_{\rcnt{m}\in\vcntsrest{n}}}
    L_{\ocnt{m}}(\rcnt{m})g(\rcnt{m})\bern{\ocnt{m}+\rcnt{m}}\notag\\
    &=\smashoperator[r]{\sum_{\rcnt{m}\in\vcntsrest{n}}}
    g(\rcnt{m})\bern{\ocnt{m}}\bern{\rcnt{m}}
    =\bern{\ocnt{m}}
    \smashoperator[r]{\sum_{\rcnt{m}\in\vcntsrest{n}}}
    g(\rcnt{m})\bern{\rcnt{m}}
    =\bern{\ocnt{m}}\vcmultrest{n}(g),
  \end{align}
  where the second equality follows from Eq.~\eqref{eq:count-reduction}, and the third from Eqs.~\eqref{eq:likelihood} and~\eqref{eq:bernstein}.
  The first and last equalities go back to Eq.~\eqref{eq:cmult}.
  \par
  Conversely, consider any $g'$ in $\gambles(\vcnts{\obs{n}+\rest{n}})$ such that $\vcmult{\obs{n}+\rest{n}}(g')=\bern{\ocnt{m}}\vcmultrest{n}(g)$.
  Since $\bern{\ocnt{m}}\vcmultrest{n}(g)$ is a polynomial of degree at most $\obs{n}+\rest{n}$, we know from the discussion in Appendix~\ref{app:bernstein} that there is one and only one such $g'$, as it represents the coefficients of the unique expansion of the polynomial $\bern{\ocnt{m}}\vcmultrest{n}(g)$ in the multivariate Bernstein basis of degree $\obs{n}+\rest{n}$.
  Since we have seen in the first part of the proof that $\vcmult{\obs{n}+\rest{n}}(+_{\ocnt{m}}\bigl(L_{\ocnt{m}}g)\bigr)=\bern{\ocnt{m}}\vcmultrest{n}(g)$, it follows that $g'=+_{\ocnt{m}}(L_{\ocnt{m}}g)$.
\end{proof}

\begin{proof}[\bfseries Proof of Proposition~\ref{prop:maximal-iid}]
  That $\lpr_{\vtdesirs}$ is a linear functional that dominates the $\min$ functional follows from Lemma~\ref{lem:maximal-iid-linearity}.
  We now show that $\lpr_{\vtdesirs}(p)=p(\vartheta)$ for all $p\in\vpoly$, where $\vartheta_z\coloneqq\lpr_{\vtdesirs}(\bern{e_z})$ for all $z\in\values$.
  Consider any $p\in\vpoly$ and $n\geq\deg(p)$, then we know that $p=\sum_{\cnt{m}\in\vcnts{n}}b_p^n (\cnt{m})\bern{m}$, and therefore $\lpr_{\vtdesirs}(p)=\sum_{\cnt{m}\in\vcnts{n}}b_p^n(\cnt{m})\lpr_{\vtdesirs}(\bern{m})$, using the linearity of $\lpr_{\vtdesirs}$ [Lemma~\ref{lem:maximal-iid-linearity}].
  To find out what $\lpr_{\vtdesirs}(\bern{m})$ is, observe that we can write $\bern{m}$ as a product of simpler Bernstein basis polynomials: $\bern{m}=\binom{n}{m}\prod_{z\in\values}\bern{e_z}^{m_z}$, and therefore Lemmas~\ref{lem:maximal-iid-linearity} and~\ref{lem:maximal-iid-products} tell us that $\lpr_{\vtdesirs}(\bern{m})=\binom{n}{m}\prod_{z\in\values}\lpr_{\vtdesirs}(\bern{e_z})^{m_z}=\binom{n}{m}\prod_{z\in\values}\vartheta^{m_z}=\bern{m}(\vartheta)$.
  Hence indeed $\lpr_{\vtdesirs}(p)=\sum_{\cnt{m}\in\vcnts{n}}b_p^n(\cnt{m})\bern{m}(\vartheta)=p(\vartheta)$.
  \par
  To complete the proof, consider any gamble $g$ on $\vcnts{n}$.
  Clearly,
  \begin{align*}
    \lpr_{\vcdesirs{n}}(g)
    &=\sup\set{\alpha}{g-\alpha\in\vcdesirs{n}}
    =\sup\set{\alpha}{g-\alpha\in(\vcmult{n})^{-1}(\vtdesirs)}\\
    &=\sup\set{\alpha}{\vcmult{n}(g-\alpha)\in\vtdesirs}
    =\sup\set{\alpha}{\vcmult{n}(g)-\alpha\in\vtdesirs}\\
    &=\lpr_{\vtdesirs}(\vcmult{n}(g)).
  \end{align*}
  The rest of the proof is now immediate.
\end{proof}

\begin{lemma}\label{lem:maximal-iid-linearity}
  Consider any maximal element $\vtdesirs$ of\/ $\allberndesirs(\vsimplex)$ that satisfies either of the equivalent conditions~\eqref{eq:iid-one} or~\eqref{eq:iid-two}.
  Then $\lpr_{\vtdesirs}$ is a linear functional that dominates the $\min$ functional.
\end{lemma}

\begin{proof}
  It follows from the Bernstein coherence of $\vtdesirs$ that $\lpr_{\vtdesirs}$ is super-additive [use B\ref{item:Bcmb}] and positively homogeneous [use B\ref{item:Bscl}].
  For any polynomial $p$:
  \begin{align*}
    \lpr_{\vtdesirs}(p)
    =\sup\set{\alpha}{p-\alpha\in\vtdesirs}
    &=\inf\set{\beta}{p-\beta\notin\vtdesirs}\\
    &=\inf\set{\beta}{\beta-p\in\vtdesirs}
    =\upr_{\vtdesirs}(p);
  \end{align*}
  the second equality follows from the fact that $\set{\alpha}{p-\alpha\in\vtdesirs}$ is a down-set [use~B\ref{item:Bapg} and~B\ref{item:Bcmb}], and the third equality follows from the maximality of $\vtdesirs$ and Proposition~\ref{prop:maximality}.
  This shows that $\lpr_{\vtdesirs}$ is self-conjugate, which together with the super-additivity and positive homogeneity readily implies that $\lpr_{\vtdesirs}$ is additive and homogeneous, and therefore a linear functional.
  To show that $\lpr_{\vtdesirs}$ dominates $\min$, consider any polynomial $p\in\vpoly$ and any $n\geq\deg(p)$.
  Then there are two possibilities.
  If $p$ is a constant, then $p=\min b^n_p$ and therefore $p-\alpha\in\vtdesirs[]\asa\alpha<\min b^n_p$, so  $\lpr_{\vtdesirs}(p)=\min b^n_p$.
  If $p$ is not constant, then we infer from Eq.~\eqref{eq:bernstein-coefficients-1} in Appendix~\ref{app:bernstein} that $p-\min b_p^n>0$ and therefore $p-\min b_p^n\in\vtdesirs$, by B\ref{item:Bapg}.
  Hence  $\lpr_{\vtdesirs}(p)\geq\min b^n_p$.
  So we infer that this inequality holds for all $p$ and all ${n\geq\deg(p)}$, whence indeed $\lpr_{\vtdesirs}(p)\geq\sup_{n\geq\deg(p)}\min b^n_p=\min p$, where the equality follows from Proposition~\ref{prop:ranges} in Appendix~\ref{app:bernstein}.
\end{proof}

\begin{lemma}\label{lem:maximal-iid-products}
  Consider any maximal element $\vtdesirs$ of\/ $\allberndesirs(\vsimplex)$ that satisfies either of the equivalent conditions~\eqref{eq:iid-one} or~\eqref{eq:iid-two}.
  Then $\lpr_{\vtdesirs}(\bern{m}p)=\lpr_{\vtdesirs}(\bern{m})\lpr_{\vtdesirs}(p)$ for all $p\in\vpoly$ and all count vectors $\cnt{m}$.
\end{lemma}

\begin{proof}
  Observe that for all real $\alpha$ and $\beta$:
  \begin{equation}\label{eq:maximal-iid}
    \bern{m}p-\alpha=\bern{m}(p-\beta)+(\beta\bern{m}-\alpha).
  \end{equation}
  First, consider any $\alpha<\lpr_{\vtdesirs}(\bern{m}p)$ and $\beta>\lpr_{\vtdesirs}(p)$.
  Then $\bern{m}p-\alpha\in\vtdesirs$ and $p-\beta\notin\vtdesirs$.
  If we take into account the maximality of $\vtdesirs$ and Proposition~\ref{prop:maximality}, the latter leads to $\beta-p\in\vtdesirs$, and therefore $\bern{m}(\beta-p)\in\vtdesirs$, using condition~\eqref{eq:iid-two}.
  But then Eq.~\eqref{eq:maximal-iid} and B\ref{item:Bcmb} lead to the conclusion that $\beta\bern{m}-\alpha\in\vtdesirs$.
  Hence $\lpr_{\vtdesirs}(\beta\bern{m}-\alpha)\geq0$, whence $\beta\lpr_{\vtdesirs}(\bern{m})\geq\alpha$, using the linearity of $\lpr_{\vtdesirs}$ [see Lemma~\ref{lem:maximal-iid-linearity}].
  Since this inequality holds for all $\alpha<\lpr_{\vtdesirs}(\bern{m}p)$ and $\beta>\lpr_{\vtdesirs}(p)$, we infer that $\lpr_{\vtdesirs}(p)\lpr_{\vtdesirs}(\bern{m})\geq\lpr_{\vtdesirs}(\bern{m}p)$.
  \par
  To prove the converse inequality,  consider any $\alpha>\lpr_{\vtdesirs}(\bern{m}p)$ and $\beta<\lpr_{\vtdesirs}(p)$.
  Then $\bern{m}p-\alpha\notin\vtdesirs$ and $p-\beta\in\vtdesirs$.
  If we take into account the maximality of $\vtdesirs$ and Proposition~\ref{prop:maximality}, the former leads to $\alpha-\bern{m}p\in\vtdesirs$, and the latter to $\bern{m}(p-\beta)\in\vtdesirs$, using condition~\eqref{eq:iid-two}.
  But then Eq.~\eqref{eq:maximal-iid} and B\ref{item:Bcmb} lead to the conclusion that $\alpha-\beta\bern{m}\in\vtdesirs$.
  Hence $\lpr_{\vtdesirs}(\alpha-\beta\bern{m})\geq0$, whence $\beta\lpr_{\vtdesirs}(\bern{m})\leq\alpha$, using the linearity of $\lpr_{\vtdesirs}$ [see Lemma~\ref{lem:maximal-iid-linearity}].
  Since this inequality holds for all $\alpha>\lpr_{\vtdesirs}(\bern{m}p)$ and $\beta<\lpr_{\vtdesirs}(p)$, we infer that $\lpr_{\vtdesirs}(p)\lpr_{\vtdesirs}(\bern{m})\leq\lpr_{\vtdesirs}(\bern{m}p)$.
 \end{proof}

\begin{proof}[\bfseries Proof of Theorem~\ref{theo:bernstein-natex}]
  This is an instance of Theorem~\ref{theo:natex} with linear space $\subspace\coloneqq\vpoly$ and cone ${\cone\coloneqq\vnnegpoly}$.
\end{proof}

\begin{proof}[\bfseries Proof of Proposition~\ref{prop:ranges}]
  Eq.~\eqref{eq:bern-ranges} follows from the fact that the $\bexp{p}{n}$ converge uniformly to the polynomial $p$ as $n\to\infty$; see for instance \citet{trump1996}.
  Alternatively, it can be shown \citep[see][Section~11.9]{prautzsch2002} that for $n\geq r$ and $\cnt{M}\in\vcnts{n}$:
  \begin{equation}
    \bexp{p}{n}(\cnt{M})
    =\smashoperator{\sum_{\cnt{m}\in\vcnts{r}}}\bexp{p}{r}(\cnt{m})\bern{\cnt{m}}(\frac{\cnt{M}}{n})
    +O(\frac{1}{n})
    =p(\frac{\cnt{M}}{n})+O(\frac{1}{n}).
  \end{equation}
  Hence $\min \bexp{p}{n}\geq\min p+O(\frac{1}{n})$ for any $n\geq r$, and as a consequence $\lim_{n\rightarrow \infty, n\geq r}\min\bexp{p}{n}\geq\min p$.
  If we now use
  Equation~\eqref{eq:bernstein-coefficients-4}, we see that
  $\lim_{n\to\infty,n\geq r}\min \bexp{p}{n}=\min p$.
  The proof of the other equality is analogous.
\end{proof}

\section{Multivariate Bernstein basis polynomials}\label{app:bernstein}
With any $n\geq0$ and $\cnt{m}\in\vcnts{n}$ there corresponds a Bernstein (basis) polynomial of degree~$n$ on $\vsimplex$, given by $\bern{m}(\btheta)=\abs{\batom{m}}\prod_{x\in\values}\theta_x^{m_x}$, $\btheta\in\vsimplex$.
These polynomials have a number of very interesting properties, see for instance \citet[Chapters~10 and~11]{prautzsch2002}, which we list here:
\begin{compactenum}[\upshape BP1.]
\item\label{item:bern-linindep} The set $\set{\bern{m}}{\cnt{m}\in\vcnts{n}}$ of all Bernstein basis polynomials of fixed degree $n$ is linearly independent: if $\sum_{\cnt{m}\in\vcnts{n}}\lambda_{\cnt{m}}\bern{m}=0$, then $\lambda_{\cnt{m}}=0$ for all $\cnt{m}$ in $\vcnts{n}$.
\item\label{item:bern-partition} The set $\set{\bern{m}}{\cnt{m}\in\vcnts{n}}$ of all Bernstein basis polynomials of fixed degree $n$ forms a partition of unity: $\sum_{\cnt{m}\in\vcnts{n}}\bern{m}=1$.
\item\label{item:bern-positive} All Bernstein basis polynomials are non-negative, and strictly positive in the interior of~$\vsimplex$.
\item\label{item:bern-basis} The set $\set{\bern{m}}{\cnt{m}\in\vcnts{n}}$ of all Bernstein basis polynomials of fixed degree $n$ forms a basis for the linear space of all polynomials whose degree is at most $n$.
\end{compactenum}
Property BP\ref{item:bern-basis} follows from BP\ref{item:bern-linindep} and BP\ref{item:bern-partition}.
It follows from BP\ref{item:bern-basis} that:
\begin{enumerate}[\upshape BP1.]\addtocounter{enumi}{4}
\item\label{item:bern-expansion} Any polynomial $p$ of degree $r$ has a unique expansion in terms of the Bernstein basis polynomials of fixed degree $n\geq r$,
\end{enumerate}
or in other words, there is a unique gamble $\bexp{p}{n}$ on $\vcnts{n}$ such that
\begin{equation}\label{eq:bernstein-decomposition}
  p
  =\smashoperator{\sum_{\cnt{m}\in\vcnts{n}}}\bexp{p}{n}(\cnt{m})\bern{m}
  =\vcmult{n}(\bexp{p}{n}).
\end{equation}
This tells us [also use BP\ref{item:bern-partition} and B\ref{item:bern-positive}] that each $p(\btheta)$ is a convex combination of the Bernstein coefficients $\bexp{p}{n}(\cnt{m})$, $\cnt{m}\in\vcnts{n}$ whence for all $\btheta\in\vsimplex$
\begin{equation}\label{eq:bernstein-coefficients-1}
  \min\bexp{p}{n}
  \leq\min p
  \leq p(\btheta)
  \leq\max p
  \leq\max\bexp{p}{n}.
\end{equation}
It follows from a combination of BP\ref{item:bern-partition} and BP\ref{item:bern-basis} that for all $k\geq0$ and all $\cnt{M}$ in $\vcnts{n+k}$,
\begin{equation}\label{eq:bernstein-coefficients-2}
  \bexp{p}{n+k}(\cnt{M})
  =\smashoperator{\sum_{\cnt{m}\in\vcnts{n}}}
  \frac{\abs{\batom{m}}\,\abs{\batom{M-m}}}{\abs{\batom{M}}}\bexp{p}{n}(\cnt{m}),
\end{equation}
or in other words
\begin{equation}\label{eq:zhou}
  \bexp{p}{n+k}=\venl{n}{n+k}(\bexp{p}{n}).
\end{equation}
This is \emph{Zhou's formula} \citep[see][Section~11.9]{prautzsch2002}.
Hence [let $p=1$ and use~BP\ref{item:bern-partition}] we find that for all $k\geq0$ and all $\cnt{M}$ in $\vcnts{n+k}$,
\begin{equation}\label{eq:bernstein-coefficients-3}
  \smashoperator{\sum_{\cnt{m}\in\vcnts{n}}}
  \frac{\abs{\batom{m}}\,\abs{\batom{M-m}}}{\abs{\batom{M}}}=1.
\end{equation}
The expressions~\eqref{eq:bernstein-coefficients-2} and~\eqref{eq:bernstein-coefficients-3} also imply that each $\smash[b]{\bexp{p}{n+k}}(\cnt{M})$ is a convex combination of the $\bexp{p}{n}(\cnt{m})$, and therefore $\min\bexp{p}{n+k}\geq\min\bexp{p}{n}$ and $\max\bexp{p}{n+k}\leq\max\bexp{p}{n}$.
Combined with the inequalities in~\eqref{eq:bernstein-coefficients-1}, this leads to:
\begin{equation}\label{eq:bernstein-coefficients-4}
  [\min p,\max p]
  \subseteq[\min\bexp{p}{n+k},\max\bexp{p}{n+k}]
  \subseteq[\min\bexp{p}{n},\max\bexp{p}{n}]
\end{equation}
for all $n\geq m$ and $k\geq0$.
This means that the non-decreasing sequence $\min \bexp{p}{n}$ converges to some real number not greater than $\min p$, and, similarly, the non-increasing sequence $\max\bexp{p}{n}$ converges to some real number not smaller than $\max p$.
The following proposition strengthens this.

\begin{proposition}\label{prop:ranges}
  For any polynomial $p$ on $\vsimplex$ of degree up to $r$,
  \begin{equation}\label{eq:bern-ranges}
    \lim_{\substack{n\to\infty\\n\geq r}}[\min \bexp{p}{n},\max\bexp{p}{n}]
    =[\min p,\max p]=p(\vsimplex).
  \end{equation}
\end{proposition}

\end{document}